\documentclass[11pt]{article}
\usepackage[a4paper,top=1.05in,bottom=1.05in,left=1in,right=1in]{geometry}
\usepackage{amsmath,amssymb,amsthm,mathtools,enumitem}
\allowdisplaybreaks[4]
\usepackage[T1]{fontenc}\usepackage{newtxtext,newtxmath,microtype}
\usepackage{titlesec}
\titleformat{\section}{\normalfont\Large\bfseries}{\thesection}{0.75em}{}
\titleformat{\subsection}{\normalfont\large\bfseries}{\thesubsection}{0.75em}{}
\titlespacing*{\section}{0pt}{2.6ex plus 0.8ex minus 0.2ex}{1.2ex plus 0.2ex}
\titlespacing*{\subsection}{0pt}{2.0ex plus 0.6ex minus 0.2ex}{0.9ex plus 0.2ex}
\usepackage{fancyhdr}
\fancypagestyle{plain}{\fancyhf{}\fancyfoot[C]{\thepage}}
\usepackage[dvipsnames]{xcolor}
\definecolor{linkblue}{RGB}{20,67,128}
\definecolor{citegreen}{RGB}{25,105,80}
\definecolor{urlpurple}{RGB}{100,55,125}
\usepackage[colorlinks=true,linkcolor=linkblue,citecolor=citegreen,urlcolor=urlpurple,plainpages=false,hypertexnames=false]{hyperref}
\usepackage{bookmark}
\usepackage{tikz}
\usetikzlibrary{arrows.meta,positioning,calc}
\numberwithin{equation}{section}
\newtheorem{thm}{Theorem}[section]\newtheorem{lem}[thm]{Lemma}
\newtheorem{prop}[thm]{Proposition}\newtheorem{cor}[thm]{Corollary}
\theoremstyle{remark}
\newcommand{\R}{\mathbb R}\newcommand{\C}{\mathcal C}\newcommand{\barC}{\bar{\mathcal C}}
\newcommand{\dd}{\,\mathrm d}
\DeclareMathOperator{\dist}{dist}

\newcommand{\Reg}{\operatorname{Reg}}

\newcommand{\supp}{\operatorname{supp}}
\newcommand{\1}{\mathbf 1}
\title{Uniqueness of tangent flows in mean curvature flow (II):\\ cylinders over smooth minimal cones}
\author{Hanbing Fang\qquad and \qquad Yongheng Han}
\date{\today}
\begin{document}
	\maketitle
	\begin{abstract}
		We prove uniqueness of multiplicity-one tangent flows modeled on a
		smooth minimal cone times a Euclidean factor, under strict Hardy stability, two conditions identifying the Gaussian
		kernel, and weighted stationary rigidity.
		We also obtain
		a discrete \L{}ojasiewicz inequality,
		quantitative convergence, and rigidity among smooth self-shrinkers.
		In particular, these conclusions hold for cylinders over the Simons
		cones $\C_{p,q}$ with $p,q\ge1$ and $p+q\ge7$.
	\end{abstract}
	\begingroup\small\tableofcontents\endgroup
	
	\section{Introduction}
	\subsection{Main theorem}
	This is the second paper on the uniqueness of tangent flows
	in mean curvature flow. In Part~I \cite{FH1}, we considered static
	multiplicity-one tangent flows modeled on minimal cones with an isolated
	vertex. We now study the cylindrical models $\barC:=\C\times\R^k$,
	where $\C$ is a nonflat minimal hypersurface cone with a smooth link and
	$k\ge1$. The Euclidean factor gives rise to the singular axis
	$\{0\}\times\R^k$. Our question is whether, once one tangent flow at a
	space-time point is the static multiplicity-one flow associated with
	$\barC$, every tangent flow at that point must coincide with it.
	
	The standard examples of cylindrical tangent flows are the round shrinking
	cylinders.
	Colding and Minicozzi \cite{CM} proved uniqueness of multiplicity-one
	round cylindrical tangent flows. Bamler and Lai
	\cite{BamlerLai} gave a new proof by using their PDE--ODI principle to 
	control the leading spectral modes of the rescaled flow. Cylinders over 
	regular minimal cones provide another class of cylindrical
	tangent flows. If $\C$ is a nonflat minimal hypersurface cone that is
	smooth away from its vertex, then $\barC:=\C\times\R^k$ defines a static
	self-similar flow with singular axis $\{0\}\times\R^k$.
	Such cylinders occur as tangent flows at finite-time singularities.
	Vel\'azquez \cite{Velazquez} constructed mean curvature flows whose
	parabolic rescalings at a Type~II singularity converge to a Simons cone.
	Taking the product of such a flow with $\R^k$ gives a mean curvature flow
	whose corresponding tangent flow is the Simons cone times $\R^k$.
	
	Cylindrical cones have been studied extensively in the stationary
	setting. Simon \cite{SimonCylinderStructure,SimonCylinderUniqueness}
	developed excess estimates near such cones and proved uniqueness of
	multiplicity-one tangent cones for a class of products $\C\times\R^k$
	in a class of oriented stationary hypersurfaces that includes
	area-minimizing boundaries. The transverse cone satisfies strict
	stability, strict minimization, and hypotheses on its Jacobi fields.
	For the seven-dimensional quadratic cones $\C_{3,3}$ and $\C_{2,4}$,
	the cylinders carry additional Jacobi fields. Sz\'ekelyhidi
	\cite{SzekelyhidiCylinderUniqueness} proved uniqueness
	of the tangent cone $\C_{3,3}\times\R$ for stationary integral
	varifolds with this multiplicity-one tangent cone.
	Firester, Tsiamis, and Wang \cite{FiresterTsiamisWang}
	proved uniqueness of the multiplicity-one tangent cone
	$\C_{2,4}\times\R$ for area-minimizing hypersurfaces. Sz\'ekelyhidi \cite{SzekelyhidiCylinderExamples} also
	constructed minimal hypersurfaces with an isolated singular point
	whose tangent cone has a one-dimensional axis. Mazzeo and Parker
	\cite{MazzeoParker} constructed minimal submanifolds
	in a flat product with $S^1$ whose singular sets are smooth curves
	and whose tangent cones along those curves are cylindrical.
	
	The stationary problem at infinity is relevant to the analysis near
	the axis. Edelen and Sz\'ekelyhidi \cite{EdelenSzekelyhidiLiouville}
	classified complete embedded minimal hypersurfaces lying on one side
	of a cylindrical cone whose transverse cone is strictly stable and
	strictly minimizing. If the density at infinity is less than twice
	that of the cylinder, such a hypersurface is a product of a
	Hardt--Simon leaf with the Euclidean factor. Halavati and Spolaor
	\cite{HalavatiSpolaor} constructed a complete stationary
	hypersurface that meets both sides of $\C_{3,3}\times\R$ and has this
	cylinder as its multiplicity-one blowdown. In our argument, the
	stationary limits obtained by rescaling near the axis satisfy a
	squared-distance bound with separate transverse and axial weights.
	We apply the weighted rigidity hypothesis \eqref{eq:cyl-S4} at an
	inner stationary scale. A growth estimate at infinity then identifies
	the stationary limits obtained from the flow.
	
	We prove uniqueness for smooth embedded flows with uniform area ratios
	under strict Hardy stability, two conditions identifying the
	cylindrical Gaussian kernel, and weighted stationary rigidity.
	The Euclidean factor contributes kernel fields arising from rotations
	that mix the two factors. Near the axis, the concentration argument
	selects both a transverse scale and an axial center. An estimate for
	the axial component of the normal, together with the time evolution
	of the distance, controls the Gaussian $L^2$ distance outside large
	balls after normalization. These estimates
	lead to a discrete \L{}ojasiewicz inequality, a rate of convergence
	of the rescaled flow, and rigidity among smooth self-shrinkers.
	
	\paragraph{Main result.}
	Let $m\ge3$ and $k\ge1$, and put $n:=m+k$. Let
	$\Sigma^{m-1}\subset S^m$ be a smooth, closed, connected, embedded,
	two-sided minimal hypersurface with nonflat cone
	\[
	\C:=\{r\omega:r\ge0,\ \omega\in\Sigma\}\subset\R^{m+1}.
	\]
	Write $\barC:=\C\times\R^k\subset\R^{n+1}$, $Z:=(x,y)$, $r:=|x|$,
	and $\Reg(\barC):=(\C\setminus\{0\})\times\R^k$.
	The notation $|\barC|$ denotes its unit-multiplicity integral varifold,
	and $\mu_{\barC}:=\mathcal H^n\!\llcorner\barC$
	is its area measure. Every reference cylinder has multiplicity one. Choose a unit normal
	$\mathbf n_\Sigma$ on $\Sigma$ and denote the induced normal on $\barC$ by $\mathbf n_{\barC}$.
	For any hypersurface $N$, write its local unit normal as
	$\mathbf n_N:=(\mathbf n_x,\mathbf n_y)$ under the splitting
	$\R^{n+1}=\R^{m+1}_x\times\R^k_y$. We write $\mathbf n$ when $N$ is clear.
	We write $A_\Sigma$ for the \textbf{second fundamental form} and $H_\Sigma$
	for \textbf{scalar mean curvature}. With
	$\Delta_\Sigma:=\operatorname{div}_\Sigma\nabla_\Sigma$ nonpositive, put
	\[
	L_\Sigma:=\Delta_\Sigma+|A_\Sigma|^2+(m-1),\qquad
	K:=\ker L_\Sigma,\qquad
	\mathcal A_\Sigma:=-\Delta_\Sigma-|A_\Sigma|^2.
	\]
	The \textbf{orthogonal complement} $K^\perp$ and the projections are taken
	in the fixed $L^2(\Sigma,\dd\mu_\Sigma)$ inner product, as in
	\cite{FH1}. Put $d_K:=\dim K$. For an orthonormal basis
	$\{\psi_l\}_{l=1}^{d_K}$ of $K$, define
	\begin{equation}
		P_Kv:=\sum_{l=1}^{d_K}\left(\int_\Sigma v\psi_l\dd\mu_\Sigma\right)\psi_l,
		\qquad P_{K^\perp}v:=v-P_Kv.
	\end{equation}
	Let $(\mu_j,\phi_j)$ be an $L^2(\Sigma)$ orthonormal eigensystem of
	$\mathcal A_\Sigma$, with $\mu_0$ the lowest eigenvalue and
	$\phi_0>0$. We write
	$T_\C:=\{\langle e,\mathbf n_\Sigma\rangle:e\in\R^{m+1}\}$.
	We impose the following conditions on $\barC$.
	
	\begin{itemize}
		\item \textbf{Strict Hardy stability.}
		\begin{equation}\tag{S1}\label{eq:cyl-S1}
			\mu_0>-(m-2)^2/4 .
		\end{equation}
		Under this condition set
		$\nu_j:=((m-2)^2/4+\mu_j)^{1/2}$ and
		$\gamma_j^\pm:=-(m-2)/2\pm\nu_j$.
		
		\item \textbf{The zero-homogeneity space.}
		\begin{equation}\tag{S2}\label{eq:cyl-S2}
			\ker\mathcal A_\Sigma=T_\C.
		\end{equation}
		
		\item \textbf{Cylindrical Gaussian nonresonance.}
		\begin{equation}\tag{S3}\label{eq:cyl-S3}
			\gamma_j^+\notin\{-1,-2,\ldots\}\qquad\text{for every }j .
		\end{equation}
		
		\item \textbf{Weighted stationary rigidity.}
		For some $B>k/2$, put
		\[
		\mathcal E_{B}(V;\barC):=\int_{\R^{n+1}}
		\frac{\operatorname{dist}(Z,\barC)^2}
		{(1+|x|^2)(1+|y|^2)^B}\,\dd\|V\|(Z),
		\]
		where for a varifold $V$, we write $\rVert V\rVert$ for its weight measure.
		For every nonzero stationary integral $n$-varifold $V$ with
		finite all-center area ratios, we require
		\begin{equation}\tag{S4}\label{eq:cyl-S4}
			\left.
			\begin{gathered}
				(D_{1/R})_\# V \longrightarrow |\barC|
				\quad \text{as } R\to\infty,\\
				\mathcal E_{B}(V;\barC)<\infty
			\end{gathered}
			\right\}
			\quad\Longrightarrow\quad
			V=|\barC|.
		\end{equation}
		Here $D_{1/R}(Z):=Z/R$ and the convergence is full varifold
		convergence with multiplicity one.
	\end{itemize}
	
	The two papers use the same angular operator
	$\mathcal A_\Sigma=-L_\Sigma+(m-1)=-\Delta_\Sigma-|A_\Sigma|^2$.
	In Part~I \cite[Section~1.1]{FH1}, condition \textup{(S2)} is radial
	Gaussian nonresonance,
	which takes the form
	\[
	(1-2l)(m-1-2l)\notin\operatorname{Spec}\mathcal A_\Sigma
	\qquad\left(l\in\mathbb Z,\quad 1\le l<m/4\right).
	\]
	Under \eqref{eq:cyl-S1}, we have $\gamma_j^+>-(m-2)/2$.
	The identity
	$\mu_j=\gamma_j^+(\gamma_j^++m-2)$ therefore shows that
	\textup{(S2)} of \cite{FH1} is equivalent to
	$\gamma_j^+\notin\{-1,-3,-5,\ldots\}$ for every $j$. Thus \eqref{eq:cyl-S3} contains
	the nonresonance condition of Part~I and also excludes the negative
	even integers.
	
	The additional axial degrees explain the role of
	\eqref{eq:cyl-S2}. By \eqref{eq:cyl-gaussian-spectrum}, a Gaussian
	mode with radial degree $l$ and axial Hermite degree $|\alpha|$
	has eigenvalue zero exactly when
	\[
	\gamma_j^++2l+|\alpha|=1.
	\]
	For the isolated cone, the axial degree is zero, so the possible
	indicial roots are $1,-1,-3,\ldots$. For the cylinder, the choice
	$l=0$ and $|\alpha|=1$ also allows $\gamma_j^+=0$.
	It gives the fields $y_i\phi(\omega)$ with
	$\phi\in\ker\mathcal A_\Sigma$. By differentiating the minimal cone
	equation under ambient translations, we obtain
	$T_\C\subset\ker\mathcal A_\Sigma$. For
	$\phi(\omega)=\langle e,\mathbf n_\Sigma(\omega)\rangle$,
	the field $y_i\phi(\omega)$ is the normal velocity of the rotation
	mixing the transverse direction $e$ with the $i$th axial direction.
	Condition \eqref{eq:cyl-S2} requires that these mixed rotations
	exhaust the fields arising from the zero-homogeneity space.
	Together with \eqref{eq:cyl-S3}, it gives the Gaussian kernel
	$rK\oplus\bigoplus_{i=1}^k y_iT_\C$ identified in
	Proposition~\ref{prop:gaussian-bordered}.
	
	Conditions \eqref{eq:cyl-S1}--\eqref{eq:cyl-S3} identify the
	Gaussian Jacobi fields satisfying the relative tube bound
	in Proposition~\ref{prop:gaussian-bordered}.
	Condition \eqref{eq:cyl-S4} excludes nontrivial stationary limits
	arising from concentration at the axis. The labels \textup{(S1)--(S4)}
	refer throughout to the cylindrical conditions above.
	
	\begin{thm}[Uniqueness of tangent flows]\label{thm:main}
		Let $\Sigma$, $\C$, and $\barC=\C\times\R^k$ be as above. Suppose
		that $\barC$ satisfies \eqref{eq:cyl-S1}--\eqref{eq:cyl-S4}. Let $M_t$, $t<T$, be a smooth
		embedded boundaryless mean curvature flow, spatially and spacetime
		proper, with
		\begin{equation}\label{equ:areabd}
			\mathcal H^n(M_t\cap B_R(Z_0))\le\Lambda R^n
			\qquad(t<T,\ Z_0\in\R^{n+1},\ R>0).
		\end{equation}
		If one tangent flow at $(X_0,T)$ is the static multiplicity-one
		$\barC$, then every tangent flow there is the same cylinder in the
		original coordinates.
	\end{thm}

	We verify these hypotheses for quadratic cone cylinders in
	Section~\ref{sec:examples}. The restriction on the transverse dimension
	comes from cylindrical nonresonance: when $m=7$ the principal larger
	indicial root is $-2$, which creates additional Gaussian kernel fields.
	
	\begin{cor}\label{cor:examples}
		Let $p,q\ge1$ and $m:=p+q+1\ge8$, and put
		\[
		\C_{p,q}:=\{(x_1,x_2)\in\R^{p+1}\oplus\R^{q+1}:
		q|x_1|^2=p|x_2|^2\}.
		\]
		For every $k\ge1$, the conclusion of Theorem~\ref{thm:main} holds
		for $\barC:=\C_{p,q}\times\R^k$. The Gaussian area excess decays
		exponentially in rescaled time.
	\end{cor}
	
	As in \cite[Theorem~1.3]{FH1}, we also obtain quantitative consequences
	of the estimates used to prove uniqueness. At the given base point $(X_0,T)$, define the
	\textbf{rescaled mean curvature flow} and the \textbf{shrinker quantity} $\Phi$ by
	\begin{equation}\label{eq:flow}
		N_s:=e^{s/2}(M_{T-e^{-s}}-X_0),\qquad
		\rho:=(4\pi)^{-n/2}e^{-|Z|^2/4},\qquad
		\Phi:=\mathbf H+Z^\perp/2.
	\end{equation}
	We use $c_n:=(4\pi)^{-n/2}$ for the Gaussian normalization. Define the \textbf{Gaussian area}, the \textbf{Gaussian shrinker norm}, and the \textbf{Gaussian area excess}, respectively, by
	\begin{equation}
		F(N):=\int_N\rho\dd\mu_{N},\qquad
		G(N)^2:=\int_N|\Phi|^2\rho\dd\mu_{N},\qquad E(s):=F(N_s)-F(\barC).
	\end{equation}
	Huisken's monotonicity formula and the prescribed tangent flow give
	\begin{equation}\label{eq:energy-dissipation}
		E(s)\ge0,\qquad \lim_{s\to\infty}E(s)=0,\qquad
		E(s_1)-E(s_2)=\int_{s_1}^{s_2}G(N_\tau)^2\dd\tau
		\quad(s_1\le s_2).
	\end{equation}
	The \textbf{reference cylinders} $\barC_z:=R_{\mathsf B}(\C_a\times\R^k)$, with $z:=(a,\mathsf B)$,
	are constructed from the analytic reduced links of
	\cite[Proposition~2.5]{FH1}. Here $a\in K$ is the link parameter and
	$\mathsf B:\R^k\to\R^{m+1}$ is a small linear map determining a rotation
	between the transverse and axial factors: $R_{\mathsf B}:=\exp(S_{\mathsf B})$, where
	$S_{\mathsf B}(x,y):=(\mathsf B y,-\mathsf B^\ast x)$; see \eqref{eq:augmented-family}.
	The parameter $z=0$ gives
	$\barC$, and $z=(a,0)$ gives the unrotated cylinder $\barC_a$.
	The reference cylinders need not be minimal. We write
	\begin{equation}\label{eq:reference-distance}
		j_z(Z):=\dist(Z,\barC_z),\qquad
		D_z(N)^2:=\int_N j_z(Z)^2\rho\dd\mu_N.
	\end{equation}
	On each interval $[s,s+2]$ we choose one constant parameter by
	minimizing the mean squared positional distance, together with the
	squared link residual, on $[s+\tfrac12,s+\tfrac32]$; see
	Proposition~\ref{prop:joint-full}.
	
	The next theorem bounds the Gaussian area excess by its decrease on
	every sufficiently late interval of length two. It also controls the
	remaining Gaussian length. For related discrete \L{}ojasiewicz estimates,
	see \cite[Theorem~6.1]{CM} and \cite[Theorem~1.3]{FL}.
	
	\begin{thm}[Discrete \L{}ojasiewicz inequality]\label{thm:energy-loja}
		Under the hypotheses of Theorem~\ref{thm:main}, there are
		$s_0$, $C<\infty$, and $0<\theta\le1/2$ such that for every
		$s\ge s_0$,
		\begin{equation}\label{eq:cyl-loja}
			E(s)^{2(1-\theta)}\le C(E(s)-E(s+2)),\qquad
			\int_s^\infty G(N_\tau)\dd\tau\le CE(s)^\theta.
		\end{equation}
		If $\theta<1/2$, then
		\[
		E(s)\le C(1+s-s_0)^{-1/(1-2\theta)}.
		\]
		For $\theta=1/2$ the decay is exponential.
		The area measures $\mu_{N_s}$ converge locally weakly as Radon
		measures to $\mu_{\barC}$ as $s\to\infty$, and the rescaled flow
		converges smoothly on compact subsets of $\operatorname{Reg}(\barC)$.
	\end{thm}
	
	The first estimate in \eqref{eq:cyl-loja} is initially local: it applies
	on intervals where the fixed-reference distance $D_0$ and the
	nonnegative energy $E$ remain below fixed small thresholds.
	
	The same local estimates apply to a constant rescaled flow, for which
	the shrinker quantity vanishes. They give the following rigidity result.
	
	\begin{cor}\label{cor:shrinker-rigidity}
		If the nonflat cylinder $\barC$ satisfies the cone hypotheses of Theorem~\ref{thm:main}, there
		is no sequence of smooth, properly embedded, boundaryless normalized
		self-shrinkers converging to $|\barC|$ locally as varifolds.
	\end{cor}
	
	\subsection{Strategy and organization}
	The proof follows the analytic reduction and normalized compactness
	scheme of \cite{FH1}.
	
	\paragraph{(A) Analytic reduction and fitting the reference family.}
	We solve the complementary spherical equation and retain its
	finite-dimensional residual. By \cite[Proposition~2.5]{FH1}, we have
	\[
	k(a):=a+\eta(a),\qquad \eta(a)\in K^\perp,\qquad
	P_{K^\perp}\mathcal G(k(a))=0,
	\]
	and, with $f(a):=\mathcal A(k(a))$,
	\[
	q(a):=\nabla f(a)=\mathcal G(k(a))\in K,\qquad
	|f(a)-f(0)|^{1-\theta}\le C|q(a)|,
	\quad 0<\theta\le1/2.
	\]
	The reference mean curvature $-q(a)/(r\mathfrak w_a)$ in
	\eqref{eq:reduced-cone-curvature} enters the normalized graph equation.
	
	Proposition~\ref{prop:gaussian-bordered} identifies the Gaussian
	Jacobi fields satisfying the relative tube bound. They span
	\[
	\mathcal G_{\barC}=rK\oplus\bigoplus_{\ell=1}^k y_\ell T_\C.
	\]
	The second summand consists of mixed rotations of the transverse and
	axial factors. Thus we augment the reduced family by $R_{\mathsf B}$ and fit
	$\barC_z:=R_{\mathsf B}\barC_a$. On an interval $J:=[t,t+2]$, with $I:=[t+\tfrac12,t+\tfrac32]$, the
	fitting quantity is
	\[
	\mathcal J_I(z):=\frac1{|I|}\int_I\int_{N_s}W_z\rho\dd\mu_{N_s}\dd s
	+|q(a)|^2.
	\]
	Here $W_z$ is the squared positional function of Lemma~\ref{lem:distance-function};
	it equals squared distance in a regular normal tube and is uniformly
	comparable to it everywhere. The minimizing parameter is
	constant on the interval. Its first-variation equation supplies
	orthogonality of the normalized limit to $\mathcal G_{\barC}$.
	
	\paragraph{(B) Normalized compactness.}
	The estimate to be proved is
	\[
	\min_z\mathcal J_I(z)\le C\int_JG(N_s)^2\dd s;
	\]
	see Proposition~\ref{prop:joint-full}. If it fails, choose minimizing
	parameters $z_i\to0$ and set
	\[
	\eta_i^2:=\mathcal J_I(z_i)>0,\qquad
	\Delta_i:=\int_JG(N_s^i)^2\dd s=o(\eta_i^2).
	\]
	The reference residual satisfies $|q(a_i)|\le\eta_i$.
	On compact regular subsets, Brakke's local regularity theorem
	\cite[Section~6.11]{Brakke} gives smooth multiplicity-one graphical
	convergence on interior time intervals. The normalized graph equation in the proof
	of Proposition~\ref{prop:joint-full}, Step~2, gives, after passage to
	a subsequence, the limits in \eqref{eq:cyl-forced-limit}:
	\[
	q(a_i)/\eta_i\longrightarrow\xi\in K,\qquad
	u_i/\eta_i\longrightarrow U,\qquad L_{\barC}U=r^{-1}\xi.
	\]
	The normalized velocity tends to zero, so $U$ is independent of time.
	Convergence on regular compact sets alone does not preserve the
	normalization: squared distance might concentrate on the axis or
	escape to infinity. The next two steps exclude these possibilities.
	
	\paragraph{(C) Weighted rigidity and concentration at the axis.}
	We rule out concentration of the normalized squared-distance integral at the singular axis.
	Align the fitted axis with $\{x=0\}$ and write
	$T_{\delta,L}=\{r<\delta,\ |y|<L\}$. Our aim is the flow tube estimate \eqref{eq:cyl-relative-tube}
	of Lemma~\ref{lem:cyl-tube}:
	for $2<\sigma_+<\min\{4,2+2\nu_0\}$, each fixed $L$, and all sufficiently large $i$,
	\begin{equation}\label{eq:strategy-flow-tube}
		\frac{1}{\eta_i^2}\int_I\int_{N_s^i\cap T_{\delta,L}}
		W_{z_i}\rho\dd\mu_{N_s^i}\dd s
		\le C_L\delta^{\sigma_+}\qquad(0<\delta\le\delta_0).
	\end{equation}
	We prove \eqref{eq:strategy-flow-tube} by contradiction. If it fails,
	there are a fixed $L$ and radii $\delta_i\downarrow0$ such that
	\begin{equation}\tag{\ref*{eq:strategy-flow-tube}a}\label{eq:strategy-failure}
		\frac{1}{\delta_i^{\sigma_+}\eta_i^2}
		\int_I\int_{N_s^i\cap T_{\delta_i,L}}W_{z_i}\rho\dd\mu_{N_s^i}\dd s
		\longrightarrow\infty.
	\end{equation}
	Applying scale and center selection to \eqref{eq:strategy-failure} gives rescaled flows whose normalized
	squared-distance measures have unit mass in a fixed cylinder. When the rescaled squared-distance integral has a positive limit,
	we exclude it using the following stationary rigidity argument.
	For stationary varifolds, choose $\max\{0,2-2\nu_0\}<\sigma_-<2$. The rectangular
	estimate \eqref{htps:rectangular-bound} in
	Corollary~\ref{cor:rectangular-selection} controls the squared-distance measure at the selected center and scale.
	Let $a_i$ be the selected transverse radius and $P_i$ the integral of
	$\dist(Z,\barC)^2$ over the selected cylinder with respect to $\|V_i\|$.
	If $P_i/a_i^{n+2}\to\kappa^2>0$, the rectangular bound implies that the limiting
	stationary varifold has finite moment $\mathcal E_B$. Its unit blowdown and
	\eqref{eq:cyl-S4} force it to equal $|\barC|$, contradicting its positive height. If
	$P_i/a_i^{n+2}\to0$, the normalized graph gives an ordinary Jacobi
	field that is nonzero by the mass retained away from the axis. Its tube estimate permits testing across the axis by
	\eqref{eq:cyl-S1}. The growth bound at infinity forces every Gaussian
	moment to vanish, so the field is zero. This contradiction yields the stationary tube estimate \eqref{eq:endpoint-static-tube}
	of Lemma~\ref{lem:endpoint-static-tube}. Write $d_{\barC}(Z)=\dist(Z,\barC)$.
	Suppose $V_i$ are stationary integral varifolds with uniform all-center area ratios,
	$F_{z,h}(V_i)\le F(\barC)$ for all $z$ and $h>0$, and $\eta_i\downarrow0$. Then
	\begin{equation}\label{eq:strategy-stationary-tube}
		\int_{B_R}d_{\barC}^2\dd\|V_i\|\le C_R\eta_i^2
		\quad\text{for every fixed }R\ge1
		\quad\Longrightarrow\quad
		\int_{T_{\delta,L}}d_{\barC}^2\dd\|V_i\|
		\le C_{L,\sigma_-}\delta^{\sigma_-}\eta_i^2
	\end{equation}
	for fixed $L\ge1$, all sufficiently large $i$, and $0<\delta\le1$.
	
	Next let $V$ satisfy the hypotheses of Proposition~\ref{prop:endpoint-static-growth}, and put
	$H(R):=\int_{B_R}\dist(Z,\barC)^2\dd\|V\|$. Suppose $H$ is nonzero. If
	$H(R)=O(R^d)$ with $d<\min\{n+2,k+2+2\nu_0\}$, normalize at
	a sequence of radii going to infinity for which the larger-ball
	height is controlled by $H(R)$. After division by $\eta_i^2$, the dilations have unit squared-distance
	mass in $B_1$. Applying \eqref{eq:strategy-stationary-tube} with $L=1$,
	we choose $\delta$ so that $C_{1,\sigma_-}\delta^{\sigma_-}<1/2$.
	At least half the normalized mass then lies in $B_1\cap\{r\ge\delta\}$. Hence a nonzero
	height gives a nonzero ordinary Jacobi field. Its growth satisfies
	\eqref{eq:endpoint-static-field-growth}, so
	Lemma~\ref{lem:endpoint-static-jacobi} forces it to vanish,
	contradicting the retained mass. Thus Proposition~\ref{prop:endpoint-static-growth}
	gives the rigidity implication
	\begin{equation}\label{eq:strategy-rigidity}
		H(R)\le CR^d\quad(R\ge1),\qquad d<\min\{n+2,k+2+2\nu_0\}
		\quad\Longrightarrow\quad V=|\barC|.
	\end{equation}

	We return to the flows selected under the assumption that
	\eqref{eq:strategy-flow-tube} fails.
	In the notation of \eqref{eq:cyl-rescaled-normalization}, the measures
	$\widehat\mu_i$ satisfy $\widehat\mu_i(T_{2,1})=1$, and the squared-distance integral before normalization is
	$\kappa_i^2=\mathsf P_i/\ell_i^{n+2}$ on the rescaled cylinder.
	Since $0<\kappa_i\le C$, we can pass to a subsequence with $\kappa_i\to\kappa\ge0$.
	
	If $\kappa>0$, the unit mass gives a stationary limit $V$ with
	$H(R)>0$ for some fixed $R$. Yet $V$ has full unit blowdown, and
	\eqref{eq:cyl-micro-selection} gives $H(R)\le CR^{k+\sigma_++\varepsilon}$
	for $R\ge1$. For sufficiently small $\varepsilon>0$,
	the rigidity implication \eqref{eq:strategy-rigidity} gives $V=|\barC|$,
	so $H(R)=0$ for every $R$. This contradicts the positive height.
	
	If $\kappa=0$, the first bound in \eqref{eq:cyl-micro-selection}
	leaves at least half the unit mass away from the axis. The normalized
	graphs therefore yield a nonzero ordinary Jacobi field $U$. Its tube
	and growth bounds \eqref{eq:cyl-micro-limit-growth} imply $U=0$ by
	Lemma~\ref{lem:endpoint-static-jacobi}, contradicting this retained mass.
	Both cases are impossible. Thus \eqref{eq:strategy-flow-tube} holds: the normalized tube
	contribution tends to zero uniformly in $i$ as $\delta\downarrow0$.
	
	\paragraph{(D) Axial tilt and the Gaussian tail.}
	We rule out loss at spatial infinity by bounding the Gaussian second
	moment of the squared-distance measure by $C\eta_i^2$. To obtain the
	moment bound, we first use the tube decay from (C) to control axial tilt.
	The inequality $\sigma=\sigma_+>2$ is needed to sum the transverse Hardy
	moment. Summing \eqref{eq:strategy-flow-tube} over dyadic shells gives
	\begin{equation}\label{eq:strategy-hardy}
		\int_I\int_{N_s^i\cap\{r<\delta_0,\ |y|<L\}}
		r^{-2}j_{z_i}^2\rho\dd\mu_{N_s^i}\dd s
		\le C\eta_i^2\sum_{j\ge0}(2^{-j}\delta_0)^{\sigma-2} \le C\delta_0^{\sigma-2}\eta_i^2.
	\end{equation}
	Lemma~\ref{lem:axial-tilt-estimate} uses \eqref{eq:strategy-hardy} on a compact axial
	region and a product radial field to obtain
	\begin{equation}\label{eq:strategy-tilt}
		\int_I\int_{N_s^i}|\mathbf n_y|^2\rho\dd\mu_{N_s^i}\dd s\le C\eta_i^2.
	\end{equation}
	The product field has zero weighted divergence on each fitted reference cylinder.
	
	The evolution estimate bounds
	$\int_I\int_{N_s^i}W_{z_i}|\Phi_i|^2\rho\dd\mu_{N_s^i}\dd s$
	by $C\eta_i^2$, as in \eqref{eq:cyl-distance-function-weighted-dissipation}. In the Gaussian first variation with the field $W_{z_i}Z$,
	the axial error contains $y\cdot\mathbf n_y$. Since $|DW_{z_i}|^2\le CW_{z_i}$, we get
	\begin{equation}\label{eq:strategy-axial-error}
		|(DW_{z_i}\cdot\mathbf n)(y\cdot\mathbf n_y)|
		\le \varepsilon |Z|^2 W_{z_i}+C_\varepsilon|\mathbf n_y|^2.
	\end{equation}
	In \eqref{eq:strategy-axial-error} we absorb the first term into the
	second moment and use \eqref{eq:strategy-tilt} to bound the integral
	of the second. The weighted dissipation bound controls the
	velocity term by the same absorption, and \eqref{eq:cyl-transverse-moment-bound}
	controls the transverse error. Integrating in time as in Lemma~\ref{lem:cyl-gaussian-tightness} gives
	\begin{equation}\label{eq:strategy-tail}
		\int_I\int_{N_s^i\setminus B_R}W_{z_i}\rho\dd\mu_{N_s^i}\dd s
		\le R^{-2}\int_I\int_{N_s^i}|Z|^2W_{z_i}\rho\dd\mu_{N_s^i}\dd s
		\le CR^{-2}\eta_i^2.
	\end{equation}
	
	\paragraph{(E) Gaussian area, finite length, and uniqueness.}
	With concentration excluded both at the axis and at infinity, the
	convergence in (B) preserves the full normalization:
	$\|U\|_{L^2_\rho}^2+|\xi|^2=1$. The tube bound makes the cutoff error
	at the axis tend to zero in Lemma~\ref{lem:all-root-moment-test}.
	By testing the punctured equation against $r\xi$, we obtain $\xi=0$.
	Proposition~\ref{prop:gaussian-bordered} then gives $U\in\mathcal G_{\barC}$,
	while the fitting equation gives $U\perp\mathcal G_{\barC}$.
	This gives $U=0$, contradicting the normalization and proving the interval estimate in Proposition~\ref{prop:joint-full}.
	
	Proposition~\ref{prop:anchorarea} gives a one-sided comparison with the
	fitted reference. Together with the interval estimate, it provides a
	middle time $s_{\rm mid}\in[s+\tfrac12,s+\tfrac32]$ such that, writing
	$\Delta:=E(s)-E(s+2)$, we have
	\begin{align*}
		D_z(N_{s_{\rm mid}})^2+G(N_{s_{\rm mid}})^2+|q(a)|^2\le C\Delta,
	\end{align*}
	and
	\begin{align*}
		E(s)\le E(s_{\rm mid})+\Delta
		\le |F(\barC_z)-F(\barC)|+C\Delta.
	\end{align*}
	The reduced gradient inequality thus implies
	\[
	E(s)^{1-\theta}
	\le C\bigl(|q(a)|+\Delta^{1-\theta}\bigr)
	\le C\Delta^{1/2},
	\]
	which is the local discrete \L{}ojasiewicz inequality.
	
	For $E_j:=E(s_*+2j)$ and $\Delta_j:=E_j-E_{j+1}$,
	concavity gives
	\[
	\sqrt{\Delta_j}\le C(E_j^\theta-E_{j+1}^\theta),\qquad
	\int_{s_*}^{s_*+2N_{\rm blk}}G(N_s)\dd s
	\le\sqrt2\sum_{j=0}^{N_{\rm blk}-1}\sqrt{\Delta_j}\le CE(s_*)^\theta.
	\]
	Here $N_{\rm blk}$ is the number of complete two-unit intervals
	on which we apply the local inequality. The squared-distance evolution estimate
	\eqref{eq:cyl-fixed-distance-transport} prevents $D_0$ from reaching
	a fixed small threshold provided the Gaussian distance $D_0(N_{s_*})$
	to $\barC$ and the Gaussian area excess $E(s_*)$ are sufficiently small;
	see \eqref{eq:cyl-continuation-entry} and the continuation argument
	in \cite[Proposition~4.5]{FH1}. Finite length identifies the limit with
	the prescribed cylinder in the original coordinates, including its axis.
	
	\paragraph{Organization.}
	Section~\ref{sec:preliminaries} develops the link, spectral, and
	Gaussian area estimates. Section~\ref{sec:compactness}
	proves the measure selection, transverse nonconcentration, axial tilt,
	and relative Gaussian tail estimates. Section~\ref{sec:main-proofs}
	proves the main conclusions: Subsection~\ref{sec:discrete} establishes
	the fitted interval estimate, Subsection~\ref{sec:uniqueness} proves
	continuation and quantitative convergence, and the final subsection
	treats smooth self-shrinkers. Section~\ref{sec:examples} verifies the
	hypotheses for quadratic cone cylinders.
	
	\paragraph{Acknowledgement.} Hanbing Fang is supported by Hong Kong RGC grants No.~14300623 and No.~14304225. Yongheng Han is supported by the Project of Stable Support for Youth Team in Basic Research Field, Chinese Academy of Sciences (YSBR-001).

	\paragraph{AI disclosure.} This paper was developed with AI assistance.
	The principal AI tools were ChatGPT-6 Astra and the Danus
	system. The authors take full responsibility for the content of this
	article, including its mathematical claims, arguments, and conclusions.
	
	\section{Preliminary estimates}\label{sec:preliminaries}
	Throughout the estimates below, all smooth flows are embedded and
	boundaryless, are spatially and spacetime proper, and satisfy the
	uniform all-center area-ratio bound \eqref{equ:areabd}. Smooth
	hypersurfaces used in the spatial estimates are likewise properly
	embedded and boundaryless and satisfy this bound. For sequences,
	the constant is uniform in the sequence index.
	
	Throughout, $B_R(Z_0)$ is the ambient Euclidean ball, $B_R:=B_R(0)$,
	and $\dd\mu_N$ is the area measure on $N$. We omit subscripts when the
	measure is clear from the integral. The letters $c,C$ denote positive
	constants allowed to change from line to line. We use
	\textbf{big-$O$ notation} $f=O_\Lambda(g)$, for $g\ge0$, to mean
	$|f|\le C_\Lambda g$ uniformly in the variables under consideration;
	the subscript records additional dependence on $\Lambda$.
	The notation $f_i=o(g_i)$, for $g_i>0$, means $f_i/g_i\to0$.
	Unless an additional dependence is indicated, constants depend on
	the fixed cylinder, the area-ratio bound, and the fixed regions
	and cutoffs in the relevant estimate. For a hypersurface $N$, we use the Gaussian pairing
	$\langle u,v\rangle_\rho:=\int_Nuv\rho\dd\mu_N$ and its associated
	norm $\|u\|_{L^2_\rho(N)}$; the domain is omitted only when it is clear.
	For $a,q(a),\xi\in K$, $|\cdot|$ denotes the fixed $L^2(\Sigma)$ norm.
	For a varifold $V$, $F(V):=\int_{\R^{n+1}}\rho\dd\|V\|$.
	All compact test functions are ambient functions unless a domain is specified. We distinguish
	the ambient radius $|Z|$ and the transverse radius $r:=|x|$.
	Ambient points are denoted by $Z=(x,y)$ or $X$, and ball centers by
	$Z_0$; in the fitting argument $z$ denotes the reference parameter
	in \eqref{eq:augmented-family}. A subscript $z$ in $\barC_z$,
	$D_z(N)$, or $W_z(Z)$ records the chosen reference. In contrast,
	$D_zW_z(Z)$ denotes differentiation in that parameter at fixed $Z$;
	$DW_z$ denotes differentiation in the ambient variable.
	We use $W_z$ for the regularized squared-distance function, extending
	the notation $W_a=U_a^2$ in the proof of \cite[Lemma~2.9]{FH1}.
	
	The area-ratio bound \eqref{equ:areabd} is invariant under the
	translations and dilations used below. In particular, for every fixed
	$d\ge0$, decomposition into unit annuli gives
	\begin{equation}\label{eq:gaussian-moment-bound}
		\begin{aligned}
			\int_N(1+|Z|)^d\rho\dd\mu_{N}
			\le c_n\sum_{j=0}^\infty(j+2)^d e^{-j^2/4}
			\mu_N(B_{j+1}) \le c_n\Lambda\sum_{j=0}^\infty(j+2)^d(j+1)^n e^{-j^2/4}
			\le C_{n,d,\Lambda}.
		\end{aligned}
	\end{equation}
	The same sum restricted to $j\ge R-1$ tends to zero as $R\to\infty$,
	uniformly in $N$. On each bounded
	axial region the area of $\{r<\delta\}$ is $O(\delta^m)$:
	cover the region by $O(\delta^{-k})$ axial balls of radius
	$\delta$ and use the area bound on the corresponding ambient
	balls of radius $2\delta$. By summing over unit axial cubes against
	the Gaussian, we obtain
	\begin{equation}\label{eq:raw-tube-area}
		\int_{N\cap\{r<\delta\}}\rho\dd\mu_{N}\le C\delta^m
		\qquad(0<\delta\le1).
	\end{equation}
	Since $m>2$, the corresponding Hardy moment is finite:
	\begin{equation}
		\begin{aligned}
			\int_Nr^{-2}\rho\dd\mu_{N}
			\le\int_{N\cap\{r\ge1\}}\rho\dd\mu_{N}
			+\sum_{j=0}^\infty 2^{2j+2}
			\int_{N\cap\{r<2^{-j}\}}\rho\dd\mu_{N} \le C+C\sum_{j=0}^\infty 2^{-j(m-2)}\le C.
		\end{aligned}
	\end{equation}
	The set $r=0$ has zero area, also by \eqref{eq:raw-tube-area}.
	These are absolute area estimates. The relative distance estimates
	needed after division by a vanishing normalization are proved in
	Section~\ref{sec:compactness}.
	
	%	\subsection{Analytic reduction of the link and the Gaussian Jacobi operator}
	
	Let $\C^m\subset\R^{m+1}$ be the transverse cone and write
	$\barC:=\C\times\R^k$, $n:=m+k$, $Z=(x,y)$, and $r:=|x|$.
	On the link $\Sigma$ use the sign convention that $\Delta_\Sigma$
	is nonpositive, and put
	\[
	\mathcal A_\Sigma:=-\Delta_\Sigma-|A_\Sigma|^2,\qquad
	L_\Sigma:=-\mathcal A_\Sigma+(m-1),\qquad
	K:=\ker L_\Sigma.
	\]
	The inner product on $K$ is the fixed, unweighted $L^2(\Sigma)$
	inner product. We will also use the \textbf{translation space}
	\[
	T_\C:=\{\langle e,\mathbf n_\Sigma\rangle:e\in\R^{m+1}\}.
	\]
	Indeed, we differentiate the minimal cone equation under ambient translation
	to obtain $(\Delta_\C+|A_\C|^2)\langle e,\mathbf n_\C\rangle=0$. The field
	is homogeneous of degree zero, so $T_\C\subset\ker\mathcal A_\Sigma$.
	
	For a small $w\in C^{2,\alpha}(\Sigma)$, $0<\alpha<1$, set
	\[
	\mathcal P_w(\omega)=\frac{\omega+w(\omega)\mathbf n_\Sigma(\omega)}
	{\sqrt{1+w(\omega)^2}},\qquad
	\mathcal A(w)=\mathcal H^{m-1}(\mathcal P_w(\Sigma)).
	\]
	Write $A_\Sigma X=-D_X\mathbf n_\Sigma$, $B_w=I-wA_\Sigma$,
	$s_w=\sqrt{1+w^2}$ and $v_w=(\mathbf n_\Sigma-w\omega)/s_w$.
	Here $v_w$ is a unit tangent vector to $S^m$ at $\mathcal P_w$,
	orthogonal to $T_\omega\Sigma$. Tangential differentiation gives
	\[
	d\mathcal P_w(X)=s_w^{-1}B_w X+s_w^{-2}(Xw)v_w .
	\]
	The induced metric is $s_w^{-2}B_w^2+s_w^{-4}\dd w\otimes\dd w$
	in an orthonormal frame of $\Sigma$. Since $B_w$ is positive
	definite for small $w$, the determinant formula for a rank-one
	perturbation gives its area Jacobian
	\[
	J_w=s_w^{-(m-1)}\det B_w\,\ell_w,\qquad
	\ell_w=\sqrt{1+s_w^{-2}|B_w^{-1}\nabla_\Sigma w|^2}.
	\]
	The unit normal of this graph, pulled back to $\Sigma$ and oriented
	continuously from $\mathbf n_\Sigma$, is
	$\mathbf n_w=(v_w-s_w^{-1}B_w^{-1}\nabla_\Sigma w)/\ell_w$.
	For a scalar graph variation $\varphi$ we have
	$D_w\mathcal P_w[\varphi]=s_w^{-2}\varphi v_w$. Take signed spherical
	mean curvature with $\mathbf H_{\Sigma_w}^{S^m}=H_{\Sigma_w}\mathbf n_w$
	at corresponding points. By the first variation formula, the gradient with respect to the fixed
	$L^2(\Sigma,\dd\mu_\Sigma)$ inner product is given by
	\begin{equation}\label{eq:cyl-link-area-gradient}
		\begin{aligned}
			D\mathcal A(w)[\varphi]
			&=\int_\Sigma\mathcal G(w)\varphi\,\dd\mu_\Sigma,\\
			\mathcal G(w)
			&=-(H_{\Sigma_w}\circ\mathcal P_w)\,\det B_w\,s_w^{-(m+1)},
			\qquad \Sigma_w=\mathcal P_w(\Sigma).
		\end{aligned}
	\end{equation}
	The multiplier of $-H_{\Sigma_w}$ is $J_w s_w^{-2}/\ell_w$, since
	the normal component of the variation is $\varphi s_w^{-2}/\ell_w$.
	
	For the transverse link $\Sigma$, we use the reduction of the spherical
	area from \cite[Proposition~2.5]{FH1}.
	
	\begin{prop}\label{prop:cyl-link-reduction}
		There is a real analytic function $\eta$ from a neighborhood of zero in
		$K$ to $C^{2,\alpha}(\Sigma)\cap K^\perp$ with
		$\eta(0)=0$ and $D\eta(0)=0$. For $a$ in this neighborhood,
		$P_{K^\perp}\mathcal G(a+\eta(a))=0$. Set
		\begin{equation}\label{eq:reduced-link-data}
			k(a)=a+\eta(a),\qquad f(a)=\mathcal A(k(a)),\qquad
			q(a)=\nabla f(a)=\mathcal G(k(a))\in K .
		\end{equation}
		Then $q(0)=0$, $Dq(0)=0$ and, after decreasing the parameter
		neighborhood, there are $C<\infty$ and $0<\theta\le1/2$ such that
		\begin{equation}\label{eq:reduced-link-gradient}
			|f(a)-f(0)|^{1-\theta}\le C|q(a)|.
		\end{equation}
		The links $\Sigma_a=\mathcal P_{k(a)}(\Sigma)$ are smooth and
		converge smoothly to $\Sigma$ as $a\to0$.
	\end{prop}
	\begin{proof}
		The analytic map $\eta$, the gradient identity \eqref{eq:reduced-link-data},
		and the inequality \eqref{eq:reduced-link-gradient} follow from
		\cite[Proposition~2.5]{FH1}. Since $\Sigma$ is minimal, we have
		$q(0)=\mathcal G(0)=0$. Differentiating $q=\mathcal G\circ k$ at zero
		and using $D\eta(0)=0$, we obtain
		\[
		Dq(0)[\psi]=-L_\Sigma\psi=0\qquad(\psi\in K).
		\]
		As $a\to0$, $k(a)\to0$ in $C^{2,\alpha}(\Sigma)$ and $q(a)\to0$
		in every $C^j$ norm, since $K$ is a finite-dimensional space of smooth
		functions. Schauder bootstrapping in $\mathcal G(k(a))=q(a)$ gives
		$k(a)\to0$ smoothly, hence $\Sigma_a\to\Sigma$ smoothly.
	\end{proof}
	
	Define $\C_a=\{r\mathcal P_{k(a)}(\omega):r\ge0,\ \omega\in\Sigma\}$
	and $\barC_a=\C_a\times\R^k$. For an unrotated reference set
	$j_a(Z)=\dist(Z,\barC_a)$ and $D_a(N)=D_{(a,0)}(N)$ as in
	\eqref{eq:reference-distance}. Since $\mathcal P_{k(a)}(\omega)$
	lies on the unit sphere, $r=|x|$ at
	$(x,y)=(r\mathcal P_{k(a)}(\omega),y)$. The radial and axial
	principal curvatures of the cylinder vanish. Its other principal
	curvatures are those of $\Sigma_a\subset S^m$ divided by $r$.
	Consequently \eqref{eq:cyl-link-area-gradient} gives the signed curvature
	\begin{equation}\label{eq:reduced-cone-curvature}
		H_{\barC_a}(r\mathcal P_{k(a)}(\omega),y)
		=-\frac{q(a)(\omega)}{r\mathfrak w_a(\omega)},\qquad
		\mathfrak w_a=\det(I-k(a)A_\Sigma)(1+k(a)^2)^{-(m+1)/2}.
	\end{equation}
	In particular $\mathfrak w_a>0$ and $\mathfrak w_a\to1$ smoothly.
	
	The normal projection of the position vector vanishes, and the axial
	Gaussian has integral one. Fubini and polar integration on $\C_a$ give
	\begin{equation}\label{eq:cyl-product-functionals}
		F(\barC_a)=F(\C_a)=\frac{f(a)}{\sigma_{m-1}},\qquad
		G(\barC_a)=G(\C_a),\qquad \sigma_{m-1}=\mathcal H^{m-1}(S^{m-1}).
	\end{equation}
	The constant in the first formula follows by carrying out the same
	radial integral on an $m$-plane. For the second, the squared quantity is
	\[
	G(\C_a)^2=(4\pi)^{-m/2}
	\left(\int_0^\infty r^{m-3}e^{-r^2/4}\dd r\right)
	\int_\Sigma |q(a)|^2\,\mathfrak w_a^{-2}J_{k(a)}\,\dd\mu_\Sigma .
	\]
	The radial integral is finite for $m>2$. The last two factors in
	the angular density are uniformly bounded above and below. Hence
	$C^{-1}|q(a)|\le G(\barC_a)\le C|q(a)|$, and
	\eqref{eq:reduced-link-gradient} controls
	$|F(\barC_a)-F(\barC)|^{1-\theta}$ by $C|q(a)|$. These Gaussian
	identities are invariant under a simultaneous rotation of the reference.

	Write $c_d:=(4\pi)^{-d/2}$, $\rho_d(z):=c_d e^{-|z|^2/4}$ in dimension $d$,
	$\rho:=\rho_m\rho_k$, and $\dd g:=\rho\dd\mu_{\barC}$. On
	$\Reg(\barC)$, put
	\[
	J=\Delta_{\barC}+|A_{\barC}|^2,\qquad
	D=Z\cdot\nabla_{\barC}=r\partial_r+y\cdot\nabla_y,\qquad
	L_{\barC}=J-\tfrac12D+\tfrac12 .
	\]
	We use a real orthonormal angular eigenbasis $\{\phi_j\}$ of
	$\mathcal A_\Sigma$, with corresponding eigenvalues $\mu_j$. Set
	\[
	\alpha_0=(m-2)/2,\qquad \nu_j=\sqrt{\alpha_0^2+\mu_j},
	\qquad \gamma_j^+=-\alpha_0+\nu_j .
	\]
	By strict Hardy stability, $\nu_j\ge\nu_0>0$ for all $j$.
	The identity $\gamma_j^+(\gamma_j^++m-2)=\mu_j$ follows from
	the definition. We use the tube threshold
	\[
	\tau_0:=\max\{0,2-2\nu_0\}<2.
	\]
	
	\begin{prop}\label{prop:spectrum}
		For $l\ge0$ and $\alpha\in\mathbb Z_{\ge0}^k$, define
		\begin{equation}\label{eq:cyl-gaussian-spectrum}
			\begin{aligned}
				e_{jl\alpha}(r,\omega,y)
				&=r^{\gamma_j^+}L_l^{\nu_j}(r^2/4)\phi_j(\omega)H_\alpha(y),\\
				L_{\barC}e_{jl\alpha}&=\lambda_{jl\alpha}e_{jl\alpha}
				\quad\text{on }\Reg(\barC),\qquad
				\lambda_{jl\alpha}=\tfrac12(1-\gamma_j^+-2l-|\alpha|).
			\end{aligned}
		\end{equation}
		Here $L_l^\nu$ is the generalized Laguerre polynomial and $H_\alpha$
		is the product of the monic Hermite polynomials for $\rho_k$.
		Both $e_{jl\alpha}$ and $De_{jl\alpha}$ belong to $L^2(\dd g)$.
		The complex bilinear moments
		\[
		[U,e_{jl\alpha}]_\rho=\int_{\barC} Ue_{jl\alpha}\dd g
		\]
		determine each $U\in L^2(\dd g;\mathbb C)$ uniquely. For each fixed $j$, the radial and axial factors give an orthogonal
		basis of its angular component.
	\end{prop}
	\begin{proof}
		The monic Hermite polynomials are
		$h_s(t)=(-2)^s e^{t^2/4}(\dd/\dd t)^s e^{-t^2/4}$ and
		$H_\alpha(y)=\prod_{i=1}^k h_{\alpha_i}(y_i)$. Their differential
		equation and the product rule give
		\begin{equation}
			(\Delta_y-\tfrac12y\cdot\nabla_y)H_\alpha
			=-\tfrac12|\alpha|H_\alpha .
		\end{equation}
		They form a complete orthogonal family for $\rho_k\dd y$, by
		the one-dimensional Hermite expansion and tensor products
		\cite[Section~5.5 and Theorem~5.7.1]{Szego}.
		For the radial factor, set $\gamma=\gamma_j^+$ and $t=r^2/4$.
		On $r^\gamma P(t)\phi_j$, the transverse part of $L_{\barC}$ acts
		by applying
		\[
		P\longmapsto tP''+(1+\nu_j-t)P'+\tfrac12(1-\gamma)P .
		\]
		This follows by cancelling the $r^{\gamma-2}$ coefficient with the
		indicial equation. The polynomial
		\[
		L_l^\nu(t)=
		\sum_{h=0}^l\frac{(-1)^h}{h!}
		\frac{\Gamma(l+\nu+1)}{(l-h)!\,\Gamma(h+\nu+1)}t^h
		\]
		solves $tP''+(1+\nu-t)P'=-lP$. Its leading coefficient is
		$(-1)^l/l!$. This proves \eqref{eq:cyl-gaussian-spectrum}.
		Near $r=0$, the squared radial norm has weight
		$r^{m-1+2\gamma}=r^{1+2\nu_j}$.
		The radial Euler derivative preserves this order of vanishing, and
		the axial Euler derivative acts on a polynomial. Gaussian decay gives
		integrability at infinity in both cases.
		
		Multiplication by $r^\gamma$ and the substitution
		$t=r^2/4$ identify the radial norm, up to a positive constant, with
		that of $L^2(t^{\nu_j}e^{-t}\dd t)$. This
		identification is onto by division by $r^\gamma$ for $r>0$.
		Polynomials are dense in this space. Indeed, if $v$ is orthogonal to
		all powers and $w(t)=t^{\nu_j}e^{-t}$, the transform
		$\int_0^\infty v(t)e^{zt}w(t)\dd t$ is holomorphic for
		$\operatorname{Re}z<1/2$ by the Cauchy--Schwarz inequality. Its derivatives at zero
		all vanish. Fourier uniqueness on the imaginary axis gives $vw=0$.
		The Laguerre polynomials span the polynomials by their leading
		coefficients. Laguerre orthogonality with weight
		$t^{\nu_j}e^{-t}$ follows by integration by parts in Rodrigues'
		formula. We tensor the radial families with the Hermite family and use
		the real angular basis to show that the span of the tests is dense. Complex conjugation
		preserves density, so vanishing of all bilinear moments forces $U=0$.
	\end{proof}
	
	\begin{lem}\label{lem:all-root-moment-test}
		Let $U$ be smooth on $\Reg(\barC)$ and belong to $L^2(\dd g)$.
		Suppose that for some $p>\tau_0$ and each $R\ge1$,
		\begin{equation}\label{eq:gaussian-moment-tube}
			\int_{\barC\cap\{r<\delta,\ |Z|<R\}}|U|^2\dd g
			\le C_R\delta^p\qquad(0<\delta<1).
		\end{equation}
		If $L_{\barC}U=f$ on $\Reg(\barC)$, where $f\in L^2(\dd g)$,
		then for every index $j,l,\alpha$ in \eqref{eq:cyl-gaussian-spectrum},
		\[
		\lambda_{jl\alpha}[U,e_{jl\alpha}]_\rho=[f,e_{jl\alpha}]_\rho .
		\]
	\end{lem}
	\begin{proof}
		Let $\chi_\delta(r)$ vanish for $r\le\delta$ and equal one for
		$r\ge2\delta$, with derivatives of order $h=1,2$ bounded by
		$C_h\delta^{-h}$. Put $\zeta_R=\zeta(|Z|^2/R^2)$, equal to one on
		$B_R$ and zero outside $B_{2R}$. For $e=e_{jl\alpha}$, test the
		punctured equation with $\chi_\delta\zeta_Re$ and integrate by
		parts using the Gaussian measure and the bilinear pairing. For either
		cutoff $\chi$ the commutator is
		\begin{equation}\label{eq:cyl-cutoff-commutator}
			[L_{\barC},\chi]v=2\nabla_{\barC}\chi\cdot\nabla_{\barC}v+
			(\Delta_{\barC}\chi-\tfrac12D\chi)v .
		\end{equation}
		At fixed $R$, on the inner shell
		$|[L_{\barC},\chi_\delta](\zeta_Re)|\le
		C_{R,e}\delta^{\gamma_j^+-2}$.
		The Cauchy--Schwarz inequality, \eqref{eq:gaussian-moment-tube}, and the
		$O_R(\delta^m)$ cone measure of that shell give
		\begin{equation}
			\left|\int_{\barC}U[L_{\barC},\chi_\delta](\zeta_Re)\dd g\right|
			\le C_{R,e}\delta^{\gamma_j^+-2+(m+p)/2}
			=C_{R,e}\delta^{\nu_j+(p-2)/2}\longrightarrow0.
		\end{equation}
		Indeed $\nu_j+(p-2)/2\ge\nu_0+(p-2)/2>0$.
		For the outer commutator the derivative term is
		$4R^{-2}\zeta'(|Z|^2/R^2)De$, since the position vector is tangent
		to the regular cylinder. Also $\Delta_{\barC}|Z|^2=2n$, so
		$|\Delta_{\barC}\zeta_R-\tfrac12D\zeta_R|\le C$ on its support.
		After sending $\delta$ to zero, the remaining error is bounded by
		\[
		C\|U\|_{L^2(\dd g;\,|Z|>R)}
		\bigl(R^{-2}\|De\|_{L^2(\dd g)}+\|e\|_{L^2(\dd g)}\bigr).
		\]
		This tends to zero as $R\to\infty$. The integrals containing $Ue$
		and $fe$ converge by the Cauchy--Schwarz inequality, which proves the identity.
	\end{proof}
	
	Denote the space generated by the reduced link parameters and
	the rotations mixing the two factors by
	\begin{equation}\label{eq:cyl-gaussian-kernel}
		\mathcal G_{\barC}:=rK\oplus\bigoplus_{i=1}^k y_iT_\C .
	\end{equation}
	The following proposition identifies the possible limits of normalized
	graph heights in the presence of the tube bound.
	\begin{prop}\label{prop:gaussian-bordered}
		Assume \eqref{eq:cyl-S1}, $\ker\mathcal A_\Sigma=T_\C$, and
		$\gamma_j^+\notin\{-1,-2,\ldots\}$ for every $j$.
		Let $U$ satisfy the regularity, integrability and tube hypotheses of
		Lemma~\ref{lem:all-root-moment-test}. If
		$L_{\barC}U=r^{-1}\xi$ on $\Reg(\barC)$ for some $\xi\in K$,
		then $\xi=0$ and $U\in\mathcal G_{\barC}$. In particular,
		$U\perp\mathcal G_{\barC}$ in $L^2(\dd g)$ implies $U=0$.
	\end{prop}
	\begin{proof}
		Polar coordinates give
		\begin{equation}\label{eq:residual-square-integrability}
			\|r^{-1}\xi\|_{L^2(\dd g)}^2
			=c_m\|\xi\|_{L^2(\Sigma)}^2
			\int_0^\infty r^{m-3}e^{-r^2/4}\dd r<\infty ,
		\end{equation}
		since $m>2$. For $\psi\in K$ the function $r\psi$ is a linear
		combination of the tests with $\gamma_j^+=1$ and $l=|\alpha|=0$.
		Their eigenvalues vanish. Lemma~\ref{lem:all-root-moment-test}
		therefore gives
		\[
		0=[r^{-1}\xi,r\psi]_\rho
		=c_m\left(\int_0^\infty r^{m-1}e^{-r^2/4}\dd r\right)
		\int_\Sigma\xi\psi\,\dd\mu_\Sigma .
		\]
		For real data, we take $\psi=\xi$ and obtain $\xi=0$. For complex data
		take $\psi=\overline\xi$. By applying the lemma again, we find that
		all moments for which $\lambda_{jl\alpha}\ne0$ vanish.
		
		For each angular component the Laguerre--Hermite factors are an
		orthogonal basis by Proposition~\ref{prop:spectrum}. Thus only terms
		with $2l+|\alpha|=1-\gamma_j^+$ remain. In particular, the larger
		indicial root is an integer at most one. By nonresonance it equals
		one or zero. The first case gives
		$rK$. In the second case $l=0$, $|\alpha|=1$ and
		$\phi_j\in T_\C$, giving the mixed rotations in
		\eqref{eq:cyl-gaussian-kernel}. This proves the proposition.
	\end{proof}
	
	We now enlarge the reduced family to include the mixed-rotation fields
	in \eqref{eq:cyl-gaussian-kernel}. For
	$\mathsf B\in\operatorname{Hom}(\R^k,\R^{m+1})$, define
	\[
	S_{\mathsf B}:=\begin{pmatrix}0&\mathsf B\\-\mathsf B^\ast&0\end{pmatrix},\qquad
	S_{\mathsf B}(x,y):=(\mathsf B y,-\mathsf B^\ast x),\qquad S_{\mathsf B}^\ast=-S_{\mathsf B}.
	\]
	Here $\mathsf B^\ast$ is the Euclidean adjoint (matrix transpose) of $\mathsf B$. Consequently $R_{\mathsf B}:=\exp(S_{\mathsf B})\in SO(n+1)$.
	Define $\mathscr T(\mathsf B)(\omega,y):=\langle \mathsf B y,\mathbf n_\Sigma(\omega)\rangle$
	as a linear map into smooth functions on $\Sigma\times\R^k$. Thus
	\[
	\ker\mathscr T=\{\mathsf B:\langle\mathsf B y,\mathbf n_\Sigma(\omega)\rangle=0
	\text{ for all }(\omega,y)\in\Sigma\times\R^k\}.
	\]
	Fix a linear subspace $\mathscr E$ such that
	\[
	\operatorname{Hom}(\R^k,\R^{m+1})=\ker\mathscr T\oplus \mathscr E.
	\]
	We take $\mathsf B\in \mathscr E$, so that $\mathscr T|_{\mathscr E}$ is injective. For
	$z=(a,\mathsf B)\in K\times \mathscr E$ near zero, set
	\begin{equation}\label{eq:augmented-family}
		R_{\mathsf B}:=\exp(S_{\mathsf B}),\qquad \barC_z:=R_{\mathsf B}\barC_a,\qquad q(z):=q(a).
	\end{equation}
	A parametrization of the regular reference is
	\[
	X_z(r,\omega,y):=R_{\mathsf B}(r\mathcal P_{k(a)}(\omega),y).
	\]
	At $z=0$ its unit normal is $\mathbf n_0=(\mathbf n_\Sigma,0)$.
	Since $Dk(0)[\psi]=\psi$ for $\psi\in K$, differentiation gives
	\begin{equation}\label{eq:parameter-normal-velocities}
		\begin{aligned}
			\langle D_aX_0[\psi],\mathbf n_0\rangle&=r\psi(\omega),\\
			\langle D_{\mathsf B}X_0[\mathsf B],\mathbf n_0\rangle
			&=\langle S_{\mathsf B}(r\omega,y),(\mathbf n_\Sigma,0)\rangle
			=\langle \mathsf B y,\mathbf n_\Sigma\rangle
			=\sum_{\ell=1}^k y_\ell\langle \mathsf B e_\ell,\mathbf n_\Sigma\rangle.
		\end{aligned}
	\end{equation}
	The image of the second map is
	$\bigoplus_{\ell=1}^k y_\ell T_\C$; restricting to $\mathscr E$ does not
	change this image. Thus the parameter velocities span
	$\mathcal G_{\barC}$ in \eqref{eq:cyl-gaussian-kernel}. For
	$\mathsf A\in\mathfrak{so}(m+1)$, differentiation along the minimal links
	$e^{t\mathsf A}\Sigma$ gives
	\[
	L_\Sigma\langle \mathsf A\omega,\mathbf n_\Sigma\rangle=0,
	\qquad \langle \mathsf A\omega,\mathbf n_\Sigma\rangle\in K.
	\]
	Rotations within the axial factor have zero normal component and
	preserve the cylinder.
	
	We verify that this family is locally injective. Choose a basis
	$\psi_1,\ldots,\psi_{d_K}$ of $K$ and a basis $\mathsf B_1,\ldots,\mathsf B_{d_E}$
	of $\mathscr E$. Write $z:=(z^1,\ldots,z^{d_{\rm par}})$ in these bases,
	where $d_{\rm par}:=d_K+d_E$ and $d_E:=\dim\mathscr E$.
	The corresponding normal velocities are
	\[
	g_j:=r\psi_j\quad(1\le j\le d_K),\qquad
	g_{d_K+\beta}:=\langle \mathsf B_\beta y,\mathbf n_\Sigma\rangle
	\quad(1\le\beta\le d_E).
	\]
	Choose $0<r_1<r_2$ and $L>0$, and a nonnegative
	$\chi\in C_c^\infty(\Reg(\barC))$ supported in
	$\Omega:=\{r_1<r<r_2,\ \omega\in\Sigma,\ |y|<L\}$,
	positive on a smaller product set containing the full link.
	Define the $d_{\rm par}\times d_{\rm par}$ matrix
	\[
	M_{ij}:=\int_{\barC}\chi g_i g_j\rho\dd\mu_{\barC}.
	\]
	For a coefficient vector $c$, let $\psi\in K$ and $\mathsf B\in \mathscr E$ be the
	corresponding linear combinations. Then
	\begin{equation}\label{eq:parameter-gram-matrix}
		c^\top Mc
		=\int_{\barC}\chi\bigl(r\psi(\omega)
		+\langle \mathsf B y,\mathbf n_\Sigma(\omega)\rangle\bigr)^2\rho\dd\mu_{\barC}
		\ge0.
	\end{equation}
	If equality holds, continuity makes the integrand vanish wherever
	$\chi>0$. We evaluate at $y=0$ to obtain $\psi=0$ on $\Sigma$.
	We then differentiate in $y_\ell$ to obtain
	$\langle \mathsf B e_\ell,\mathbf n_\Sigma\rangle=0$ for every $\ell$ and $\omega$.
	Thus $\mathsf B\in \mathscr E\cap\ker\mathscr T=\{0\}$, and $c=0$.
	Therefore $M$ is positive definite.
	
	For small $z$, write $\barC_z$ as a normal graph of height $h_z$
	over a neighborhood of $\supp\chi$ in $\barC$, and put
	\[
	\mathcal F_i(z):=\int_{\barC}\chi h_zg_i\rho\dd\mu_{\barC},
	\qquad 1\le i\le d_{\rm par}.
	\]
	The graph parametrization and \eqref{eq:parameter-normal-velocities}
	give $h_0=0$, $\partial_{z^j}h_z|_{z=0}=g_j$, and hence
	\[
	\partial_{z^j}\mathcal F_i(0)
	=\int_{\barC}\chi g_jg_i\rho\dd\mu_{\barC}=M_{ij}.
	\]
	The inverse function theorem gives a neighborhood of zero on which
	$\mathcal F$ is injective. If $\barC_z=\barC_{z'}$ for parameters
	in this neighborhood, uniqueness of normal projection implies
	$h_z=h_{z'}$ on $\supp\chi$, so $\mathcal F(z)=\mathcal F(z')$
	and $z=z'$.
	
	Orthogonal invariance of the Gaussian and the induced measure gives
	\begin{equation}
		F(\barC_z)
		=\int_{\barC_a}(4\pi)^{-n/2}e^{-|R_{\mathsf B}Z|^2/4}\dd\mu_{\barC_a}(Z)
		=F(\barC_a)=F(\C_a).
	\end{equation}
	With the normal transported by $R_{\mathsf B}$, scalar mean curvature satisfies
	\begin{equation}
		H_{\barC_z}(X_z(r,\omega,y))
		=H_{\barC_a}(r\mathcal P_{k(a)}(\omega),y)
		=-\frac{q(a)(\omega)}{r\mathfrak w_a(\omega)}.
	\end{equation}
	Since $K$ is finite dimensional and $\mathfrak w_a$ is uniformly
	positive with uniformly bounded derivatives, for every integer $j\ge0$, we have
	\[
	\left\|\nabla_\Sigma^j\bigl(q(a)/\mathfrak w_a\bigr)
	\right\|_{C^0(\Sigma)}\le C_j|q(a)|.
	\]
	In particular, suppose $z_i\to0$, $b_i>0$, and $q(a_i)/b_i\to\xi$
	in $K$. For every compact parameter set $\mathcal K$ with $r>0$, we have
	\begin{equation}
		\left\|b_i^{-1}H_{\barC_{z_i}}\circ X_{z_i}
		+r^{-1}\xi\right\|_{C^j(\mathcal K)} \le C_{j,\mathcal K}
		\left(\|q(a_i)/b_i-\xi\|_{C^j(\Sigma)}
		+\|\xi\|_{C^j(\Sigma)}
		\|\mathfrak w_{a_i}^{-1}-1\|_{C^j(\Sigma)}\right)
		\longrightarrow0.
	\end{equation}
	This is the reference-curvature limit used in the normalized graph equation.
	
	In the next lemma, we extend the construction for a single cone in \cite[Lemma~2.3]{FH1}
	to the family of cylinders $\barC_z$ by pulling back the transverse
	function under the rotation $R_{\mathsf B}$.
	\begin{lem}\label{lem:distance-function}
		There are nonnegative functions $W_z$ uniformly comparable to
		$\operatorname{dist}(\cdot,\barC_z)^2$ and equal to that squared
		distance in a fixed narrow conical normal tube of $\barC_z$.
		In coordinates aligned with $\barC_z$, $W_z$ depends only on the
		transverse variable, is smooth off the reference axis, and
		extends as a $C^{1,1}$ function across the axis. Uniformly for small $z$, we have
		\[
		W_z\le Cr^2,\qquad |DW_z|^2\le CW_z,\qquad
		|D_zW_z|\le C|Z|\sqrt{W_z},
		\]
		where $r$ in the first estimate is the aligned transverse
		radius. On a fixed compact regular normal graph of height $u$
		over $\barC_z$, differentiation with respect to $z^j$ at a fixed ambient point gives
		\[
		\partial_{z^j}W_z=-2u\,g_{z,j}+O(u^2)
		\]
		with $g_{z,j}$ the normal velocity of the reference
		parametrization in the $z^j$ direction.
	\end{lem}
	\begin{proof}
		Let $\delta_a$ denote signed spherical distance to $\Sigma_a$
		on a uniform tubular neighborhood. Choose a smooth signed defining
		function $\varphi_a$ on $S^m$, depending smoothly on $a$, with
		\[
		\varphi_a^{-1}(0)=\Sigma_a,\qquad
		\varphi_a=\sin\delta_a\text{ near }\Sigma_a.
		\]
		To obtain a common family, we carry a fixed defining function to the
		nearby links by a smooth family of diffeomorphisms and interpolate
		with $\sin\delta_a$ inside a slightly larger tubular neighborhood.
		The two functions have the same sign on each side, so this
		interpolation introduces no additional zeros. Define
		\begin{equation}\label{eq:positional-function-definition}
			U_a(r\omega):=r \varphi_a(\omega),\qquad U_a(0):=0,
			\qquad \pi_x(x,y):=x,
			\qquad W_z(Z):=U_a(\pi_xR_{\mathsf B}^{-1}Z)^2.
		\end{equation}
		In the smaller normal neighborhood, elementary cone geometry gives
		$U_a=\operatorname{sdist}(\cdot,\C_a)$.
		Compactness of the sphere, transversality, and homogeneity give
		\begin{equation}\label{eq:homogeneous-distance-bounds}
			\begin{gathered}
				c\dist(x,\C_a)\le |U_a(x)|\le C\dist(x,\C_a)\le C|x|,\\
				|DU_a(x)|\le C,\qquad
				|D^2U_a(x)|\le C|x|^{-1},\qquad
				|\partial_aU_a(x)|\le C|x|\quad(x\ne0).
			\end{gathered}
		\end{equation}
		Consequently, for $x\ne0$, we have
		\[
		|D(U_a^2)|=2|U_aDU_a|\le C|x|,\qquad
		|D^2(U_a^2)|\le2|DU_a|^2+2|U_a||D^2U_a|\le C.
		\]
		The value and gradient extend to zero with gradient zero.
		The uniform Hessian bound, integrated along line segments (splitting
		at zero if necessary), shows that this extension is $C^{1,1}$.
		Thus $W_z$ has the asserted regularity across its reference axis.
		
		Write $x_{\mathsf B}:=\pi_xR_{\mathsf B}^{-1}Z$ and $r:=|x_{\mathsf B}|$. Since $R_{\mathsf B}$ is orthogonal,
		\eqref{eq:homogeneous-distance-bounds} gives
		\[
		W_z\le Cr^2,\qquad |DW_z|^2\le C W_z.
		\]
		Smoothness of the matrix exponential on the fixed parameter
		neighborhood gives $|D_{\mathsf B}x_{\mathsf B}|\le C|Z|$. Therefore we obtain
		\begin{equation}\label{eq:distance-parameter-bounds}
			\begin{aligned}
				|D_aW_z|&\le2|U_a(x_{\mathsf B})||D_aU_a(x_{\mathsf B})|
				\le Cr\sqrt{W_z}\le C|Z|\sqrt{W_z},\\
				|D_{\mathsf B}W_z|&\le2|U_a(x_{\mathsf B})||DU_a(x_{\mathsf B})||D_{\mathsf B}x_{\mathsf B}|
				\le C|Z|\sqrt{W_z}.
			\end{aligned}
		\end{equation}
		At $x_{\mathsf B}=0$, a parameter increment $\dot z$ satisfies
		$W_{z+\dot z}(Z)\le C|\dot z|^2|Z|^2$, so $D_zW_z(Z)=0$
		there, and \eqref{eq:distance-parameter-bounds} remains valid.
		
		Finally, let $l_z$ be signed Euclidean distance to $\barC_z$
		in a fixed compact regular normal neighborhood. There $W_z=l_z^2$.
		If $X_z(p)$ parametrizes the reference and $\mathbf n_z$ is its unit normal,
		then we have
		\[
		l_z(X_z(p))=0,\qquad Dl_z(X_z(p))=\mathbf n_z(p),\qquad
		g_{z,j}(p):=\langle\partial_{z^j}X_z(p),\mathbf n_z(p)\rangle.
		\]
		By differentiating the first identity with respect to $z^j$, we obtain
		\[
		(\partial_{z^j}l_z)(X_z(p))
		=-Dl_z(X_z(p))\cdot\partial_{z^j}X_z(p)=-g_{z,j}(p).
		\]
		For a graph point $X:=X_z(p)+u(p)\mathbf n_z(p)$, we have $l_z(X)=u(p)$.
		Smoothness of $l_z$ on a fixed larger normal neighborhood and Taylor's
		formula, with parameter differentiation at fixed ambient $X$, give
		\[
		\begin{aligned}
			|(\partial_{z^j}l_z)(X)+g_{z,j}(p)|&\le C|u(p)|,\\
			|(\partial_{z^j}W_z)(X)+2u(p)g_{z,j}(p)|
			&=2|u(p)|\,|(\partial_{z^j}l_z)(X)+g_{z,j}(p)|
			\le C|u(p)|^2.
		\end{aligned}
		\]
		This proves the claimed parameter expansion.
	\end{proof}

	Let $Y$ be a fixed smooth reference hypersurface on a
	relatively compact coordinate domain, with metric $g$, unit normal
	$\mathbf n_Y$, and second fundamental form $A_Y$, with components
	$A_{ij}:=\langle D_iD_jY,\mathbf n_Y\rangle$. For $X_u:=Y+u \mathbf n_Y$, put
	$\mathsf T_u:=I-uA_Y^\sharp$ and
	$\mathsf v_u:=(1+|\mathsf T_u^{-1}\nabla_Yu|^2)^{1/2}$. For a sufficiently small
	$C^1$ graph, direct differentiation gives
	\begin{equation}
		\begin{aligned}
			\partial_iX_u&=\mathsf T_u\partial_iY+(\partial_i u) \mathbf n_Y,\\
			(g_u)_{ij}&=g_{ij}-2uA_{ij}
			+u^2 A_i{}^lA_{lj}+(\partial_i u)(\partial_j u),\\
			\mathbf n_u&=\mathsf v_u^{-1}(\mathbf n_Y-\mathsf T_u^{-1}\nabla_Yu),\qquad
			\langle \mathbf n_Y,\mathbf n_u\rangle=\mathsf v_u^{-1},\\
			H(u)&=(g_u)^{ij}\langle D_iD_jX_u,\mathbf n_u\rangle.
		\end{aligned}
	\end{equation}
	Here $H(u)$ is scalar mean curvature with respect to $\mathbf n_u$.
	For a variation $u=t\varphi$, the metric and normal derivatives at zero are
	\[
	\left.\frac{\dd}{\dd t}\right|_0(g_{t\varphi})^{ij}
	=2\varphi A^{ij},\qquad
	\left.\frac{\dd}{\dd t}\right|_0 \mathbf n_{t\varphi}=-\nabla_Y\varphi.
	\]
	Using $\langle D_iD_j\mathbf n_Y,\mathbf n_Y\rangle=-A_i{}^lA_{lj}$ and
	$(D_iD_jY)^\top=\Gamma_{ij}^l\partial_lY$, we obtain
	\begin{equation}\label{eq:graph-mean-curvature-linearization}
		DH(0)[\varphi]
		=2\varphi|A_Y|^2
		+g^{ij}(\partial_{ij}\varphi-\varphi A_i{}^lA_{lj})
		-g^{ij}\Gamma_{ij}^l\partial_l\varphi =\Delta_Y\varphi+|A_Y|^2\varphi.
	\end{equation}
	
	Consider a flow whose normal velocity is
	$\mathbf H+\tfrac12(v+\epsilon X)^\perp$, where $v$ is an ambient
	vector and $\epsilon$ a scalar, independent of the spatial variable.
	They may depend smoothly on time. Define
	\[
	\mathcal D_{v,\epsilon}(u)
	:=\tfrac12\langle v+\epsilon X_u,\mathbf n_u\rangle.
	\]
	The same normal variation gives
	\begin{equation}\label{eq:graph-drift-linearization}
		D\mathcal D_{v,\epsilon}(0)[\varphi]
		=\tfrac12\left(\epsilon\varphi
		-\langle(v+\epsilon Y)^\top,\nabla_Y\varphi\rangle\right).
	\end{equation}
	Since $Y$ is fixed in time, the normal component of $\partial_sX_u$
	is $\mathsf v_u^{-1}\partial_su$. Thus the exact scalar graph equation is
	\begin{equation}\label{eq:normal-graph-flow}
		\partial_su=\mathscr F_{v,\epsilon}(u),\qquad
		\mathscr F_{v,\epsilon}(u)
		:=\mathsf v_u\bigl(H(u)+\mathcal D_{v,\epsilon}(u)\bigr).
	\end{equation}
	We have $\mathsf v_0=1$ and $D\mathsf v_0[\varphi]=0$. Therefore
	\eqref{eq:graph-mean-curvature-linearization} and
	\eqref{eq:graph-drift-linearization} yield
	\begin{equation}\label{eq:graph-flow-linearization}
		D\mathscr F_{v,\epsilon}(0)[\varphi]
		=\Delta_Y\varphi+|A_Y|^2\varphi
		-\tfrac12\langle(v+\epsilon Y)^\top,\nabla_Y\varphi\rangle
		+\tfrac\epsilon2\varphi.
	\end{equation}
	If $Y$ is conical and $v$ is tangent to $Y$, then
	$\mathcal D_{v,\epsilon}(0)=0$ and $\mathscr F_{v,\epsilon}(0)=H_Y$.
	The fundamental theorem of calculus in the graph variable gives
	\begin{equation}\label{eq:averaged-graph-equation}
		\partial_su-H_Y
		=\mathscr F_{v,\epsilon}(u)-\mathscr F_{v,\epsilon}(0)
		=\int_0^1D\mathscr F_{v,\epsilon}(tu)[u]\dd t.
	\end{equation}
	Define $\mathscr L_u:=\int_0^1D\mathscr F_{v,\epsilon}(tu)\dd t$.
	Then \eqref{eq:averaged-graph-equation} reads
	$\partial_su-H_Y=\mathscr L_u u$. On each fixed regular coordinate domain,
	$\mathscr L_u=a^{ij}_u\partial_{ij}+b^i_u\partial_i+c_u$ is linear
	in its argument once $u$ is fixed. Its principal coefficients are
	$a^{ij}_u:=\int_0^1(g_{tu})^{ij}\dd t$, so for small $\|u\|_{C^1}$, we have
	\[
	c\,g^{ij}\zeta_i\zeta_j
	\le a^{ij}_u\zeta_i\zeta_j
	\le C\,g^{ij}\zeta_i\zeta_j.
	\]
	The lower-order coefficients are smooth functions of the reference
	geometry, $u$, and its derivatives. Hence, when the references converge
	smoothly and $u\to0$ smoothly, these coefficients converge on compact
	regular subsets to those in \eqref{eq:graph-flow-linearization}.
	For $(v,\epsilon)=(0,1)$ and $Y=\barC$, that operator is
	$L_{\barC}$. For the dilations used in Section~\ref{sec:compactness},
	$v$ is axial and thus tangent to the reference cylinder, so
	\eqref{eq:averaged-graph-equation} retains precisely the reference
	mean curvature as its inhomogeneous term.

	Let $\barC_a:=\C_a^m\times\mathbb R^k\subset\mathbb R^{m+1}\times\mathbb R^k$,
	$m>2$, $k\ge1$, be a reference cylinder over a small deformation of the two-sided regular cone.
	The reference need not be minimal. Its links have uniform smooth
	geometry, and their spherical mean curvature is bounded by $Cb$, where $b:=|q(a)|$, in
	the fixed smooth norms. Write $n:=m+k$, $Z=(x,y)$, $r:=|x|$, $r_y:=|y|$,
	and write $\mathbf n_N:=(\mathbf n_x,\mathbf n_y)$ for a local unit
	normal to $N$, with $\mathbf n_y\in\R^k$ its axial component.
	We use $\rho$, $\Phi$, $D_a(N)$, and $G(N)$ as defined in
	\eqref{eq:flow}--\eqref{eq:reference-distance}.
	By \eqref{eq:gaussian-moment-bound}, all Gaussian moments of $N$
	are uniformly bounded.
	An infinite $G(N)^2$ makes the inequalities below automatic.
	
	\begin{lem}\label{lem:weighted-normal-angle}
		Fix the area-ratio bound $\Lambda$. For every sufficiently small reduced
		reference $\barC_a$, define
		$\mathcal T_{a,\tau}:=\{Z=(x,y):r>0,\ j_a(Z)<\tau r\}$ for $\tau>0$.
		There are uniform constants $\varepsilon,C>0$ such that signed distance
		$l_a$ to $\barC_a$ is smooth on $\mathcal T_{a,2\varepsilon}$.
		For any smooth properly embedded boundaryless hypersurface $N$ with
		this area-ratio bound, put $\mathbf n_a:=D l_a$,
		$\vartheta_a^2:=1-(\mathbf n_N\cdot\mathbf n_a)^2$, and $w:=r^2/(1+r^2)$.
		Then
		\begin{equation}\label{eq:cyl-weighted-caccioppoli}
			\int_{N\cap\mathcal T_{a,\varepsilon}}
			w\vartheta_a^2\rho\dd\mu_{N}
			\le C\bigl(D_a(N)^2+G(N)^2+|q(a)|^2\bigr).
		\end{equation}
	\end{lem}
	\begin{proof}
		The Gaussian moments are uniformly bounded by the
		area-ratio assumption. We may assume that $G(N)$ is finite.
		
		Write $\mathbf n:=\mathbf n_N$. By the signed-distance identities, we have
		\begin{equation}\label{eq:cyl-tilt-distance-identities}
			|l_a|=j_a,\qquad Z\cdot\mathbf n_a=l_a,\qquad
			D\mathbf n_a\,\mathbf n_a=0,\qquad |\nabla_N l_a|^2=\vartheta_a^2.
		\end{equation}
		At the normal projection of $Z$, let $\kappa_j$ be the eigenvalues
		of $D\mathbf n_a$ restricted to the tangent space. The base radius
		is comparable to $r$, and the parallel-surface formula gives
		\begin{equation}\label{eq:cyl-tilt-parallel-curvature}
			\begin{aligned}
				|D\mathbf n_a|&\le C\max_j\frac{|\kappa_j|}{|1+l_a\kappa_j|}\le C/r,\\
				|\operatorname{div}\mathbf n_a|
				&=\left|\sum_j\frac{\kappa_j}{1+l_a\kappa_j}\right|
				\le\left|\sum_j\kappa_j\right|
				+|l_a|\sum_j\frac{\kappa_j^2}{|1+l_a\kappa_j|}
				\le C|q(a)|/r+C|l_a|/r^2.
			\end{aligned}
		\end{equation}
		Here $|l_a\kappa_j|\le1/2$ after decreasing $\varepsilon$;
		the last inequalities use \eqref{eq:reduced-cone-curvature}.
		By \eqref{eq:cyl-tilt-distance-identities} and
		\eqref{eq:cyl-tilt-parallel-curvature},
		\begin{equation}\label{eq:cyl-tilt-tubular-divergence}
			\begin{aligned}
				\operatorname{div}_N(l_a\mathbf n_a)
				&=\vartheta_a^2+l_a\bigl(\operatorname{div}\mathbf n_a
				-\langle D\mathbf n_a\,\mathbf n,\mathbf n\rangle\bigr)\\
				&\ge\vartheta_a^2-|l_a|\bigl(C|q(a)|/r+C|l_a|/r^2+C\vartheta_a^2/r\bigr)\\
				&\ge(1-C|l_a|/r)\vartheta_a^2-Cl_a^2/r^2-C|q(a)||l_a|/r.
			\end{aligned}
		\end{equation}
		In the first inequality, $D\mathbf n_a\,\mathbf n_a=0$ gives
		$|\langle D\mathbf n_a\,\mathbf n,\mathbf n\rangle|
		\le |D\mathbf n_a|\vartheta_a^2$. This proves
		\eqref{eq:cyl-tilt-tubular-divergence}.
		
		Choose a smooth cutoff $\zeta$ of $l_a/r$ equal to one for
		$|l_a|/r<\varepsilon$ and zero for $|l_a|/r\ge2\varepsilon$, with
		$\varepsilon$ fixed so that $2C\varepsilon<1/2$ in
		\eqref{eq:cyl-tilt-tubular-divergence}. Extend tests containing
		$\zeta$ by zero outside $\mathcal T_{a,2\varepsilon}$. Put
		$w:=r^2/(1+r^2)$. For the field $w\zeta^2l_a\mathbf n_a$,
		Young's inequality and the Cauchy--Schwarz inequality give
		\begin{equation}\label{eq:cyl-tilt-absorption-terms}
			\begin{aligned}
				|l_a\zeta^2\nabla_Nw\cdot\nabla_N l_a|
				&\le\epsilon w\zeta^2\vartheta_a^2
				+C_\epsilon l_a^2\frac{|Dw|^2}{w}
				\le\epsilon w\zeta^2\vartheta_a^2+C_\epsilon l_a^2,\\
				2|w\zeta l_a\nabla_N\zeta\cdot\nabla_N l_a|
				&\le Cw\mathbf1_{\{\varepsilon r\le |l_a|\le2\varepsilon r\}}
				\le Cl_a^2,\\
				wl_a^2/r^2&\le l_a^2,\qquad
				\tfrac12w\zeta^2l_a\,Z\cdot\mathbf n_a
				=\tfrac12w\zeta^2l_a^2\le\tfrac12l_a^2,\\
				w|q(a)|\frac{|l_a|}{r}
				&\le\tfrac12\frac{w}{r^2}l_a^2+\tfrac12w|q(a)|^2
				\le C(l_a^2+|q(a)|^2),\\
				\left|\int_Nw\zeta^2l_a\,\Phi\cdot\mathbf n_a\rho\dd\mu_{N}\right|
				&\le\left(\int_Nw^2\zeta^4l_a^2\rho\dd\mu_{N}\right)^{1/2}G(N)\\
				&\le D_a(N)G(N)\le\tfrac12(D_a(N)^2+G(N)^2).
			\end{aligned}
		\end{equation}
		The first line uses $|Dw|^2/w=4(1+r^2)^{-3}\le4$;
		the second uses $|D\zeta|\le C/r$ on its support.
		Gaussian first variation, \eqref{eq:cyl-tilt-tubular-divergence},
		and \eqref{eq:cyl-tilt-absorption-terms} yield
		\[
		(\tfrac12-\epsilon)\int_Nw\zeta^2\vartheta_a^2\rho\dd\mu_{N}
		\le C_\epsilon\bigl(D_a(N)^2+G(N)^2+|q(a)|^2\bigr).
		\]
		For $\epsilon<1/4$, $\zeta=1$ on $\mathcal T_{a,\varepsilon}$ gives
		\eqref{eq:cyl-weighted-caccioppoli}.
		To justify first variation, insert cutoffs at $r=\delta$ and $|Z|=R$.
		On $\delta<r<2\delta$, the field has size $O(r^3)$, so its
		cutoff error is bounded by $C\delta^2\int_N\rho\dd\mu_N\to0$.
		The errors at infinity tend to zero by \eqref{eq:gaussian-moment-bound}.
		
	\end{proof}
	
	The following estimate is the cylindrical analogue of \cite[Proposition~2.13]{FH1}. It allows a nonminimal reduced reference and retains a positive axial-normal term; accordingly, the conclusion is a one-sided area bound.
	\begin{prop}\label{prop:anchorarea}
		Fix the all-center area-ratio bound $\Lambda$.
		There are constants $0<r_1<r_2$ and $L_0>0$, a compact set
		\[
		\mathcal K:=\{(x,y)\in\R^{m+1}\times\R^k:
		r_1\le |x|\le r_2,\ |y|\le L_0\},
		\]
		a relatively compact open set $\mathcal U$ disjoint from
		the axis with $\mathcal K\Subset\mathcal U$, and constants
		$\epsilon,c,C>0$ with the following property. Suppose that $N$ is a
		smooth, properly embedded, boundaryless hypersurface satisfying this
		area-ratio bound and containing a normal graph of a function
		$u$ over an open set $\Omega_a\Subset\Reg(\barC_a)$, with
		$\|u\|_{C^1(\Omega_a)}\le\epsilon$. Suppose that
		$\overline{\mathcal U}\cap\barC_a\Subset\Omega_a$, with a uniform
		positive distance to $\partial\Omega_a$, and that
		$N\cap\mathcal U=\operatorname{graph}(u)\cap\mathcal U$.
		Then
		\begin{equation}\label{eq:anchorarea-statement}
			F(N)-F(\barC_a)+c\int_N|\mathbf n_y|^2\rho\dd\mu_{N}
			\le C\bigl(D_a(N)^2+G(N)^2+|q(a)|^2\bigr).
		\end{equation}
		The constants are uniform for the reduced references. The same
		statement holds for $\barC_z$ after aligning its axis by the fixed
		rotation $R_{\mathsf B}^{-1}$.
	\end{prop}
	\begin{proof}
		\noindent\textbf{Step 1. The radial fields and the exact identity.}
		Choose nonnegative smooth functions $g_x(r)$ and $g_y(r_y)$, with $g_x$
		supported in a fixed annulus $r_1<r<r_2$ and $g_y$ supported in $[0,L_0]$.
		More precisely, choose $g_y(r_y):=\widehat g_y(r_y^2)$ with
		$\widehat g_y\in C_c^\infty(\R)$ nonnegative and nonzero on
		$[0,L_0^2)$, and vanishing for arguments at least $L_0^2$.
		Then $y\mapsto g_y(|y|)=\widehat g_y(|y|^2)$ is smooth also at
		$y=0$. Multiplying $g_x$ and $g_y$ by positive constants, normalize them by
		\[
		\int_0^\infty(1-g_x(r))r^{m-1}e^{-r^2/4}\dd r=0,\qquad
		\int_0^\infty(1-g_y(r_y))r_y^{k-1}e^{-r_y^2/4}\dd r_y=0.
		\]
		Define
		\[
		f_x(r):=r^{-m}e^{r^2/4}\int_0^r(1-g_x(t))t^{m-1}e^{-t^2/4}\dd t,
		\quad
		f_y(r_y):=r_y^{-k}e^{r_y^2/4}\int_0^{r_y}(1-g_y(t))t^{k-1}e^{-t^2/4}\dd t.
		\]
		They solve
		\[
		mf_x+rf_x'-r^2f_x/2=1-g_x,\qquad
		kf_y+r_y f_y'-r_y^2f_y/2=1-g_y.
		\]
		Their extensions at zero are smooth radial functions; in particular
		$f_x(0)=1/m$ and $f_x'/r$ is bounded there. The zero-mean normalizations
		give $f_x=-2r^{-2}+O(r^{-4})$ and $f_y=-2r_y^{-2}+O(r_y^{-4})$ at infinity,
		with corresponding derivative bounds. In particular
		$|f_x|\le C/(1+r^2)$, $|f_x'/r|\le C/(1+r^2)^2$,
		and $r_y|f_y|$, $r_y|f_y'|$, and $f_y'/r_y$ are bounded.
		
		Put $\mathcal X:=f_x(r)x+g_x(r)f_y(r_y)y$. This is a bounded smooth ambient field.
		For a local choice of unit normal $\mathbf n:=(\mathbf n_x,\mathbf n_y)$, set
		\[
		\mathfrak a:=x\cdot\mathbf n_x,\qquad \mathfrak b:=y\cdot\mathbf n_y,\qquad \mathfrak t:=|\mathbf n_y|^2.
		\]
		For an ambient vector field $\mathcal X$, we use the convention
		\begin{equation}
			\operatorname{div}_{N,\rho}\mathcal X
			:=\operatorname{div}_N\mathcal X-\tfrac12\langle Z,\mathcal X\rangle.
		\end{equation}
		For the field $\mathcal X$, direct differentiation gives
		\begin{equation}\label{eq:area-product-divergence}
			\operatorname{div}_{N,\rho}\mathcal X=1-g_xg_y+Q_{\rm area},
		\end{equation}
		where
		\begin{equation}\label{lat:Q}
			Q_{\rm area}:=(f_x-g_xf_y)\mathfrak t-\frac{f_x'}r\mathfrak a^2
			-g_x\frac{f_y'}r_y \mathfrak b^2-\frac{g_x'}r f_y\mathfrak a\mathfrak b.
		\end{equation}
		Indeed,
		\[
		\begin{aligned}
			\operatorname{div}_{N,\rho}(f_xx)
			&=1-g_x+f_x\mathfrak t-(f_x'/r)\mathfrak a^2,\\
			\operatorname{div}_{N,\rho}(g_xf_yy)
			&=g_x(1-g_y)-g_xf_y\mathfrak t
			-g_x(f_y'/r_y)\mathfrak b^2-(g_x'/r)f_y\mathfrak a\mathfrak b,\\
			\operatorname{div}_{N,\rho}\mathcal X
			&=1-g_x+g_x(1-g_y)+Q_{\rm area}=1-g_xg_y+Q_{\rm area}.
		\end{aligned}
		\]
		
		Let $\chi_R$ be an ambient cutoff equal to one on $B_R$,
		zero outside $B_{2R}$, and satisfying $|D\chi_R|\le C/R$.
		Gaussian first variation with the compactly supported field $\chi_R\mathcal X$ and $\operatorname{div}_N\mathcal X
		=\operatorname{div}_N(\mathcal X^\top)-\langle\mathbf H_N,\mathcal X\rangle$ give
		\begin{equation}\label{eq:area-cutoff-first-variation}
			\int_N\chi_R\operatorname{div}_{N,\rho}\mathcal X\,\rho\dd\mu_N
			=-\int_N\chi_R\Phi\cdot \mathcal X\,\rho\dd\mu_N
			-\int_N\nabla_N\chi_R\cdot \mathcal X^\top\rho\dd\mu_N.
		\end{equation}
		Since $\mathcal X$ is bounded, the last term is bounded in absolute value by
		\[
		\frac{C}{R}\int_{N\setminus B_R}\rho\dd\mu_N\longrightarrow0
		\]
		by \eqref{eq:gaussian-moment-bound}. The source is integrable by
		the Cauchy--Schwarz inequality, since $G(N)<\infty$ and $F(N)<\infty$.
		The explicit coefficients in \eqref{lat:Q} also give $|Q_{\rm area}|\le C$.
		Passing to the limit in \eqref{eq:area-cutoff-first-variation} by
		dominated convergence and substituting \eqref{eq:area-product-divergence}
		and \eqref{lat:Q}, we obtain
		\[
		F(N)-\int_N g_xg_y\rho\dd\mu_N
		+\int_N Q_{\rm area}\rho\dd\mu_N
		=-\int_N\Phi\cdot \mathcal X\rho\dd\mu_N.
		\]
		By rearranging, we get
		\begin{equation}\label{lat:exact}
			F(N)-F(\barC_a)
			=\mathcal K_a-\int_N Q_{\rm area}\rho\dd\mu_N-\int_N\Phi\cdot \mathcal X\rho\dd\mu_N,
		\end{equation}
		where
		\begin{align*}
			\mathcal K_a
			:=\int_Ng_xg_y\rho\dd\mu_N-\int_{\barC_a}g_xg_y\rho\dd\mu_{\barC_a}.
		\end{align*}
		To identify the reference integral, let $\Sigma_a$ be the link of $\C_a$.
		Polar integration and the two zero-mean normalizations give
		\[
		\begin{aligned}
			\int_{\barC_a}g_xg_y\rho\dd\mu_{\barC_a}
			&=(4\pi)^{-n/2}\mathcal H^{m-1}(\Sigma_a)
			\mathcal H^{k-1}(S^{k-1})\\
			&\quad\cdot\left(\int_0^\infty g_x(r)r^{m-1}e^{-r^2/4}\dd r\right)
			\left(\int_0^\infty g_y(r_y)r_y^{k-1}e^{-r_y^2/4}\dd r_y\right)\\
			&=F(\barC_a).
		\end{aligned}
		\]
		
		In a fixed narrow normal tube of $\barC_a$, let $l_a$ denote signed
		Euclidean distance to $\barC_a$, set $\mathbf n_a:=D l_a$, and put
		$\vartheta_a^2:=1-(\mathbf n\cdot\mathbf n_a)^2$.
		We use the \textbf{radial weight}
		\[
		w(r):=\frac{r^2}{1+r^2},\qquad r:=|x|.
		\]
		It satisfies $0\le w\le1$, $w(r)\asymp r^2$ for $0<r\le1$,
		and $w(r)\asymp1$ for $r\ge1$. Thus the integral of
		$w\vartheta_a^2$ controls the normal-angle error with a quadratic
		weight near the axis and a uniformly positive weight on each
		fixed region separated from the axis.
		
		\noindent\textbf{Step 2. The compact graphical contribution.}
		The support of $g_x(r)g_y(r_y)$ is contained in the compact set $\mathcal K$ from the statement. Assume that on an ambient neighborhood of
		$\mathcal K$, $N$ is precisely one sufficiently small $C^1$ normal
		graph over $\Omega_a$. Choose
		$\Omega_a^0\Subset\Omega_a$ containing
		$\overline{\mathcal U}\cap\barC_a$, with the distance from
		$\Omega_a^0$ to $\partial\Omega_a$ uniformly positive.
		
		For a base point $Z\in\Omega_a$, the graph point is $Z+u(Z)\mathbf n_a(Z)$. On this graph, $|x(Z+u\mathbf n_a)|^2=|x(Z)|^2+u^2$, the axial coordinate
		is unchanged, and $\rho(Z+u\mathbf n_a)=\rho(Z)e^{-u^2/4}$.
		Thus the radial cutoffs and the Gaussian weight have no term linear in $u$.
		On $\Omega_a^0$ the graph area factor and the Gaussian satisfy
		\[
		J_u=1-H_{\barC_a}u+O(u^2+|\nabla u|^2),\qquad
		e^{-u^2/4}=1+O(u^2),\qquad |H_{\barC_a}|\le Cb.
		\]
		All constants are uniform because the domain is separated from the
		axis. The compact cutoffs can be pulled back to $\Omega_a^0$; their changes are $O(u^2)$ because the transverse normal
		is perpendicular to the base radial direction. Hence we obtain
		\begin{equation}\label{lat:compact}
			\begin{aligned}
				|\mathcal K_a|&\le C\int_{\Omega_a^0}
				(u^2+|\nabla u|^2+b|u|)\rho\dd\mu_{\barC_a}\\
				&\le C\int_{\Omega_a^0}u^2\rho\dd\mu_{\barC_a}
				+C\int_{\operatorname{graph}(u|_{\Omega_a^0})}w\vartheta_a^2\rho\dd\mu_N
				+Cb^2\int_{\Omega_a^0}\rho\dd\mu_{\barC_a}\\
				&\le C(D_a(N)^2+G(N)^2+b^2).
			\end{aligned}
		\end{equation}
		The second line uses Young's inequality for $b|u|$, graph measure
		comparability, and $|\nabla u|\le C\vartheta_a$ on a sufficiently small
		graph. For the last line, apply the weighted normal-angle estimate
		\eqref{eq:cyl-weighted-caccioppoli} of Lemma~\ref{lem:weighted-normal-angle}:
		\[
		\int_{N\cap\{|l_a|<\varepsilon r\}}w\vartheta_a^2\rho\dd\mu_N
		\le C\bigl(D_a(N)^2+G(N)^2+|q(a)|^2\bigr)
		=C(D_a(N)^2+G(N)^2+b^2).
		\]
		The compact graph lies in this tube when its $C^1$ norm is
		sufficiently small, and $w\ge c>0$ on the fixed domain.
		
		\noindent\textbf{Step 3. The near-axis term and absorption.}
		Choose a fixed $r_c<r_1$ so small that $f_x\ge c_0>0$ on $r<r_c$.
		There $g_x=0$, so \eqref{lat:Q} gives
		\[
		Q_{\rm area}=f_x\mathfrak t-(f_x'/r)\mathfrak a^2\ge c_0\mathfrak t-C\mathfrak a^2.
		\]
		Inside the normal tube, we have
		\[
		|\mathfrak a|\le |l_a|+Cr\vartheta_a,\qquad |\mathfrak b|\le r_y\vartheta_a.
		\]
		Hence the negative $\mathfrak a^2$ term integrates to at most
		$C(D_a(N)^2+G(N)^2+b^2)$ on $r<r_c$, by \eqref{eq:cyl-weighted-caccioppoli}.
		Off the tube, $j_a\ge cr$, so $\mathfrak a^2\le r^2\le Cj_a^2$.
		For $r\ge r_c$, the explicit coefficients in \eqref{lat:Q} give
		$|Q_{\rm area}|\le C(l_a^2+\vartheta_a^2)$ inside the tube. Since $w\ge c(r_c)>0$
		there, \eqref{eq:cyl-weighted-caccioppoli} controls its integral. Outside the tube,
		$Q_{\rm area}$ is bounded and $j_a\ge cr_c$. We conclude
		\begin{equation}\label{lat:Q-lower}
			\int_N Q_{\rm area}\rho\dd\mu_N\ge c_0\int_{N\cap\{r<r_c\}}\mathfrak t\rho\dd\mu_N-C(D_a(N)^2+G(N)^2+b^2).
		\end{equation}
		
		We next estimate the source term. Near the axis,
		$|\mathcal X^\perp|^2=f_x^2\mathfrak a^2\le C \mathfrak a^2$. Away from the axis its explicit
		formula gives $|\mathcal X^\perp|^2\le C(l_a^2+\vartheta_a^2)$ in the tube.
		Off the tube, $|\mathcal X|^2\le Cr^2$ near zero and is bounded elsewhere,
		so it is bounded by $Cj_a^2$. Therefore
		\begin{equation}\label{eq:area-source-bound}
			\begin{aligned}
				\int_N|\mathcal X^\perp|^2\rho\dd\mu_N
				&\le C\bigl(D_a(N)^2+G(N)^2+b^2\bigr),\\
				\left|\int_N\Phi\cdot\mathcal X\rho\dd\mu_N\right|
				&\le G(N)\left(\int_N|\mathcal X^\perp|^2\rho\dd\mu_N\right)^{1/2}\\
				&\le CG(N)\sqrt{D_a(N)^2+G(N)^2+b^2}
				\le C\bigl(D_a(N)^2+G(N)^2+b^2\bigr).
			\end{aligned}
		\end{equation}
		By combining \eqref{eq:area-source-bound} with \eqref{lat:exact},
		\eqref{lat:compact}, and \eqref{lat:Q-lower}, we obtain
		\[
		F(N)-F(\barC_a)+c_0\int_{N\cap\{r<r_c\}}|\mathbf n_y|^2\rho\dd\mu_N
		\le C(D_a(N)^2+G(N)^2+b^2).
		\]
		For $r\ge r_c$, the cone normal has no axial component, so
		$|\mathbf n_y|^2\le\vartheta_a^2$ in the tube, where $w\ge r_c^2/(1+r_c^2)$.
		Outside the tube, $j_a\ge\varepsilon r_c$. Thus
		\[
		\int_{N\cap\{r\ge r_c\}}|\mathbf n_y|^2\rho\dd\mu_N
		\le C\int_{N\cap\mathcal T_{a,\varepsilon}}w\vartheta_a^2\rho\dd\mu_N
		+(\varepsilon r_c)^{-2}D_a(N)^2
		\le C(D_a(N)^2+G(N)^2+b^2),
		\]
		where the last inequality is \eqref{eq:cyl-weighted-caccioppoli}.
		We add a sufficiently small positive multiple of this estimate to obtain
		\begin{equation}
			F(N)-F(\barC_a)+c\int_N|\mathbf n_y|^2\rho\dd\mu_{N}
			\le C(D_a(N)^2+G(N)^2+b^2).
		\end{equation}
		Since $b=|q(a)|$, this is \eqref{eq:anchorarea-statement}.
	\end{proof}
	
	\section{Compactness and weighted stationary rigidity}\label{sec:compactness}
	In this section, we first prove the simultaneous scale and center selection in
	Theorem~\ref{thm:measure-selection}, and use it to establish the
	transverse tube estimate \eqref{eq:cyl-relative-tube} in
	Lemma~\ref{lem:cyl-tube}. We first prove a stationary tube estimate with exponent below two
	and apply it to normalized dilations of a stationary varifold at infinity.
	Proposition~\ref{prop:endpoint-static-growth} then allows an exponent
	$\sigma>2$ in the flow estimate. This permits the dyadic summation in
	\eqref{eq:cyl-relative-hardy}, which controls the local inverse-radius term
	in \eqref{eq:cyl-tilt-local-to-global}.
	The product first-variation identity \eqref{eq:cyl-tilt-product-divergence}
	then gives the axial-tilt estimate of
	Lemma~\ref{lem:axial-tilt-estimate}. Finally,
	Lemma~\ref{lem:cyl-gaussian-tightness} yields the moment and tail bounds
	\eqref{eq:cyl-relative-tail}, which prevent loss of normalized
	Gaussian squared distance at spatial infinity.
	
	\subsection{The transverse measure selection argument}
	\label{gsel:section}
\begin{figure}[p]
		\centering
		\begin{tikzpicture}[
			font=\small,
			box/.style={draw=linkblue!65!black,fill=linkblue!3,
				rounded corners=2pt,line width=.6pt,align=center,
				text width=13.6cm,inner sep=7pt},
			result/.style={box,text width=6.25cm,minimum height=3.15cm,
				fill=citegreen!4,draw=citegreen!70!black},
			arrow/.style={-{Latex[length=2.5mm]},line width=.8pt,draw=linkblue!80!black}
			]
			\node[box] (input) {
				\textbf{Mass and error measures}\\[3pt]
				$\mu_{i,t}$: mass measure\qquad $\eta_{i,t}$: error measure\\[3pt]
				$\Delta_i/\mathsf W_i\to0$; small scales in weighted time measure
				\eqref{eq:ms-radius-probability}
			};
			\node[box,below=7mm of input] (scale) {
				\textbf{1. Choose a transverse scale at each time}\\[3pt]
				$\mathfrak r_i(t)$ nearly maximizes
				$\mathfrak r^{-\sigma}\mu_{i,t}(T_{\mathfrak r})$\\[3pt]
				$\displaystyle
				\mu_{i,t}(T_s)\le C\bigl(s/\mathfrak r_i(t)\bigr)^\sigma
				\mathsf M_i(t)$\qquad\eqref{eq:ms-1}--\eqref{eq:ms-2}
			};
			\node[box,below=7mm of scale] (time) {
				\textbf{2. Choose a time with small normalized error}\\[3pt]
				$\displaystyle \dd\pi_i(t):=\frac{Q_i(t)}{\mathsf W_i}\dd t,
				\qquad \int_I b_i(t)\dd \pi_i(t)\le C\frac{\Delta_i}{\mathsf W_i}\to0$\\[3pt]
				Choose $t_i$ with $\mathfrak r_i(t_i)\to0$ and $b_i(t_i)\to0$.
				\quad\eqref{eq:ms-16}
			};
			\node[box,below=7mm of time] (center) {
				\textbf{3. Choose an axial center with positive $\mu_{i,t_i}$-mass}\\[3pt]
				Choose $Y_i\in\Omega$ with $\mathfrak f(Y_i)>0$,
				$S(Y_i)\le C\mathfrak f(Y_i)$, and
				$\mathscr B(Y_i)\le o(1)\mathfrak f(Y_i)$.\quad\eqref{eq:ms-21}\\[4pt]
				$\ell_i:=\mathfrak r_i(t_i)h(Y_i)$,\qquad
				$\mathsf P_i:=\mu_{i,t_i}(\mathscr Q_{2,1}(Y_i))>0$
			};
			\node[result,below=10mm of center.south west,anchor=north west] (position) {
				\textbf{Mass bounds}\\[5pt]
				$\displaystyle
				\frac{\mu_{i,t_i}(\mathscr Q_{\delta,R}(Y_i))}{\mathsf P_i}
				\le C\delta^{\sigma-\varepsilon}R^{k+\varepsilon}$\\[5pt]
				$\displaystyle
				\frac{\mu_{i,t_i}(\mathscr Q_{R,R}(Y_i))}{\mathsf P_i}
				\le CR^{k+\sigma+\varepsilon}$\\[4pt]
				\eqref{eq:ms-17}--\eqref{eq:ms-18}
			};
			\node[result,below=10mm of center.south east,anchor=north east] (error) {
				\textbf{Vanishing normalized error}\\[5pt]
				$\displaystyle
				\frac{\ell_i^4}{\mathsf P_i}\,
				\eta_{i,t_i}(\mathscr Q_{R,R}(Y_i))
				\le\epsilon_iR^{k+3}$\\[5pt]
				$\epsilon_i\to0$ for every fixed $R$.\\[4pt]
				Slice estimate: \eqref{eq:ms-19}\\[2pt]
				Spacetime estimate: \eqref{eq:ms-20}
			};
			\coordinate (join) at ($(position.south)!0.5!(error.south)$);
			\node[box,below=11mm of join] (use) {
				\textbf{Later application: rescaling about $(0,Y_i,t_i)$}\\[3pt]
				$Z\mapsto\ell_i^{-1}(Z-(0,Y_i))$,\qquad
				$t\mapsto\ell_i^{-2}(t-t_i)$\\[3pt]
				Positive normalized mass and vanishing normalized error\\[3pt]
				$\displaystyle
				\frac{\operatorname{dist}(Y_i,\partial\Omega)}{\ell_i}
				\ge\frac1{\mathfrak r_i(t_i)}\to\infty$\\[3pt]
				Used in Lemma~\ref{lem:cyl-tube}, with \eqref{htps:rectangular-bound}.
			};
			\draw[arrow] (input)--(scale);
			\draw[arrow] (scale)--(time);
			\draw[arrow] (time)--(center);
			\draw[arrow] (center.south -| position.north)--(position.north);
			\draw[arrow] (center.south -| error.north)--(error.north);
			\draw[arrow] (position.south)--(position.south |- use.north);
			\draw[arrow] (error.south)--(error.south |- use.north);
		\end{tikzpicture}
		\caption{The simultaneous choice in Theorem~\ref{thm:measure-selection}.
			Here $\mathscr Q_{s,L}(Y_i)=\{|x|<s\ell_i,\ |y-Y_i|<L\ell_i\}$.}
		\label{fig:measure-selection}
	\end{figure}
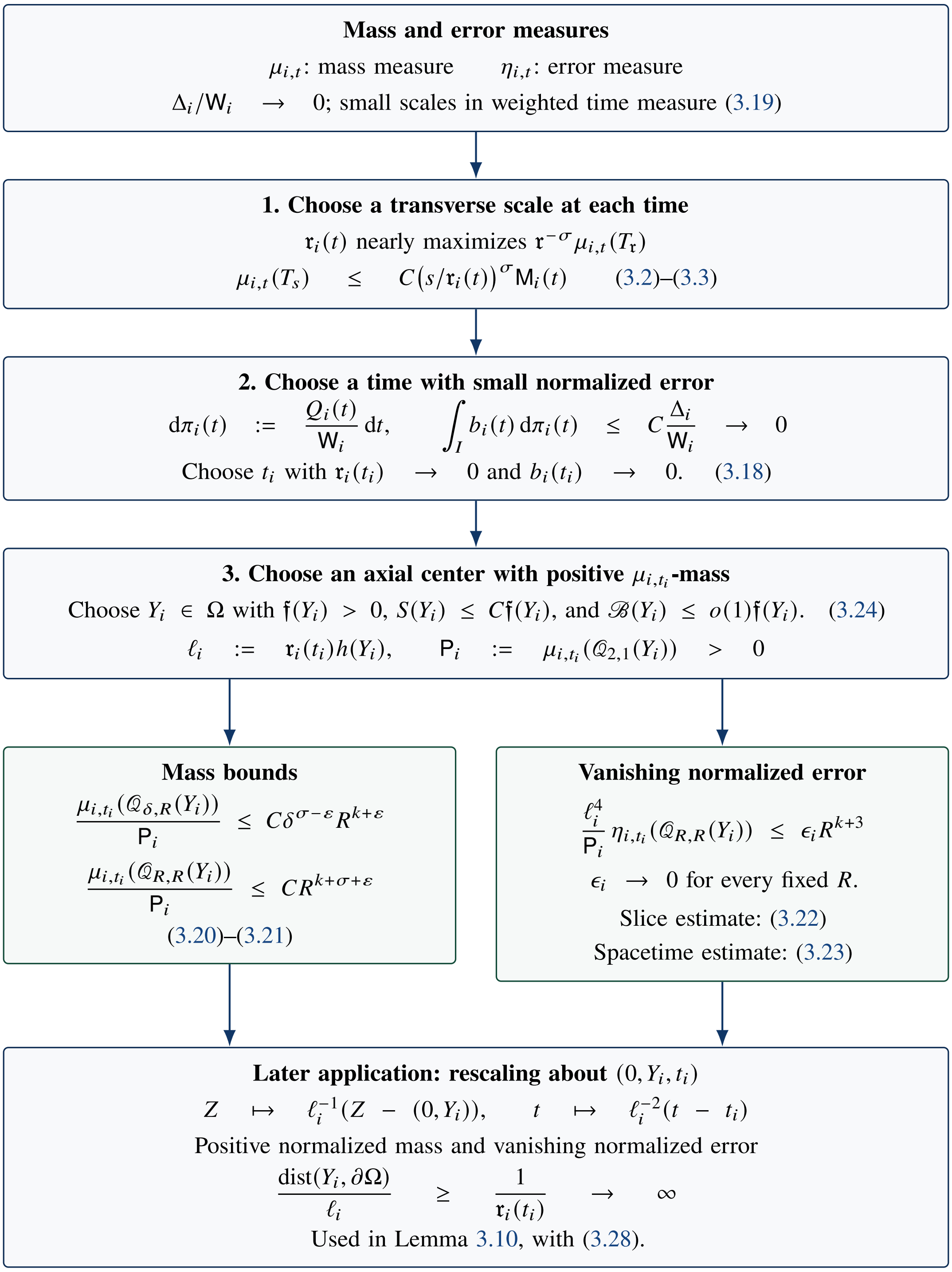

	In this subsection, we aim to choose a time, a transverse scale,
	and an axial center for which the squared-distance integral is positive and
	the integral of $|\Phi_i|^2$ is small relative to the squared-distance integral.
	The main conclusion is Theorem~\ref{thm:measure-selection}: the same
	choice satisfies the squared-distance bounds \eqref{eq:ms-17}--\eqref{eq:ms-18}
	and the error estimates \eqref{eq:ms-19}--\eqref{eq:ms-20}.
	Lemma~\ref{lem:small-selected-radii} supplies a sufficient condition
	for the chosen scales to tend to zero. These results, together with
	the two-parameter refinement \eqref{htps:rectangular-bound}, are used
	in Step~1 of the proof of Lemma~\ref{lem:cyl-tube} to obtain the
	transverse tube bound \eqref{eq:cyl-relative-tube}.
	Figure~\ref{fig:measure-selection} summarizes the choices and their estimates.

	This subsection is purely measure-theoretic. Fix transverse cone dimension
	$m>2$ and axial dimension $k\ge1$. It supplies the selection used in Lemma~\ref{lem:cyl-tube}. The abstract selection estimates below apply to $0<\sigma<4$.
	In the stationary argument below the exponent is between
	$\tau_0=\max\{0,2-2\nu_0\}$ and $2$. In the flow argument it is
	between $2$ and $\min\{4,2+2\nu_0\}$. The measure selection is the
	same for the two ranges.
	
	Let $\Omega\subset\R^k$ be bounded and open, and put
	\[
	h(y):=\min\{1,\operatorname{dist}(y,\R^k\setminus\Omega)\}.
	\]
	Thus $h>0$ on $\Omega$, $h\leq1$, and $h$ is $1$-Lipschitz.  Write a point in
	$\R^{m+1}\times\R^k$ as $(x,y)$.  For $0<{\mathfrak r}\leq {\mathfrak r}_0$, define the \textbf{variable transverse tube}
	\[
	T_{\mathfrak r}:=\{(x,y):y\in\Omega,\ |x|<{\mathfrak r} h(y)\}.
	\]
	Fix an interval $I$, a slightly larger interval
	$J\supset I+[-c_0^2,c_0^2]$, and a small constant $0<c_0<1/100$.
	We take $0<{\mathfrak r}_0<c_0/100$.  The supremum defining $Q_i$ below is taken over
	$\{2^{-j}{\mathfrak r}_0:j\in\mathbb Z_{\ge0}\}$.
	
	For each index $i$ and time $t\in J$, let $\mu_{i,t}$ and
	$\eta_{i,t}$ be nonnegative Radon measures on $\R^{m+1}\times\Omega$.
	We assume that these are \textbf{measurable families of measures}:
	for every Borel set $E\subset\R^{m+1}\times\Omega$, the functions
	$t\mapsto\mu_{i,t}(E)$ and $t\mapsto\eta_{i,t}(E)$ are Lebesgue
	measurable on $J$. Set
	\[
	\mathcal D_i(t):=\eta_{i,t}(\{|x|<1\}\times\Omega),\qquad
	\Delta_i:=\int_J \mathcal D_i(t)\dd t<\infty.
	\]
	Extend $\eta_{i,t}$ by zero outside $J$.
	For a smooth rescaled flow $N_t^i$, with shrinker quantity
	$\Phi_i$ as in \eqref{eq:flow}, the intended application is
	\[
	\dd\mu_{i,t}:=\operatorname{dist}((x,y),\barC)^2\dd\mu_{N_t^i},
	\qquad \dd\eta_{i,t}:=|\Phi_i|^2\dd\mu_{N_t^i},
	\quad \barC:=\C\times\R^k.
	\]
	
	Fix $\sigma>0$, and define
	\begin{equation}\label{eq:ms-maximal-mass}
		Q_i(t):=\sup_{\mathfrak r\in\{2^{-j}\mathfrak r_0:j\in\mathbb Z_{\ge0}\}}{\mathfrak r}^{-\sigma}\mu_{i,t}(T_{\mathfrak r}),
		\qquad \mathsf W_i:=\int_I Q_i(t)\dd t.
	\end{equation}
	Assume $0<\mathsf W_i<\infty$. Redefine $Q_i$ to be zero on the
	null subset of $I$ where it is infinite. On $\{0<Q_i(t)<\infty\}$
	choose measurably a \textbf{dyadic near maximizer}
	${\mathfrak r}_i(t)$, and abbreviate
	\begin{equation}
		{\mathfrak r}:={\mathfrak r}_i(t),\qquad \mathsf M_i(t):=\mu_{i,t}(T_{\mathfrak r}),\qquad
		\tfrac12 {\mathfrak r}^\sigma Q_i(t)\leq \mathsf M_i(t)\leq {\mathfrak r}^\sigma Q_i(t).
		\label{eq:ms-1}
	\end{equation}
	Measurability follows by taking the first radius in a countable enumeration
	which satisfies the near-maximality inequality.  In particular, whenever
	$s\leq {\mathfrak r}_0$, we have
	\begin{equation}
		\mu_{i,t}(T_s)\leq C_\sigma(s/{\mathfrak r})^\sigma \mathsf M_i(t) . \label{eq:ms-2}
	\end{equation}
	For every $0<s\le\mathfrak r_0$, choose
	$\widehat s\in\{2^{-j}\mathfrak r_0:j\in\mathbb Z_{\ge0}\}$ with
	$s\le\widehat s<2s$. Then we have
	\begin{equation}\label{eq:ms-arbitrary-radius}
		Q_i(t)\ge\widehat s^{-\sigma}\mu_{i,t}(T_{\widehat s})
		\ge2^{-\sigma}s^{-\sigma}\mu_{i,t}(T_s).
	\end{equation}
	Together with \eqref{eq:ms-1}, this proves \eqref{eq:ms-2} also for non-dyadic $s$.
	We call a time $t$ an \textbf{admissible time} if $Q_i(t)>0$ and $3{\mathfrak r}_i(t)\leq {\mathfrak r}_0$.
	Set the nonnegative averaging functions below equal to zero at nonadmissible times.
	The positive lower bounds apply only at admissible times.
	When selected radii tend to zero in $Q_i\dd t/\mathsf W_i$ probability, admissible
	times have probability tending to one.
	
	Fix an index $i$ and an admissible time $t$. The functions below
	also depend on $i,t$; we suppress these arguments while retaining
	them in $\mathsf M_i(t)$. For $Y\in\Omega$, set $\ell(Y):={\mathfrak r} h(Y)$, and use the cylinders
	\[
	\mathscr Q_{s,L}(Y)
	:=\{(x,y): |x|<s \ell(Y),\ |y-Y|<L \ell(Y)\}.
	\]
	The \textbf{normalization function} is
	\begin{equation}
		\mathfrak f(Y):=\frac{\mu_{i,t}(\mathscr Q_{2,1}(Y))}{\mathsf M_i(t)\,\ell(Y)^k}.
	\end{equation}
	The factor $2$ in the transverse radius is essential.  The same expression
	with transverse radius $\ell(Y)$ need not have an averaged lower bound, because
	the original measure may lie arbitrarily close to the boundary of $T_{\mathfrak r}$.
	
	For $0<\delta\leq1$ and $R\geq1$ satisfying $R{\mathfrak r}\leq c_0$, put
	\begin{equation}
		F_{\delta,R}(Y)
		:=\frac{\mu_{i,t}(\mathscr Q_{\delta,R}(Y))}{\mathsf M_i(t)\,\ell(Y)^k},
		\qquad
		H_R(Y)
		:=\frac{\mu_{i,t}(\mathscr Q_{R,R}(Y))}{\mathsf M_i(t)\,\ell(Y)^k}.
	\end{equation}
	All the following constants are independent of $i,t,{\mathfrak r},\mathsf M_i(t)$.
	
	\begin{lem}\label{lem:measure-averaging}
		At every admissible time,
		\begin{align}
			c\leq\int_\Omega \mathfrak f(Y)\dd Y&\leq C, \label{eq:ms-5}\\
			\int_\Omega F_{\delta,R}(Y)\dd Y&\leq C\delta^\sigma R^k,
			\label{eq:ms-6}\\
			\int_\Omega H_R(Y)\dd Y&\leq CR^{k+\sigma}. \label{eq:ms-7}
		\end{align}
		Define
		\begin{align*}
			\mathcal D_R^{\rm slice}(Y)&:=\frac{\ell(Y)^{4-k}}{\mathsf M_i(t)}\,
			\eta_{i,t}(\mathscr Q_{R,R}(Y)),\\
			\mathcal D_R^{\rm time}(Y)&:=\frac{\ell(Y)^{2-k}}{\mathsf M_i(t)}
			\int_{|u-t|<R^2\ell(Y)^2}
			\eta_{i,u}(\mathscr Q_{R,R}(Y))\dd u.
		\end{align*}
		Then
		\begin{align}
			\int_\Omega \mathcal D_R^{\rm slice}(Y)\dd Y
			&\leq CR^k\frac{{\mathfrak r}^4}{\mathsf M_i(t)}\mathcal D_i(t), \label{eq:ms-8}\\
			\int_\Omega \mathcal D_R^{\rm time}(Y)\dd Y
			&\leq CR^k\frac{{\mathfrak r}^2}{\mathsf M_i(t)}
			\int_{|u-t|<R^2{\mathfrak r}^2}\mathcal D_i(u)\dd u. \label{eq:ms-9}
		\end{align}
		Here \eqref{eq:ms-6}--\eqref{eq:ms-9} are used only when $R{\mathfrak r}\leq c_0$.  In \eqref{eq:ms-7}, the additional
		restriction $2R{\mathfrak r}\leq {\mathfrak r}_0$ is also imposed when necessary to use \eqref{eq:ms-2}.
	\end{lem}
	
	\begin{proof}
		If $|y-Y|<R{\mathfrak r}h(Y)$ and $R{\mathfrak r}\leq c_0$, the Lipschitz bound for $h$ gives
		\begin{equation}
			(1-c_0)h(Y)\leq h(y)\leq(1+c_0)h(Y).
			\label{eq:ms-10}
		\end{equation}
		In particular, every cylinder in question lies inside $\Omega$ in the axial
		direction.  For a fixed $y$, the set of its possible centers $Y$ has volume
		at most $C(R{\mathfrak r}h(y))^k$, and on this set
		$\ell(Y)^{-k}\leq C({\mathfrak r}h(y))^{-k}$.  Therefore we get
		\[
		\int_\Omega F_{\delta,R}\dd Y
		\leq\frac{CR^k}{\mathsf M_i(t)}\mu_{i,t}(T_{2\delta {\mathfrak r}})
		\leq C\delta^\sigma R^k.
		\]
		The identical argument with $\delta=R$ gives \eqref{eq:ms-7}.  It gives the upper bound
		in \eqref{eq:ms-5} by using $T_{3{\mathfrak r}}$.
		
		For the lower bound in \eqref{eq:ms-5}, fix $(x,y)\in T_{\mathfrak r}$.  Every
		$Y\in B_{c {\mathfrak r}h(y)}(y)$, for a sufficiently small fixed $c>0$, belongs to
		$\Omega$ and satisfies
		\[
		|y-Y|<{\mathfrak r}h(Y),\qquad |x|<2{\mathfrak r}h(Y),\qquad h(Y)\asymp h(y).
		\]
		Its contribution to the $Y$ integral of $\ell(Y)^{-k}$ is bounded below by a
		positive constant. We integrate over $T_{\mathfrak r}$ to obtain \eqref{eq:ms-5}.
		
		For \eqref{eq:ms-8}, Fubini and \eqref{eq:ms-10} give
		\[
		\int_{\{Y\in\Omega:|y-Y|<R{\mathfrak r}h(Y)\}}\ell(Y)^{4-k}\dd Y
		\leq CR^k {\mathfrak r}^4h(y)^4\leq CR^k {\mathfrak r}^4.
		\]
		The transverse condition only reduces this integral.  This proves \eqref{eq:ms-8}.
		For \eqref{eq:ms-9}, similarly use
		\[
		\int_{\{Y\in\Omega:|y-Y|<R{\mathfrak r}h(Y)\}}\ell(Y)^{2-k}\dd Y
		\leq CR^k {\mathfrak r}^2h(y)^2\leq CR^k {\mathfrak r}^2
		\]
		and enlarge the time interval using $\ell(Y)\leq {\mathfrak r}$.
	\end{proof}
	
	\begin{lem}
		Let $0<\sigma<4$, and let ${\mathfrak r}(t)$ take dyadic values in $(0,{\mathfrak r}_0]$.  For every
		$R\geq1$ and every $u$,
		\begin{equation}
			\int_I {\mathfrak r}(t)^{2-\sigma}
			\1_{\{|u-t|<R^2{\mathfrak r}(t)^2\}}\dd t
			\leq C_\sigma R^2{\mathfrak r}_0^{4-\sigma}.
			\label{eq:ms-11}
		\end{equation}
		For the selected radii $\mathfrak r(t)=\mathfrak r_i(t)$ of Lemma~\ref{lem:measure-averaging},
		\begin{align}
			\int_I Q_i(t)\int_\Omega \mathcal D_R^{\rm slice}(Y)\dd Y\dd t
			&\leq CR^k{\mathfrak r}_0^{4-\sigma}\Delta_i, \label{eq:ms-12}\\
			\int_I Q_i(t)\int_\Omega \mathcal D_R^{\rm time}(Y)\dd Y\dd t
			&\leq CR^{k+2}{\mathfrak r}_0^{4-\sigma}\Delta_i. \label{eq:ms-13}
		\end{align}
		For \eqref{eq:ms-12}--\eqref{eq:ms-13}, the nonnegative functions
		$\mathcal D_R^{\rm slice}(Y)$ and $\mathcal D_R^{\rm time}(Y)$ are set equal to zero at any time $t$ when
		the corresponding admissibility restriction on $R\mathfrak r(t)$ fails.
	\end{lem}
	
	\begin{proof}
		Put ${\mathfrak r}_j:=2^{-j}{\mathfrak r}_0$. For fixed $u$, the times for which
		${\mathfrak r}(t)={\mathfrak r}_j$ and $|u-t|<R^2{\mathfrak r}_j^2$ form a set of length at most
		$2R^2{\mathfrak r}_j^2$. Hence we obtain
		\[
		\begin{aligned}
			\int_I{\mathfrak r}(t)^{2-\sigma}\1_{\{|u-t|<R^2{\mathfrak r}(t)^2\}}\dd t
			&\le\sum_{j\ge0}{\mathfrak r}_j^{2-\sigma}(2R^2{\mathfrak r}_j^2)\\
			&=2R^2{\mathfrak r}_0^{4-\sigma}\sum_{j\ge0}2^{-j(4-\sigma)}
			=\frac{2R^2{\mathfrak r}_0^{4-\sigma}}{1-2^{-(4-\sigma)}}.
		\end{aligned}
		\]
		This proves \eqref{eq:ms-11}. Next, \eqref{eq:ms-1} gives
		$Q_i(t)/\mathsf M_i(t)\le2{\mathfrak r}^{-\sigma}$, so \eqref{eq:ms-8} implies
		\[
		\int_IQ_i(t)\int_\Omega \mathcal D_R^{\rm slice}\dd Y\dd t
		\le CR^k\int_I {\mathfrak r}(t)^{4-\sigma}\mathcal D_i(t)\dd t \le CR^k{\mathfrak r}_0^{4-\sigma}\int_J\mathcal D_i(t)\dd t
		=CR^k{\mathfrak r}_0^{4-\sigma}\Delta_i.
		\]
		For the spacetime error, Tonelli's theorem, \eqref{eq:ms-9}, and
		\eqref{eq:ms-11} yield
		\[
		\begin{aligned}
			\int_IQ_i(t)\int_\Omega \mathcal D_R^{\rm time}\dd Y\dd t
			&\le CR^k\int_I {\mathfrak r}(t)^{2-\sigma}
			\int_{|u-t|<R^2{\mathfrak r}(t)^2}\mathcal D_i(u)\dd u\dd t\\
			&\le CR^k\int_J\mathcal D_i(u)
			\left(\int_I{\mathfrak r}(t)^{2-\sigma}
			\1_{\{|u-t|<R^2{\mathfrak r}(t)^2\}}\dd t\right)\dd u\\
			&\le CR^{k+2}{\mathfrak r}_0^{4-\sigma}\Delta_i.
		\end{aligned}
		\]
		The admissibility restriction ensures the enlarged time interval is
		contained in $J$. These are \eqref{eq:ms-12} and \eqref{eq:ms-13} and thus conclude the proof.
	\end{proof}

	Take $\varepsilon>0$ such that
	\begin{equation}
		0<\sigma-\varepsilon<\sigma+\varepsilon<4.
		\label{eq:ms-14}
	\end{equation}
	Let $\delta_j:=2^{-j}$ and $R_\ell:=2^\ell$.  Define
	\begin{align*}
		S(Y):=\sum_{j,\ell\geq0}
		\delta_j^{-\sigma+\varepsilon}R_\ell^{-k-\varepsilon}
		F_{\delta_j,R_\ell}(Y) +\sum_{\ell\geq0}
		R_\ell^{-k-\sigma-\varepsilon}H_{R_\ell}(Y).
	\end{align*}
	Only admissible terms are included.  From \eqref{eq:ms-6}--\eqref{eq:ms-7}, we get
	\begin{equation}
		\int_\Omega S(Y)\dd Y\leq C_\varepsilon. \label{eq:ms-15}
	\end{equation}
	Define the \textbf{combined error function} by
	\[
	\mathscr B(Y):=\sum_{\ell\geq0}R_\ell^{-k-3}
	\bigl(\mathcal D_{R_\ell}^{\rm slice}(Y)+\mathcal D_{R_\ell}^{\rm time}(Y)\bigr),
	\qquad b_i(t):=\int_\Omega \mathscr B(Y)\dd Y.
	\]
	Equations \eqref{eq:ms-12}--\eqref{eq:ms-13} imply
	\begin{equation}
		\int_I Q_i(t)b_i(t)\dd t\leq C\Delta_i. \label{eq:ms-16}
	\end{equation}
	
	\begin{thm}\label{thm:measure-selection}
		Assume \eqref{eq:ms-14} and $\Delta_i/\mathsf W_i\to0$. Assume also that the
		selected radii $\mathfrak r_i(t)$ tend to zero in probability with
		respect to the probability measure $Q_i(t)\dd t/\mathsf W_i$ on $I$;
		that is, for every fixed $s>0$,
		\begin{equation}\label{eq:ms-radius-probability}
			\frac1{\mathsf W_i}\int_{\{t\in I:\mathfrak r_i(t)\ge s\}}Q_i(t)\dd t
			\longrightarrow0.
		\end{equation}
		There are times $t_i\in I$ and axial
		centers $Y_i\in\Omega$ such that, with
		\[
		{\mathfrak r}_i:={\mathfrak r}_i(t_i)\longrightarrow0,\qquad
		h_i:=h(Y_i),\qquad \ell_i:={\mathfrak r}_i h_i,\qquad
		\mathsf P_i:=\mu_{i,t_i}(\mathscr Q_{2,1}(Y_i))>0,
		\]
		the following bounds hold for each fixed $R\geq1$ and $0<\delta\leq1$ and all sufficiently large $i$:
		\begin{align}
			\frac{\mu_{i,t_i}(\{|x|<\delta \ell_i,\ |y-Y_i|<R\ell_i\})}{\mathsf P_i}
			&\leq C\delta^{\sigma-\varepsilon}R^{k+\varepsilon}, \label{eq:ms-17}\\
			\frac{\mu_{i,t_i}(\{|x|<R\ell_i,\ |y-Y_i|<R\ell_i\})}{\mathsf P_i}
			&\leq CR^{k+\sigma+\varepsilon}. \label{eq:ms-18}
		\end{align}
		Moreover, there is a sequence $\epsilon_i\to0$ such that, for each fixed $R\ge1$, the following bounds hold for all sufficiently large $i$:
		\begin{align}
			\frac{\ell_i^4}{\mathsf P_i}
			\eta_{i,t_i}(\{|x|<R\ell_i,\ |y-Y_i|<R\ell_i\})
			&\leq \epsilon_iR^{k+3}, \label{eq:ms-19}\\
			\frac{\ell_i^2}{\mathsf P_i}
			\int_{|u-t_i|<R^2\ell_i^2}
			\eta_{i,u}(\{|x|<R\ell_i,\ |y-Y_i|<R\ell_i\})\dd u
			&\leq \epsilon_iR^{k+3}. \label{eq:ms-20}
		\end{align}
	\end{thm}
	
	\begin{proof}
		We will use the \textbf{selection probability measure}
		$\dd\pi_i(t):=\mathsf W_i^{-1}Q_i(t)\dd t$ on $I$. By \eqref{eq:ms-16}, we obtain
		\[
		\int_I b_i\dd\pi_i\le C\Delta_i/\mathsf W_i\longrightarrow0.
		\]
		If $\Delta_i>0$, then $b_i(t)>(\Delta_i/\mathsf W_i)^{1/2}$ on the
		set in the following integral. By integrating this inequality, we obtain
		\[
		\begin{aligned}
			\left(\frac{\Delta_i}{\mathsf W_i}\right)^{1/2}
			\pi_i\!\left(\left\{t\in I:b_i(t)>
			(\Delta_i/\mathsf W_i)^{1/2}\right\}\right)
			&\le\int_{\{t\in I:b_i(t)>(\Delta_i/\mathsf W_i)^{1/2}\}}
			b_i(t)\dd\pi_i(t)\\
			&\le\int_I b_i(t)\dd\pi_i(t)
			\le C\frac{\Delta_i}{\mathsf W_i}.
		\end{aligned}
		\]
		Dividing by $(\Delta_i/\mathsf W_i)^{1/2}>0$, we obtain
		\[
		\pi_i\!\left(\left\{t\in I:b_i(t)>
		(\Delta_i/\mathsf W_i)^{1/2}\right\}\right)
		\le C\left(\frac{\Delta_i}{\mathsf W_i}\right)^{1/2}\longrightarrow0.
		\]
		If $\Delta_i=0$, then $b_i=0$ almost everywhere for $\pi_i$.
		By \eqref{eq:ms-radius-probability}, there are $\widehat\delta_i\downarrow0$
		with $\pi_i\{{\mathfrak r}_i>\widehat\delta_i\}\to0$. After excluding these two sets and the
		times that are not admissible, we are left with a set of positive probability for all large $i$.
		Choose $t_i$ in that set. Then ${\mathfrak r}_i(t_i)\to0$ and
		$b_i(t_i)\to0$.
		At this time, \eqref{eq:ms-5} and \eqref{eq:ms-15} give
		\[
		c\leq\int_\Omega \mathfrak f(Y)\dd Y\leq C,\qquad \int_\Omega S(Y)\dd Y\leq C_\varepsilon,
		\qquad \int_\Omega \mathscr B(Y)\dd Y=b_i(t_i)\longrightarrow0.
		\]
		Choose a fixed $L>4C_\varepsilon/c$.  Then we have
		\[
		\int_{\{Y\in\Omega:S(Y)>L\mathfrak f(Y)\}}\mathfrak f(Y)\dd Y\leq L^{-1}\int_\Omega S(Y)\dd Y<c/4,
		\]
		If $b_i(t_i)>0$, the same direct comparison gives
		\[
		\int_{\{Y\in\Omega:\mathscr B(Y)>\sqrt{b_i(t_i)}\mathfrak f(Y)\}}
		\mathfrak f(Y)\dd Y
		\le\frac1{\sqrt{b_i(t_i)}}\int_\Omega\mathscr B(Y)\dd Y
		=\sqrt{b_i(t_i)}\longrightarrow0.
		\]
		If $b_i(t_i)=0$, then $\mathscr B=0$ almost everywhere on $\Omega$,
		so the same excluded set has zero $\mathfrak f(Y)\dd Y$ measure.
		Thus there is a center $Y_i$ with $\mathfrak f(Y_i)>0$ and
		\begin{equation}
			S(Y_i)\leq L \mathfrak f(Y_i),\qquad
			\mathscr B(Y_i)\leq\sqrt{b_i(t_i)}\mathfrak f(Y_i).
			\label{eq:ms-21}
		\end{equation}
		Since $\mathsf P_i=\mathsf M_i(t_i) \ell_i^k \mathfrak f(Y_i)$, divide each individual summand of $S$ and $\mathscr B$
		by $\mathfrak f(Y_i)$ to obtain \eqref{eq:ms-17}--\eqref{eq:ms-20}.  A fixed $R$ is admissible for all large
		$i$ because ${\mathfrak r}_i\to0$. Rounding $R$ and $\delta$ up to dyadic values
		extends the bounds to arbitrary fixed $R\ge1$ and $0<\delta\le1$, after increasing
		$C$ and multiplying $\epsilon_i$ by a fixed constant. Finally, $\ell_i/h_i={\mathfrak r}_i\to0$, so the axial boundary is
		sent to infinity even if $Y_i$ approaches $\partial\Omega$.
	\end{proof}
	
	\begin{lem}\label{lem:small-selected-radii}
		Put $p_i:=\int_I\mu_{i,t}(T_{{\mathfrak r}_0})\dd t$.  Suppose there are numbers $\mathsf b_i>0$ such that
		\begin{equation}
			p_i\leq C \mathsf b_i,\qquad \mathsf W_i/\mathsf b_i\longrightarrow\infty.
			\label{eq:ms-22}
		\end{equation}
		Then ${\mathfrak r}_i(t)\to0$ in $Q_i(t)\dd t/\mathsf W_i$ probability.  Moreover,
		$\Delta_i=o(\mathsf b_i)$ implies $\Delta_i=o(\mathsf W_i)$.
		In particular, if there are radii $s_i\leq {\mathfrak r}_0$ with
		\[
		\frac{\int_I\mu_{i,t}(T_{s_i})\dd t}{s_i^\sigma \mathsf b_i}
		\longrightarrow\infty,
		\]
		then $\mathsf W_i/\mathsf b_i\to\infty$.
	\end{lem}
	
	\begin{proof}
		If ${\mathfrak r}_i(t)\geq s$, \eqref{eq:ms-1} gives
		$Q_i(t)\leq2s^{-\sigma}\mu_{i,t}(T_{{\mathfrak r}_0})$.  Hence we get
		\begin{equation}\label{eq:ms-large-radius-probability}
			\frac1{\mathsf W_i}\int_{\{t\in I:{\mathfrak r}_i(t)\geq s\}}Q_i(t)\dd t
			\leq\frac{2s^{-\sigma}p_i}{\mathsf W_i}
			\leq\frac{Cs^{-\sigma}\mathsf b_i}{\mathsf W_i}\longrightarrow0.
		\end{equation}
		Moreover,
		\begin{equation}\label{eq:ms-small-radius-consequences}
			\frac{\Delta_i}{\mathsf W_i}=\frac{\Delta_i}{\mathsf b_i}\frac{\mathsf b_i}{\mathsf W_i}\longrightarrow0,
			\qquad
			\mathsf W_i=\int_I Q_i(t)\dd t
			\ge 2^{-\sigma}s_i^{-\sigma}\int_I\mu_{i,t}(T_{s_i})\dd t.
		\end{equation}
		The lower bound for $\mathsf W_i$ follows from \eqref{eq:ms-arbitrary-radius}.
	\end{proof}

	\subsection{Concentration at the axis}
	
	We use the rectangular estimate proved in this subsection in
	Lemma~\ref{lem:cyl-tube}, which controls concentration of the squared-distance
	measure near the axis. We first refine Theorem~\ref{thm:measure-selection}
	to obtain independent bounds in the transverse and axial directions.
	The resulting rectangular estimate \eqref{htps:rectangular-bound}
	gives the product-weighted integrability in the stationary tube argument.
	
	For a nonnegative Radon measure $\mu$ assigning zero mass to the axis,
	we call $\int_{E\cap\{r>0\}}r^{-2}\dd\mu$ its
	\textbf{transverse Hardy moment} on $E$, where $r:=|x|$.
	For the distance measure $\dd\mu:=j_a^2\rho\dd\mu_N$, this is
	$\int_{N\cap E\cap\{r>0\}}r^{-2}j_a^2\rho\dd\mu_N$.
	The relative tube estimate in Lemma~\ref{lem:cyl-tube} will control
	this moment by the dyadic calculation in \eqref{eq:cyl-relative-hardy}.
	
	To establish the rectangular estimate, use the radius-time selection
	notation from Subsection~\ref{gsel:section}.
	Fix $i$ and an admissible time $t$, and write $\mu:=\mu_{i,t}$.
	In particular
	$\Omega\subset\mathbb R^k$ is bounded,
	$h(y):=\min\{1,\operatorname{dist}(y,\mathbb R^k\setminus\Omega)\}$,
	$T_{\mathfrak r}:=\{|x|<{\mathfrak r}h(y),\ y\in\Omega\}$, and at an admissible time
	\[
	\mathsf M_i(t):=\mu(T_{\mathfrak r})>0,\qquad
	\mu(T_s)\le C(s/{\mathfrak r})^\sigma \mathsf M_i(t)\quad(0<s\le {\mathfrak r}_0),
	\qquad \ell(Y):={\mathfrak r}h(Y).
	\]
	The near-maximizing radius ${\mathfrak r}$ is dyadic. Along the
	sequence of times selected in Theorem~\ref{thm:measure-selection},
	these radii tend to zero.
	For $R,L\ge1$, define
	\[
	\mathscr Q_{R,L}(Y)
	:=\{|x|<R\ell(Y),\ |y-Y|<L\ell(Y)\},\qquad
	H_{R,L}(Y):=\frac{\mu(\mathscr Q_{R,L}(Y))}{\mathsf M_i(t)\ell(Y)^k}.
	\]
	A pair $(R,L)$ is called \textbf{admissible} if $L{\mathfrak r}<c_0$ and $CR{\mathfrak r}<{\mathfrak r}_0$, where $C$ is a fixed
	large geometric constant. Along that selected sequence, every fixed
	pair $(R,L)$ is admissible for all sufficiently large indices.
	
	\begin{lem}
		\label{htps:average}
		At every admissible time and for every admissible pair $(R,L)$,
		\[
		\int_\Omega H_{R,L}(Y)\dd Y\le C R^\sigma L^k.
		\]
	\end{lem}
	\begin{proof}
		For a fixed source point $(x,y)$, contributing centers satisfy
		$|y-Y|<L{\mathfrak r}h(Y)$. Since $L{\mathfrak r}<c_0$, the Lipschitz estimate gives
		$h(Y)\asymp h(y)$. Their volume is at most
		$C(L{\mathfrak r}h(y))^k$, and on that set $\ell(Y)^{-k}\le C({\mathfrak r}h(y))^{-k}$.
		Note that $|x|<R{\mathfrak r}h(Y)$ implies $|x|<CR{\mathfrak r}h(y)$.
		Therefore the transverse near-maximizer estimate \eqref{eq:ms-2} gives
		\[
		\int_\Omega H_{R,L}(Y)\dd Y
		\le\frac{CL^k}{\mathsf M_i(t)}\mu(T_{CR{\mathfrak r}}) \le C L^kR^\sigma.
		\]
	\end{proof}
	
	\begin{cor}\label{cor:rectangular-selection}
		Under the hypotheses of Theorem~\ref{thm:measure-selection}, the same
		$t_i,Y_i,\mathfrak r_i,\ell_i$, and $\mathsf P_i$ can be chosen to satisfy
		\eqref{eq:ms-17}--\eqref{eq:ms-20} and, for every fixed $R,L\ge1$
		and all sufficiently large $i$,
		\begin{equation}\label{htps:rectangular-bound}
			\frac{\mu_{i,t_i}(\{|x|<R\ell_i,\ |y-Y_i|<L\ell_i\})}{\mathsf P_i}
			\le C R^{\sigma+\varepsilon}L^{k+\varepsilon}.
		\end{equation}
		The constant is independent of $i,R,L$.
	\end{cor}
	\begin{proof}
		Let $\mathfrak f(Y):=\mu(\mathscr Q_{2,1}(Y))/(\mathsf M_i(t)\ell(Y)^k)$.
		Lemma~\ref{lem:measure-averaging} gives
		$c\le\int_\Omega \mathfrak f(Y)\dd Y\le C$.
		For admissible dyadic pairs $R,L\ge1$, add the nonnegative function
		\[
		S_{\rm rect}(Y):=\sum_{R,L}
		R^{-\sigma-\varepsilon}L^{-k-\varepsilon}H_{R,L}(Y)
		\]
		to the selection in Theorem~\ref{thm:measure-selection}. Lemma~\ref{htps:average}
		implies
		\[
		\int_\Omega S_{\rm rect}(Y)\dd Y\le
		C\sum_{j,\ell\ge0}2^{-\varepsilon(j+\ell)}<\infty.
		\]
		More precisely, define
		$E_{\rm rect}:=\{Y\in\Omega:S_{\rm rect}(Y)>C_*\mathfrak f(Y)\}$.
		Then
		\begin{equation}
			\int_{E_{\rm rect}}\mathfrak f(Y)\dd Y
			\le C_*^{-1}\int_{E_{\rm rect}}S_{\rm rect}(Y)\dd Y
			\le \frac{C_\varepsilon}{C_*}<\frac c4
		\end{equation}
		provided $C_*>4C_\varepsilon/c$, where $c$ is the lower bound in
		\eqref{eq:ms-5}. In the proof of Theorem~\ref{thm:measure-selection},
		the original growth exclusions remove less than $c/4$, and the error
		exclusion removes less than $c/4$ for all large $i$.
		Thus, after also excluding $E_{\rm rect}$, we retain $\mathfrak f(Y)\dd Y$ mass
		at least $c/4$, and \eqref{eq:ms-21} still holds. At the
		remaining center, $H_{R,L}(Y)\le C R^{\sigma+\varepsilon}L^{k+\varepsilon}\mathfrak f(Y)$,
		which gives \eqref{htps:rectangular-bound} after division by $\mathfrak f(Y)$.
		We can round each radius up to a dyadic radius with only a change in the constant.
		The radius--time error sum is still
		$\sum_j(2^{-j}{\mathfrak r}_0)^{4-\sigma}$, so it remains valid for every
		$0<\sigma<4$.
		The diagonal bound \eqref{eq:ms-18} is retained separately, with exponent
		$k+\sigma+\varepsilon$ instead of the less precise exponent obtained
		by setting $R=L$ in \eqref{htps:rectangular-bound}.
		
	\end{proof}
	
	\subsection{Stationary concentration and growth at infinity}
	For a varifold $V$, a center $z$, and a variance $h>0$, write
	$F_{z,h}(V)=(4\pi h)^{-n/2}\int_{\R^{n+1}} e^{-|Z-z|^2/(4h)}\,\dd\|V\|(Z)$.
	The following two lemmas use the ordinary Jacobi equation and
	Euclidean stationarity.
	
	\begin{lem}\label{lem:endpoint-static-jacobi}
		Assume \eqref{eq:cyl-S1}, and put $\tau_0=\max\{0,2-2\nu_0\}$.
		Let $U$ be a smooth solution of $JU=0$ on $\Reg(\barC)$. Suppose that,
		for some $p>\tau_0$ and $q<2+2\nu_0$,
		\begin{equation}\label{eq:endpoint-static-field-growth}
			\int_{\{r<\delta,\ |y|<R\}} U^2\dd\mu_{\barC}\le C_R\delta^p,
			\qquad \int_{B_R}U^2\dd\mu_{\barC}\le C R^{k+q}
			\quad(0<\delta<1,\ R\ge1).
		\end{equation}
		Then $U=0$.
	\end{lem}
	\begin{proof}
		For $s\in\R$, we write $U_s(Z)=U(e^sZ)$. The homogeneity of $J$ gives
		$JU_s=0$ and hence $\partial_s U_s=DU_s=(1-2L_{\barC})U_s$ on the
		regular part. Let $e=e_{j\ell\alpha}$ have Gaussian eigenvalue $\lambda$
		in \eqref{eq:cyl-gaussian-spectrum}, and put $M(s)=\int_{\barC} U_s e\,\dd g$.
		Fix a compact interval $K_0\subset\R$. By dilation and the ball bound
		in \eqref{eq:endpoint-static-field-growth}, for $s\in K_0$ and all sufficiently large $T$, we have
		\[
		\int_{B_T}U_s^2\dd\mu_{\barC}
		=e^{-ns}\int_{B_{e^sT}}U^2\dd\mu_{\barC}\le C_{K_0}T^{k+q}.
		\]
		Summing on the shells $2^hR<|Z|\le2^{h+1}R$ gives
		\begin{equation}\label{eq:static-dilate-gaussian-tail}
			\sup_{s\in K_0}\int_{\barC\setminus B_R}|U_s|^2\dd g
			\le C_{K_0}\sum_{h\ge0}(2^{h+1}R)^{k+q}e^{-4^hR^2/4}
			\longrightarrow0\qquad(R\to\infty).
		\end{equation}
		For $\psi\in C_c^\infty(\R)$ with support in $K_0$, put
		$v=\chi_\delta\zeta_R e$, using the cutoffs of Lemma~\ref{lem:all-root-moment-test}.
		Integration by parts in $s$ and Gaussian integration by parts on the regular cylinder give
		\[
		-\int_{\R}\psi'[U_s,v]_\rho\dd s
		=(1-2\lambda)\int_{\R}\psi[U_s,v]_\rho\dd s
		-2\int_{\R}\psi[U_s,[L_{\barC},\chi_\delta](\zeta_Re)
		+\chi_\delta[L_{\barC},\zeta_R]e]_\rho\dd s .
		\]
		Fix $R$ and set $A_{\delta,R}=\barC\cap B_{2R}\cap\{\delta<r<2\delta\}$.
		The tube bound in \eqref{eq:endpoint-static-field-growth}, after dilation,
		gives $\int_{A_{\delta,R}}|U_s|^2\dd g\le C_{K_0,R}\delta^p$ for small $\delta$.
		By \eqref{eq:cyl-cutoff-commutator} and the Cauchy--Schwarz inequality,
		uniformly for $s\in K_0$,
		\[
		\begin{aligned}
			|[U_s,[L_{\barC},\chi_\delta](\zeta_Re)]_\rho|
			&\le C_{R,e}\delta^{\gamma_j^+-2}
			\left(\int_{A_{\delta,R}}|U_s|^2\dd g\right)^{1/2}
			\left(\int_{A_{\delta,R}}1\,\dd g\right)^{1/2}\\
			&\le C_{K_0,R,e}\delta^{\gamma_j^+-2+(p+m)/2}
			=C_{K_0,R,e}\delta^{\nu_j+(p-2)/2}\longrightarrow0 .
		\end{aligned}
		\]
		Here $\nu_j+(p-2)/2\ge\nu_0+(p-2)/2>0$. After letting $\delta\downarrow0$, we bound the remaining
		commutator term by a constant times
		\[
		\sup_{s\in\operatorname{spt}\psi}\|U_s\|_{L^2(\dd g;\,|Z|>R)}
		\bigl(\|e\|_{L^2(\dd g)}+R^{-2}\|De\|_{L^2(\dd g)}\bigr).
		\]
		By \eqref{eq:static-dilate-gaussian-tail}, we can let $R\to\infty$ to obtain $M'=(1-2\lambda)M$ in distributions. The same tube
		and tail bounds and smooth convergence on compact subsets of the regular
		part give continuity of $M$. It follows that
		\[
		M(s)=e^{(\gamma_j^++2\ell+|\alpha|)s}M(0).
		\]
		For $s\ge0$, put $t=e^s$ and $A_h=B_{2^{h+1}t}\setminus B_{2^ht}$.
		Changing variables $X=tZ$ and using \eqref{eq:endpoint-static-field-growth}, we obtain
		\[
		\begin{aligned}
			\|U_s\|_{L^2(\dd g)}^2
			&=c_n t^{-n}\int_{\barC}|U(X)|^2e^{-|X|^2/(4t^2)}\dd\mu_{\barC}(X)\\
			&\le c_n t^{-n}\left(\int_{B_t}U^2\dd\mu_{\barC}
			+\sum_{h\ge0}e^{-4^h/4}\int_{A_h}U^2\dd\mu_{\barC}\right)\\
			&\le C t^{k+q-n}\left(1+\sum_{h\ge0}2^{(h+1)(k+q)}e^{-4^h/4}\right)
			\le C e^{(q-m)s}.
		\end{aligned}
		\]
		Put $b=\gamma_j^++2\ell+|\alpha|$. Since $b\ge-(m-2)/2+\nu_0>(q-m)/2$,
		the Cauchy--Schwarz inequality gives
		\[
		|M(0)|=e^{-bs}|M(s)|
		\le e^{-bs}\|U_s\|_{L^2(\dd g)}\|e\|_{L^2(\dd g)}
		\le C_e e^{((q-m)/2-b)s}\longrightarrow0 .
		\]
		All Gaussian moments vanish, so Proposition~\ref{prop:spectrum} gives $U=0$.
	\end{proof}
	
	\begin{lem}
		\label{lem:endpoint-static-tube}
		Assume \eqref{eq:cyl-S1} and \eqref{eq:cyl-S4} with $B>k/2$. Write $d(Z)=\dist(Z,\barC)$ and $\Theta=F(\barC)$. Let
		$V_i$ be stationary integral $n$-varifolds in $\R^{n+1}$ with uniform
		all-center area ratios.
		Suppose $F_{z,h}(V_i)\le\Theta$ for every $z,h$, and that
		$\eta_i\downarrow0$ with, for every fixed $R\ge1$,
		\[
		\int_{B_R} d^2\dd\|V_i\|\le C_R\eta_i^2.
		\]
		For each $\tau_0<\sigma<2$ and $L\ge1$, for all sufficiently large $i$,
		\begin{equation}\label{eq:endpoint-static-tube}
			\int_{\{r<\delta,\ |y|<L\}}d^2\dd\|V_i\|
			\le C_{L,\sigma}\delta^\sigma\eta_i^2\qquad(0<\delta\le1).
		\end{equation}
	\end{lem}
	\begin{proof}
		Suppose \eqref{eq:endpoint-static-tube} fails for some fixed $L,\sigma$.
		Write $\mu_i=d^2\|V_i\|$ and $T_{\delta,L}=\{r<\delta,\ |y|<L\}$.
		After passing to a subsequence, there are $0<\delta_i\le1$ such that
		\[
		\frac{\mu_i(T_{\delta_i,L})}{\delta_i^\sigma\eta_i^2}\longrightarrow\infty .
		\]
		Since $T_{\delta_i,L}\subset B_{L+1}$, the local height hypothesis gives
		$\mu_i(T_{\delta_i,L})\le\mu_i(B_{L+1})\le C_{L+1}\eta_i^2$.
		The displayed ratio is therefore at most $C_{L+1}\delta_i^{-\sigma}$,
		so $\delta_i\to0$. Take $\Omega=B_{L+2}^k$ in the measure selection,
		so $h=1$ on $B_L^k$. We use the time-independent measures $\mu_{i,t}=\mu_i$
		on a fixed interval $J$ and take the error measures identically zero.
		Choose a unit interval $I$ with $I+[-c_0^2,c_0^2]\subset J$.
		These are inputs to Theorem~\ref{thm:measure-selection}, which concerns
		Radon measures. On a fixed bounded axial cylinder the area-ratio bound and
		a cover by $O_L(a^{-k})$ balls of radius $2a$ give
		$\|V_i\|\{r<a,\ y\in\Omega\}\le C_La^m$. Since $d\le r$, the
		supremum in \eqref{eq:ms-maximal-mass} is finite, bounded at small
		scales by $\sup_{a\le1}C_La^{m+2-\sigma}$. A dyadic radius between
		$\delta_i$ and $2\delta_i$ shows $\mathsf W_i/\eta_i^2\to\infty$.
		The mass of the largest working tube is $O(\eta_i^2)$ by the local
		height bound. Lemma~\ref{lem:small-selected-radii} gives the required
		smallness of the selected radii.
		
		Fix $0<\varepsilon<\min\{\sigma-\tau_0,2-\sigma,2B-k\}$, and put
		$p=\sigma-\varepsilon$, $q=\sigma+\varepsilon$.
		Theorem~\ref{thm:measure-selection} and Corollary~\ref{cor:rectangular-selection}
		give axial centers $z_i=(0,Y_i)$, scales $a_i\downarrow0$ and
		$P_i=\mu_i\{r<2a_i,\ |y-Y_i|<a_i\}>0$. Set
		\[
		W_i=(Z\mapsto (Z-z_i)/a_i)_\#V_i,\qquad
		\kappa_i^2=P_i/a_i^{n+2},\qquad \lambda_i=\kappa_i^{-2}d^2\|W_i\|.
		\]
		Axial translation invariance and homogeneity of $d$ give
		$\lambda_i\{r<2,|y|<1\}=1$. The area bound and $d\le r$ give
		$0<\kappa_i\le C$. We then obtain
		\begin{align}
			\lambda_i\{r<b,\ |y|<L'\}&\le C b^p (L')^{k+\varepsilon}
			&&(0<b\le1,\ L'\ge1),\label{eq:static-selected-inner}\\
			\lambda_i\{r<R,\ |y|<L'\}&\le C R^q(L')^{k+\varepsilon}
			&&(R,L'\ge1),\label{eq:static-selected-rectangle}\\
			\lambda_i(B_R)&\le C R^{k+q} &&(R\ge1).\label{eq:static-selected-ball}
		\end{align}
		For each fixed choice of radii these bounds hold for all sufficiently large $i$.
		The ball bound \eqref{eq:static-selected-ball} follows from \eqref{eq:ms-18} directly. By stationary
		integral-varifold compactness \cite{Allard}, we pass to a subsequence along which $W_i\to W$ and
		$\kappa_i\to\kappa\ge0$. The normalizing cylinder has positive $d^2\|W_i\|$ measure, so for each $i$ we can choose
		a point in this cylinder at which $W_i$ has density at least one. Stationary monotonicity gives a uniform positive area on
		a fixed ball about that point, contained in $B_5$. A compact cutoff
		equal to one on $B_5$ preserves this lower bound in the limit. Under
		translation and dilation $F_{z,h}(W_i)=F_{z_i+a_i z,a_i^2h}(V_i)\le\Theta$.
		For fixed $z,h$, let $0\le\chi_R\uparrow1$ be compactly supported
		cutoffs equal to one on $B_R$. Local convergence and monotone convergence give
		\[
		F_{z,h}(W)=\lim_{R\to\infty}\lim_{i\to\infty}
		\int_{\R^{n+1}}\chi_R\rho_{z,h}\dd\|W_i\|\le\Theta .
		\]
		This bounds the multiplicity in the constancy argument below.
		Indeed $F_{0,1}((D_{1/R})_\#W)=F_{0,R^2}(W)\le\Theta$.
		If a blowdown is $j|\barC|$, the same cutoff argument gives
		$j\Theta=F_{0,1}(j|\barC|)\le\Theta$, hence $j\le1$.
		
		If $\kappa>0$, then $\lambda_i$ converges locally to
		$\kappa^{-2}d^2\|W\|$, which is nonzero by a compact cutoff equal to
		one on the normalizing cylinder. To estimate its weighted height, partition each radius into $[0,2)$
		and $[2^j,2^{j+1})$, $j\ge1$, and denote the resulting rectangles
		in $(r,|y|)$ by $E_{j\ell}$, indexing the first interval by zero.
		The rectangular bound \eqref{eq:static-selected-rectangle} passes to the limiting measure by testing with
		cutoffs supported in $\{r<2^{j+2},\ |y|<2^{\ell+2}\}$. Hence we have
		\[
		\begin{aligned}
			\int_{E_{j\ell}}\frac{d^2}{(1+r^2)(1+|y|^2)^B}\dd\|W\| \le C2^{-2j-2B\ell}\int_{E_{j\ell}}d^2\dd\|W\|\le C 2^{-j(2-q)}2^{-\ell(2B-k-\varepsilon)} .
		\end{aligned}
		\]
		Both exponents are positive by the choice of $\varepsilon$, so summing
		over $j,\ell\ge0$ gives $\mathcal E_B(W;\barC)<\infty$. We next verify the unit
		blowdown. For $W_R=(D_{1/R})_\#W$, the diagonal estimate gives, for
		fixed $L'\ge1$,
		\[
		\int_{B_{L'}}d^2\dd\|W_R\|
		\le C_{L'}R^{k+q-n-2}\longrightarrow0.
		\]
		Every subsequential blowdown is therefore supported on $\barC$.
		Its mass on the axis is zero by the area bound and the same covering
		argument as above. Constancy on the connected regular part gives
		$j|\barC|$. Monotonicity centered at a positive-density point of $W$
		gives positive limiting mass on a ball after rescaling, so $j\ge1$.
		The inherited bound $F_{0,1}(j|\barC|)\le\Theta$ forces $j=1$.
		Compactness now gives full multiplicity-one blowdown, so \eqref{eq:cyl-S4} gives
		$W=|\barC|$, contradicting the nonzero limiting height measure.
		
		If $\kappa=0$, the local height bounds imply $\operatorname{spt}\|W\|
		\subset\barC$. The same constancy and Gaussian bound give $W=|\barC|$.
		On compact regular subsets Allard's multiplicity-one regularity
		theorem \cite{Allard} gives single normal graphs over $\barC$; the
		minimal graph equation and interior elliptic estimates give smooth
		convergence to zero. On any compact set disjoint from $\barC$ there
		can be no extra support, since monotonicity would preserve positive
		mass in a ball disjoint from the limit. Write $u_i$ for the graph
		heights. In the regular normal tube $d=|u_i|$ and the graph Jacobians
		converge uniformly to one on compact sets. The local height bound thus
		gives uniform local $L^2$ bounds for $u_i/\kappa_i$. We subtract the
		minimal graph equation at zero and integrate its linearization along
		$t u_i$, $0\le t\le1$. This yields a homogeneous uniformly elliptic
		linear equation for $u_i/\kappa_i$ on each compact regular patch,
		by \eqref{eq:averaged-graph-equation} with $v=0$ and $\epsilon=0$.
		As $u_i\to0$ smoothly there, the coefficients converge smoothly to those of
		$J=\Delta_{\barC}+|A_{\barC}|^2$ by \eqref{eq:graph-mean-curvature-linearization}. Interior
		elliptic estimates and a diagonal subsequence give smooth convergence
		to $U$ with $JU=0$. In particular $\lambda_i\to U^2\mu_{\barC}$
		on compact subsets away from the axis. Choose $b>0$ with $Cb^p<1/2$
		in the first selected estimate at $L'=1$. Then
		for a compactly supported cutoff $0\le\chi\le1$ away from the axis,
		equal to one on $\{b\le r\le2,\ |y|\le1\}$, we have
		\[
		\int_{\barC}\chi U^2\dd\mu_{\barC}
		=\lim_i\int_{\R^{n+1}}\chi\dd\lambda_i
		\ge\liminf_i\bigl(\lambda_i(T_{2,1})-\lambda_i(T_{b,1})\bigr)
		\ge1-Cb^p>\tfrac12 .
		\]
		In particular $U\ne0$.
		Then we exhaust the regular part of each tube and ball by compact subsets.
		The bounds \eqref{eq:static-selected-inner} and \eqref{eq:static-selected-ball} then give
		\eqref{eq:endpoint-static-field-growth}, with the displayed $p,q$.
		They obey $p>\tau_0$ and $q<2<2+2\nu_0$. Lemma~\ref{lem:endpoint-static-jacobi}
		gives $U=0$, the final contradiction. This concludes the proof.
	\end{proof}
	
	\begin{prop}\label{prop:endpoint-static-growth}
		Assume \eqref{eq:cyl-S1} and \eqref{eq:cyl-S4}. Let $V$ be a nonzero
		stationary integral $n$-varifold with finite all-center area ratios and
		full multiplicity-one blowdown $|\barC|$. Write
		\[
		H(R):=\int_{B_R(0)}\dist(Z,\barC)^2\,\dd\|V\|(Z).
		\]
		If $H(R)\le C R^d$ for all $R\ge1$ and some $d<\min\{n+2,k+2+2\nu_0\}$, then $V=|\barC|$.
	\end{prop}
	\begin{proof}
		We first record the Gaussian area bound needed in the application of
		Lemma~\ref{lem:endpoint-static-tube}. With
		$\rho_{z,h}(Z)=(4\pi h)^{-n/2}\exp(-|Z-z|^2/(4h))$, stationarity
		and the first variation with the field $(Z-z)\rho_{z,h}$ give
		\[
		\frac{\dd}{\dd h}\int_{\R^{n+1}}\rho_{z,h}\,\dd\|V\|
		=\frac1{4h^2}\int_{\R^{n+1}\times\mathrm{Gr}(n,n+1)} |(Z-z)^{\perp_T}|^2\rho_{z,h}(Z)\,\dd V(Z,T)
		\ge 0 .
		\]
		Here $\mathrm{Gr}(n,n+1)$ denotes the Grassmannian of unoriented $n$-planes in $\R^{n+1}$.
		Indeed the tangential divergence of this field is
		$[n-|(Z-z)^{\top_T}|^2/(2h)]\rho_{z,h}$, whereas differentiation
		of the kernel uses the full squared length. The area bound controls
		the errors of outer cutoffs and permits differentiation under the
		integral on each compact subinterval of $h>0$. At fixed $z$, rescale
		by $h^{-1/2}$. The rescaled center tends to zero, the rescaled varifolds
		converge to $|\barC|$, and the uniform area ratios control the Gaussian
		tails. Consequently the increasing integral tends to $F(\barC)$.
		Thus the Gaussian bound holds at every center and variance for $V$,
		and also for every dilation of $V$.
		
		If $H(R)$ vanishes identically, the support is contained in $\barC$.
		The constancy theorem on its connected regular part gives an integer
		multiple there. The axis has zero $n$-dimensional measure since $k<n$.
		The blowdown fixes the multiplicity to be one. We can therefore assume
		that $H(R)>0$ for every sufficiently large $R$. Choose
		\[
		\max\{d,0\}<d'<\min\{n+2,k+2+2\nu_0\}.
		\]
		For each sufficiently large integer $i$, the positive numbers
		$2^{-jd'}H(2^j)$ tend to zero as $j\to\infty$. Their maximum over
		$j\ge i$ is attained. Let $R_i=2^{j_i}$ be a maximizing radius.
		Then $R_i\ge 2^i$ and, for every $L\ge1$, choosing
		$h\in\mathbb Z_{\ge0}$ with $L\le2^h<2L$ gives
		\[
		H(LR_i)\le H(2^hR_i)\le 2^{hd'}H(R_i)
		\le (2L)^{d'}H(R_i).
		\]
		Set $V_i=(D_{1/R_i})_\#V$ and
		$\eta_i^2=R_i^{-n-2}H(R_i)>0$. As $d<n+2$, we have $\eta_i\to0$.
		We may pass to a subsequence for which $\eta_i$ decreases to zero.
		Homogeneity of the distance gives
		\begin{equation}\label{eq:endpoint-outer-normalization}
			\eta_i^{-2}\int_{B_1}\dist(Z,\barC)^2\,\dd\|V_i\|=1,\qquad
			\eta_i^{-2}\int_{B_L}\dist(Z,\barC)^2\,\dd\|V_i\|
			\le (2L)^{d'}\qquad(L\ge1).
		\end{equation}
		The full blowdown gives $V_i\to|\barC|$. We apply
		Lemma~\ref{lem:endpoint-static-tube} with
		$\max\{0,2-2\nu_0\}<p<2$. For every fixed $L\ge1$, this gives
		\begin{equation}\label{eq:endpoint-outer-tube}
			\eta_i^{-2}\int_{\{r<\delta,\ |y|<L\}}\dist(Z,\barC)^2\,\dd\|V_i\|
			\le C_L\delta^p
		\end{equation}
		for all sufficiently large $i$ and every $0<\delta\le1$.
		
		Since $V_i\to|\barC|$ with multiplicity one and
		\eqref{eq:endpoint-outer-normalization} bounds the normalized height on each ball,
		the argument for $\kappa_i\to0$ in the proof of Lemma~\ref{lem:endpoint-static-tube}
		applies with $\eta_i$ in place of $\kappa_i$. Along a subsequence the
		normalized graph heights converge locally smoothly to a field $U$ with
		$J_{\barC}U=0$ on $\Reg(\barC)$, and away from the axis we have
		\[
		\eta_i^{-2}\dist(\cdot,\barC)^2\,\|V_i\|
		\rightharpoonup U^2\mu_{\barC}
		\quad\text{locally as Radon measures}.
		\]
		
		We verify the normalization and bounds for this field. Fix $\delta>0$
		with $C_2\delta^p<1/2$. The part of the normalized squared-distance
		measure of $B_1$ on $r\ge\delta$ has mass at least $1/2$ by
		\eqref{eq:endpoint-outer-normalization} and \eqref{eq:endpoint-outer-tube}.
		A smooth compactly supported cutoff equal to one on
		$\overline B_1\cap\{r\ge\delta\}$ and supported where $r>\delta/2$
		therefore has limiting integral at least $1/2$ against
		$U^2\dd\mu_{\barC}$. In particular $U\ne0$. By compact exhaustion
		away from the axis and monotone convergence, the two estimates imply
		\[
		\int_{\barC\cap B_L}U^2\,\dd\mu_{\barC}\le (2L)^{d'},\qquad
		\int_{\barC\cap\{r<\delta,\ |y|<L\}}U^2\,\dd\mu_{\barC}
		\le C_L\delta^p .
		\]
		Taking $q=d'-k$, we have $k+q=d'$ and $q<2+2\nu_0$.
		Thus these are precisely the two bounds in \eqref{eq:endpoint-static-field-growth}. Lemma~\ref{lem:endpoint-static-jacobi} now gives $U=0$,
		contradicting the retained mass.
	\end{proof}
	
	\subsection{The tube estimate for the rescaled flow}
	In this subsection, we return to the time-dependent problem. Put \(n=m+k\),
	\(\barC=\C\times\R^k\), and \(r=|x|\). The stationary growth estimate
	permits us to choose a tube exponent above two.
	
	The next lemma plays the role of the near-vertex estimates in \cite[Section~3]{FH1}. The singular set here is an axis, so transverse concentration must be estimated together with the axial distribution of the distance measure.
	\begin{lem}\label{lem:cyl-tube}
		Assume \eqref{eq:cyl-S1} and \eqref{eq:cyl-S4}. Let \(N_s^i\)
		be rescaled flows with uniform all-center area ratios on a fixed bounded interval \(J\). Suppose that the area
		measures $\mu_{N_{s}^i}$ converge locally weakly as Radon measures
		to $\mu_{\barC}$, uniformly for $s\in J$.
		Let $a_i\to0$ and $\barC_i:=\barC_{a_i}$ be fixed reduced references.
		Write $W_i:=W_{(a_i,0)}$ for their homogeneous squared-distance
		functions from Lemma~\ref{lem:distance-function}. Suppose that the reference
		mean curvatures and their fixed-order angular derivatives are
		bounded by \(C_\ell \mathfrak q_i/r\), where \(\mathfrak q_i=O(\eta_i)\) and
		\(\eta_i\downarrow0\), and that
		\begin{equation}\label{eq:cyl-tube-input}
			\int_J\int_{N_{s}^i}|\Phi_i|^2\rho\dd\mu_{N_{s}^i}\dd s=o(\eta_i^2),
			\qquad
			\int_{J'}\int_{N_{s}^i\cap B_L}W_i\rho\dd\mu_{N_{s}^i}\dd s
			\le C_{J',L}\eta_i^2
		\end{equation}
		for every \(J'\Subset J\) and fixed \(L\).
		Fix \(I\Subset J\), \(L>0\), and
		\[
		2<\sigma<\min\{4,2+2\nu_0\}.
		\]
		Then there are \(\delta_0>0\) and \(C<\infty\) such that, for all
		sufficiently large \(i\) and \(0<\delta\le\delta_0\),
		\begin{equation}\label{eq:cyl-relative-tube}
			\int_I\int_{N_{s}^i\cap\{r<\delta,\ |y|<L\}}
			W_i\rho\dd\mu_{N_{s}^i}\dd s\le C\delta^\sigma\eta_i^2.
		\end{equation}
	\end{lem}
	
	\begin{proof}
		\noindent\textbf{Step 1. Selection of a transverse scale and an axial center.}
		All Gaussian densities in the fixed larger spatial region are
		uniformly comparable to one. On each bounded axial region, the area
		bound gives
		\begin{equation}\label{eq:cyl-local-tube-area}
			\mu_{N_{s}^i}\{r<\delta,\ |y|<L'\}\le C_{L'}\delta^m.
		\end{equation}
		To prove \eqref{eq:cyl-local-tube-area}, cover the axial region by \(O_{L'}(\delta^{-k})\) balls of
		radius \(\delta\), and apply \eqref{equ:areabd} to the corresponding
		ambient balls of radius $C\delta$. We sum their areas to obtain
		$C_{L'}\delta^{-k}\delta^{m+k}=C_{L'}\delta^m$, proving
		\eqref{eq:cyl-local-tube-area}.
		
		Consequently, for
		$0<\delta_0\le1$, we have
		\begin{equation}\label{eq:cyl-local-hardy-bound}
			\int_{N_{s}^i\cap\{0<r<\delta_0,\ |y|<L'\}}r^{-2}\dd\mu_{N_{s}^i}
			\le C_{L'}\sum_{j\ge0}(2^{-j}\delta_0)^{m-2}<\infty.
		\end{equation}
		Also $W_i\le Cr^2$ gives
		$\int_{N_{s}^i\cap T_{\mathfrak r}}W_i\dd\mu_{N_{s}^i}\le C\mathfrak r^{m+2}$
		for the variable tubes used below.
		Since $\sigma<4<m+2$, the supremum defining $Q_i(s)$ is finite.
		Choose a fixed interval $J_1$ with $I\Subset J_1\Subset J$;
		choose $0<c_0<1/100$ so that $I+[-c_0^2,c_0^2]\Subset J_1$,
		and then fix $0<\mathfrak r_0<c_0/100$. In the error sums, restrict to $R\mathfrak r_i(t)\le c_0$. The corresponding time intervals
		thus lie in $J_1$, and every fixed $R$ is eventually admissible
		along the selected sequence.
		
		Choose \(\Omega\) so that the function \(h\) of the measure selection
		is bounded from below on \(\{|y|<L\}\). Use the measure \(W_i\dd\mu_{N_{s}^i}\) and the error measure
		\[
		\bigl(|\Phi_i|^2+\mathfrak q_i^2/r^2\bigr)\dd\mu_{N_{s}^i}
		\]
		in the fixed cylinder $\mathscr D:=\{|x|<1\}\times\Omega$.
		By \eqref{eq:cyl-tube-input}, \eqref{eq:cyl-local-hardy-bound}, and
		$\mathfrak q_i=O(\eta_i)$, its integral over $J_1$ satisfies
		\begin{equation}\label{eq:cyl-selection-total-defect}
			\begin{aligned}
				\Delta_i&:=\int_{J_1}\int_{N_{s}^i\cap\mathscr D}
				(|\Phi_i|^2+\mathfrak q_i^2r^{-2})\dd\mu_{N_{s}^i}\dd s\\
				&\le C\int_J\int_{N_{s}^i}|\Phi_i|^2\rho\dd\mu_{N_{s}^i}\dd s
				+C|J_1|\mathfrak q_i^2
				\le o(\eta_i^2)+C\eta_i^2.
			\end{aligned}
		\end{equation}
		Recall from \eqref{eq:ms-maximal-mass} that, for
		$\dd\mu_{i,s}:=W_i\dd\mu_{N_{s}^i}$, we have
		\[
		\mathsf W_i=\int_I\sup_{\mathfrak r\in\{2^{-j}\mathfrak r_0:j\in\mathbb Z_{\ge0}\}}
		\mathfrak r^{-\sigma}\mu_{i,s}(T_{\mathfrak r})\dd s.
		\]
		If \eqref{eq:cyl-relative-tube} fails, pass to a subsequence and choose
		$\delta_i\downarrow0$ such that
		\[
		\frac{1}{\delta_i^\sigma\eta_i^2}
		\int_I\int_{N_{s}^i\cap\{r<\delta_i,\ |y|<L\}}
		W_i\rho\dd\mu_{N_{s}^i}\dd s\longrightarrow\infty.
		\]
		The radii must tend to zero by the local position bound
		\eqref{eq:cyl-tube-input}. Set
		$h_*:=\inf_{|y|<L}h(y)>0$ and choose a dyadic radius
		$\widehat{\mathfrak r}_i$ with
		$\delta_i/h_*\le\widehat{\mathfrak r}_i<2\delta_i/h_*$.
		Then $\{r<\delta_i,\ |y|<L\}\subset T_{\widehat{\mathfrak r}_i}$, so
		\begin{equation}\label{eq:cyl-selection-divergent-mass}
			\frac{\mathsf W_i}{\eta_i^2}
			\ge\frac{\widehat{\mathfrak r}_i^{-\sigma}}{\eta_i^2}
			\int_I\mu_{i,s}(T_{\widehat{\mathfrak r}_i})\dd s \ge\frac{c}{\delta_i^\sigma\eta_i^2}
			\int_I\int_{N_{s}^i\cap\{r<\delta_i,\ |y|<L\}}
			W_i\rho\dd\mu_{N_{s}^i}\dd s\longrightarrow\infty.
		\end{equation}
		Also \eqref{eq:cyl-tube-input} gives
		$\int_I\mu_{i,s}(T_{\mathfrak r_0})\dd s\le C\eta_i^2$.
		Lemma~\ref{lem:small-selected-radii}, with $\mathsf b_i:=\eta_i^2$, therefore
		implies that the selected radii tend to zero in $\pi_i$ probability.
		Finally, \eqref{eq:cyl-selection-total-defect} and
		\eqref{eq:cyl-selection-divergent-mass} give
		\[
		0\le\frac{\Delta_i}{\mathsf W_i}
		\le C\frac{\eta_i^2}{\mathsf W_i}\longrightarrow0,
		\]
		as required in
		Theorem~\ref{thm:measure-selection}.
		
		Choose
		\[
		0<\varepsilon<\min\{\sigma-2,\ 4-\sigma,\ 2+2\nu_0-\sigma\}.
		\]
		Corollary~\ref{cor:rectangular-selection}
		gives slice times \(s_i\in I\), bounded axial centers \(Y_i\), scales
		\(\ell_i\downarrow0\), and masses \(\mathsf P_i>0\) on the
		normalizing cylinders of transverse radius \(2\ell_i\) and axial
		radius \(\ell_i\). Write \(Z_i:=(0,Y_i)\) and put
		\begin{equation}\label{eq:cyl-rescaled-normalization}
			\begin{gathered}
				\mathcal M_i(t):=\ell_i^{-1}\bigl(N_{s_i+\ell_i^2t}^i-Z_i\bigr),\qquad
				\kappa_i^2:=\frac{\mathsf P_i}{\ell_i^{n+2}},\\
				\widehat\mu_i:=\kappa_i^{-2}W_i\dd\mu_{\mathcal M_i(0)}.
			\end{gathered}
		\end{equation}
		Since \(W_i\le Cr^2\), the area bound gives \(0<\kappa_i\le C\).
		The homogeneity of the squared-distance function and its axial independence give
		\(\widehat\mu_i(\{r<2,\ |y|<1\})=1\). For each fixed \(R,L'\ge1\)
		and \(0<\delta\le1\), the selected inequalities are
		\begin{equation}\label{eq:cyl-micro-selection}
			\begin{aligned}
				\widehat\mu_i\{r<\delta,\ |y|<R\}
				&\le C\delta^{\sigma-\varepsilon}R^{k+\varepsilon},\\
				\widehat\mu_i\{r<R,\ |y|<L'\}
				&\le CR^{\sigma+\varepsilon}(L')^{k+\varepsilon},\\
				\widehat\mu_i(B_R)&\le CR^{k+\sigma+\varepsilon}.
			\end{aligned}
		\end{equation}
		
		\medskip
		\noindent\textbf{Step 2. Stationary compactness and the Gaussian density bound.}
		The normal velocity of \(\mathcal M_i\) is
		\[
		\Psi_i:=\mathbf H_{\mathcal M_i}
		+\tfrac12(\ell_i Z_i+\ell_i^2 X)^\perp.
		\]
		The selected error bounds imply, on every fixed dilated
		ball $B_R$ and bounded time interval $[-T_0,T_0]$, with $R,T_0>0$ fixed,
		\begin{equation}\label{eq:cyl-micro-defect}
			\begin{gathered}
				\kappa_i^{-2}\int_{\mathcal M_i(0)\cap B_R}
				\bigl(|\Psi_i|^2+\mathfrak q_i^2/r^2\bigr)\dd\mu_{\mathcal M_i(0)}
				\longrightarrow0,\\
				\kappa_i^{-2}\int_{-T_0}^{T_0}\int_{\mathcal M_i(t)\cap B_R}
				\bigl(|\Psi_i|^2+\mathfrak q_i^2/r^2\bigr)\dd\mu_{\mathcal M_i(t)}\dd t
				\longrightarrow0 .
			\end{gathered}
		\end{equation}
		The normalization factors are \(\ell_i^4/\mathsf P_i\) on a time slice and
		\(\ell_i^2/\mathsf P_i\) in spacetime. Both the squared shrinker quantity
		and the reference term have this scaling. The dilated weights
		\[
		\rho_i(X):=\exp\bigl(-\ell_i Z_i\cdot X/2-\ell_i^2|X|^2/4\bigr)
		\]
		tend smoothly to one on compact sets. For a compactly supported
		smooth function \(\varphi\), the rescaled first-variation formula
		\cite[Section~2.1]{FH1} gives
		\begin{equation}\label{eq:cyl-micro-transport}
			\frac{\dd}{\dd t}\int_{\mathcal M_i(t)}\varphi\rho_i\dd\mu_{\mathcal M_i(t)}
			=\int_{\mathcal M_i(t)}
			\bigl(D\varphi\cdot\Psi_i-\varphi|\Psi_i|^2\bigr)\rho_i\dd\mu_{\mathcal M_i(t)}.
		\end{equation}
		At time zero the absolute local mean curvature tends to zero in
		\(L^2\), by \eqref{eq:cyl-micro-defect} and the vanishing drift.
		Integral-varifold compactness \cite{Allard} gives a stationary integral limit \(V\).
		Let $\mathcal K_\varphi$ contain the support of $\varphi$ and fix
		$t,t'\in[-T_0,T_0]$. Set
		$\Delta_i^\varphi:=\int_{-T_0}^{T_0}\int_{\mathcal M_i(\tau)\cap \mathcal K_\varphi}|\Psi_i|^2\rho_i
		\dd\mu_{\mathcal M_i(\tau)}\dd\tau=o(\kappa_i^2)$.
		By integrating \eqref{eq:cyl-micro-transport}, we obtain
		\[
		\left|\int_{\mathcal M_i(t)}\varphi\rho_i\dd\mu_{\mathcal M_i(t)}
		-\int_{\mathcal M_i(t')}\varphi\rho_i\dd\mu_{\mathcal M_i(t')}\right|
		\le C_{\varphi,T_0}(\Delta_i^\varphi)^{1/2}+\|\varphi\|_\infty \Delta_i^\varphi\longrightarrow0.
		\]
		Here the first term follows from the Cauchy--Schwarz inequality and the uniform
		area bound on $\mathcal K_\varphi$, and $\kappa_i\le C$.
		By taking $t'=0$ and using $\rho_i\to1$ on compact sets, we conclude
		that $\mu_{\mathcal M_i(t)}$ converge locally weakly as Radon measures
		to $\|V\|$, uniformly for $t\in[-T_0,T_0]$.
		
		Next, we verify both nonvanishing and the needed Gaussian upper bound
		directly. The ordinary flows
		\[
		\widetilde{\mathcal M}_i(\tau):=\ell_i^{-1}
		\left(\sqrt{1-\ell_i^2\tau}\,
		N_{s_i-\log(1-\ell_i^2\tau)}^i-Z_i\right)
		\]
		are the exact parabolic dilations of the original mean curvature flows
		\(M_{u}^i:=\sqrt{-u}\,N_{-\log(-u)}^i\), centered at
		\((e^{-s_i/2}Z_i,-e^{-s_i})\), with spatial scale
		\(e^{-s_i/2}\ell_i\). On each fixed dilated cylinder the time
		change tends to the identity and the spatial change does as well
		because \(Z_i\) is bounded. The normalizing slice contains a point
		in a fixed bounded ball. If \(V=0\), apply Huisken's monotonicity
		formula \cite{Huisken} to \(\widetilde{\mathcal M}_i\), centered at
		this point and time zero. Its density is one, whereas the Gaussian
		integral at any fixed earlier time tends to zero by local convergence
		and the uniform Gaussian tail bound. For fixed $h>0$ and bounded centers $Z_0$, the all-center area-ratio bound gives
		\begin{equation}\label{eq:cyl-micro-gaussian-tails}
			\int_{\{|X-Z_0|>R\}}\rho_{Z_0,h}\dd\mu_{\widetilde{\mathcal M}_i(-h)}
			\le C_h\sum_{j\ge0}(2^{j+1}R)^n
			\exp\!\left(-\frac{4^jR^2}{4h}\right)
			\longrightarrow0\quad(R\to\infty),
		\end{equation}
		uniformly in $i$; here $\rho_{Z_0,h}$ is defined below.
		For the chosen points $X_i^*\in\mathcal M_i(0)$,
		Huisken's monotonicity gives
		\[
		1\le\int_{\widetilde{\mathcal M}_i(-h)}\rho_{X_i^*,h}
		\dd\mu_{\widetilde{\mathcal M}_i(-h)}\longrightarrow0
		\qquad\text{if }V=0.
		\]
		This contradiction gives $V\ne0$.
		
		Fix a dilated center \(Z_0\), center time \(\tau_*\), and
		backward time difference \(h>0\). Write
		\[
		\rho_{Z_0,h}(X):=(4\pi h)^{-n/2}\exp(-|X-Z_0|^2/(4h)).
		\]
		For a hypersurface or varifold $V'$, define
		$F_{x,h}(V'):=\int_{\R^{n+1}}\rho_{x,h}\dd\|V'\|$, with the area measure understood
		when $V'$ is a hypersurface. In the original coordinates put
		\[
		x_i^*:=e^{-s_i/2}(Z_i+\ell_i Z_0),\qquad
		t_i^*:=e^{-s_i}(-1+\ell_i^2\tau_*).
		\]
		The corresponding backward time difference is $e^{-s_i}\ell_i^2h$.
		Choose a fixed $H>0$ sufficiently small that $t_i^*-H$ remains in
		one compact subinterval of the original time interval for all large $i$;
		this is possible because $I\Subset J$.
		For each fixed $h>0$, eventually $e^{-s_i}\ell_i^2h<H$, and
		Huisken's monotonicity formula \cite{Huisken} gives
		\begin{equation}\label{eq:cyl-two-scale-gaussian}
			\begin{aligned}
				F_{Z_0,h}(\widetilde{\mathcal M}_i(\tau_*-h))
				&=F_{x_i^*,e^{-s_i}\ell_i^2h}(M^i_{t_i^*-e^{-s_i}\ell_i^2h})\\
				&\le F_{x_i^*,H}(M^i_{t_i^*-H})\longrightarrow F(\barC)=:\Theta.
			\end{aligned}
		\end{equation}
		Indeed, after passing to a subsequence, $x_i^*\to x_*$ for an axial
		point $x_*\in\{0\}\times\R^k$. The measures $\mu_{M_u^i}$ converge locally weakly as Radon
		measures to $\mu_{\barC}$, uniformly on this compact original time
		interval, because $M_u^i=\sqrt{-u}\,N_{-\log(-u)}^i$ and $\barC$
		is dilation invariant. The area-ratio bound controls the Gaussian
		tails, and axial translation invariance and dilation invariance give
		$F_{x_*,H}(\barC)=F(\barC)$. This proves the limit in
		\eqref{eq:cyl-two-scale-gaussian}.
		With $h$ fixed, we let $i\to\infty$ to obtain
		\begin{equation}\label{eq:cyl-micro-density-bound}
			\int_{\R^{n+1}}\rho_{Z_0,h}\dd\|V\|\le\Theta
			\qquad(Z_0\in\mathbb R^{n+1},\ h>0).
		\end{equation}
		
		\medskip
		\noindent\textbf{Step 3. A positive normalized distance after dilation.}
		Suppose first that \(\kappa_i\to\kappa>0\). The measures
		\(\widehat\mu_i\) converge locally to
		\(\kappa^{-2}W_0\dd\|V\|\). This measure is nonzero,
		as seen by a compact cutoff equal to one on the closure of the
		normalizing cylinder. The diagonal bound in
		\eqref{eq:cyl-micro-selection} and comparability of $W_0$ with
		squared distance give, for $R\ge1$,
		\[
		\int_{B_R}\dist(Z,\barC)^2\dd\|V\|\le C R^{k+\sigma+\varepsilon}.
		\]
		We can first pass to the limit with a continuous cutoff
		supported in $B_{2R}$. By the choice of $\varepsilon$ the exponent
		is less than $k+2+2\nu_0$ and also less than $n+2$, since
		$\sigma+\varepsilon<4<m+2$. We verify the blowdown required by
		Proposition~\ref{prop:endpoint-static-growth}.
		
		For a fixed $L>0$, homogeneity of distance and the ball estimate in \eqref{eq:cyl-micro-selection} give
		\begin{equation}
			\begin{aligned}
				\int_{B_L}\dist(X,\barC)^2\dd\|(D_{1/R})_\#V\|(X)
				&=R^{-n-2}\int_{B_{LR}}\dist(X,\barC)^2\dd\|V\|(X)\\
				&\le C_L R^{k+\sigma+\varepsilon-n-2}
				=C_LR^{\sigma+\varepsilon-m-2}\longrightarrow0.
			\end{aligned}
		\end{equation}
		Here $\sigma+\varepsilon<4<m+2$. Thus every subsequential
		blowdown is supported on $\barC$.
		Every such limit is stationary and integral by varifold compactness
		\cite{Allard}.
		Its multiplicity is a constant integer on the connected regular
		part by the constancy theorem \cite{SimonGMT}. The area bound gives no mass
		to the axis: on any fixed region its \(\delta\)-tube has mass
		\(O(\delta^m)\) by the covering estimate \eqref{eq:cyl-local-tube-area}, applied to the limiting area-ratio bound. Stationary monotonicity \cite{SimonGMT} gives a fixed
		positive density at infinity, so all subsequential blowdowns equal
		the same \(j|\barC|\), \(j\ge1\). Take $h=R^2$ in
		\eqref{eq:cyl-micro-density-bound}. Full blowdown convergence
		and the area-ratio tail bound give
		\[
		j\Theta=F(j|\barC|)
		=\lim_{R\to\infty}F((D_{1/R})_\#V)
		=\lim_{R\to\infty}F_{0,R^2}(V)\le\Theta.
		\]
		Since $\Theta>0$ and $j\ge1$, we have $j=1$.
		Thus Proposition~\ref{prop:endpoint-static-growth} gives \(V=|\barC|\), contrary to the
		nonzero squared-distance measure. This implies $\kappa_i\to 0$.
		
		\medskip
		\noindent\textbf{Step 4. A vanishing normalized distance after dilation.}
		Assume \(\kappa_i\to0\). The bounded normalized
		position on each fixed ball shows \(\operatorname{spt}\|V\|\subset \barC\).
		The constancy theorem \cite{SimonGMT} and
		\eqref{eq:raw-tube-area} give \(V=j|\barC|\).
		By taking \(Z_0\) on the axis in \eqref{eq:cyl-micro-density-bound},
		we obtain \(j\le1\), so \(V=|\barC|\). Applying
		Brakke's local regularity theorem \cite[Section~6.11]{Brakke} to the ordinary flows
		\(\widetilde{\mathcal M}_i\) on any bounded dilated interval, we see that the
		dilated flows have complete single graphs of heights \(u_i\)
		over \(\barC_i\) on each fixed regular compact patch, smoothly tending
		to zero. To see that no extra support can persist on a compact set
		disjoint from the limiting support, suppose there are points of
		the ordinary dilated flows converging to such a point at an
		interior time. The backward Gaussian integral centered at each
		point is at least one by monotonicity and smoothness. At a fixed
		small backward time these integrals converge to that of the
		static cylinder by local measure convergence and the uniform
		Gaussian tails. At the limiting center this last integral tends
		to zero with the backward time, since the distance to the
		cylinder is positive. This gives a contradiction. We claim that,
		after passing to a subsequence and identifying the reference links,
		$u_i/\kappa_i$ converges smoothly on compact regular spacetime cylinders
		to a nonzero time-independent function $U$ satisfying
		\eqref{eq:cyl-limit-jacobi-equation} and the growth bounds
		\eqref{eq:cyl-micro-limit-growth}.

		We first describe the normalized limit on these patches. A compact
		test \(\chi^2W_i\) in \eqref{eq:cyl-micro-transport} has
		\(|D(\chi^2W_i)|^2\le C_\chi\chi^2W_i\) and is uniformly bounded
		on its support. Smooth approximation is legitimate by the squared-distance function
		bounds in Lemma~\ref{lem:distance-function}. Define
		\begin{equation}
			\mathcal P_i^\chi(t):=\int_{\mathcal M_i(t)}\chi^2W_i\rho_i
			\dd\mu_{\mathcal M_i(t)},\qquad
			\mathcal D_i^\chi(t):=\int_{\mathcal M_i(t)\cap\operatorname{spt}\chi}
			|\Psi_i|^2\rho_i\dd\mu_{\mathcal M_i(t)}.
		\end{equation}
		Then \eqref{eq:cyl-micro-transport} gives
		\begin{equation}
			|(\mathcal P_i^\chi)'|\le C_\chi\sqrt{\mathcal P_i^\chi \mathcal D_i^\chi}+C_\chi \mathcal D_i^\chi.
		\end{equation}
		Fix an interval $J_0$ containing zero, set $T_0:=|J_0|$, and put
		$M_i:=\sup_{t\in J_0} \mathcal P_i^\chi(t)$ and $\Delta_i^\chi:=\int_{J_0} \mathcal D_i^\chi(t)\dd t$.
		Equations \eqref{eq:cyl-micro-selection} and \eqref{eq:cyl-micro-defect} give $\mathcal P_i^\chi(0)\le C_\chi\kappa_i^2$
		and $\Delta_i^\chi=o(\kappa_i^2)$. By integrating in either time direction, we obtain
		\begin{equation}\label{eq:cyl-local-position-absorption}
			M_i\le C_\chi\kappa_i^2
			+C_\chi\sqrt{M_i}\sqrt{T_0 \Delta_i^\chi}+C_\chi \Delta_i^\chi \le C_\chi\kappa_i^2+\tfrac12M_i+C_{\chi,T_0}\Delta_i^\chi.
		\end{equation}
		We absorb $M_i/2$ in \eqref{eq:cyl-local-position-absorption} to obtain
		$\sup_{t\in J_0}\mathcal P_i^\chi(t)=M_i\le C_{\chi,T_0}\kappa_i^2$.
		
		Choose a relatively compact ambient neighborhood $\mathcal U\Subset B_R$
		away from the axis that meets $\barC$ in positive area. Smooth graphical
		convergence gives
		$\int_{\mathcal M_i(0)\cap \mathcal U}r^{-2}\dd\mu_{\mathcal M_i(0)}\ge c_{\mathcal U}>0$, so we have
		\[
		c_{\mathcal U}\frac{\mathfrak q_i^2}{\kappa_i^2}
		\le\kappa_i^{-2}\int_{\mathcal M_i(0)\cap B_R}
		\mathfrak q_i^2r^{-2}\dd\mu_{\mathcal M_i(0)}\longrightarrow0.
		\]
		Since $c_{\mathcal U}>0$, this proves $\mathfrak q_i/\kappa_i\to0$ and hence
		that the reference mean curvature divided by $\kappa_i$ tends to zero
		on every compact regular patch.
		
		The graph equation for \(u_i/\kappa_i\), obtained from
		\eqref{eq:normal-graph-flow} under the dilation
		\eqref{eq:cyl-rescaled-normalization}, is uniformly parabolic on
		smaller regular cylinders. Its coefficients converge smoothly
		to those of \(\partial_t-\Delta_{\barC}-|A_{\barC}|^2\).
		The inhomogeneous term is the normalized reference mean curvature,
		which tends smoothly to zero by \(\mathfrak q_i/\kappa_i\to0\) and the
		angular derivative bounds. The affine drift at the reference has
		zero normal component, since \(Z_i\) is axial and \(X\) is tangent
		to \(\barC_i\). Thus the normalized reference mean curvature is the only inhomogeneous term. Interior estimates, using the local spacetime \(L^2\)
		bound \eqref{eq:cyl-local-position-absorption}, give smooth convergence on compact regular
		cylinders, along a subsequence, to a function \(U\). On a compact regular patch write
		$F_i(P,t):=P+u_i(P,t)\mathbf n_{\barC_i}(P)$ and choose the graph normal
		$\mathbf n_i$ with $\mathbf n_{\barC_i}\cdot\mathbf n_i\ge1/2$.
		The normal-velocity equation gives
		\[
		\partial_t(u_i/\kappa_i)
		=\frac{\Psi_i(F_i(P,t),t)\cdot\mathbf n_i}
		{\kappa_i\,\mathbf n_{\barC_i}\cdot\mathbf n_i}.
		\]
		For every compact regular patch $\mathcal K_{\rm reg}$ and bounded time interval $J_0$,
		the graph area elements are uniformly comparable. Choosing $R$ with
		$F_i(\mathcal K_{\rm reg}\times J_0)\subset B_R$ and using \eqref{eq:cyl-micro-defect}, we get
		\begin{equation}
			\int_{J_0}\int_{\mathcal K_{\rm reg}}|\partial_t(u_i/\kappa_i)|^2\dd\mu_{\barC_i}\dd t
			\le C\kappa_i^{-2}\int_{J_0}\int_{\mathcal M_i(t)\cap B_R}
			|\Psi_i|^2\dd\mu_{\mathcal M_i(t)}\dd t
			\longrightarrow0.
		\end{equation}
		Here $\mathcal K_{\rm reg}$ is identified with the corresponding patch of $\barC_i$.
		Smooth convergence and the limiting graph equation therefore imply
		\begin{equation}\label{eq:cyl-limit-jacobi-equation}
			\partial_tU=0,\qquad
			(\Delta_{\barC}+|A_{\barC}|^2)U=0\quad(r>0).
		\end{equation}
		The normalized slice measures converge to \(U^2\dd\mu_{\barC}\)
		away from the axis, because the squared-distance function is exactly squared distance
		in the regular tube. The first estimate in
		\eqref{eq:cyl-micro-selection} excludes any residual mass on the
		axis. To make the nonvanishing explicit, choose a fixed
		$\delta>0$ with $C\delta^{\sigma-\varepsilon}<1/2$. The normalization
		\eqref{eq:cyl-rescaled-normalization} and \eqref{eq:cyl-micro-selection} give
		\[
		\widehat\mu_i\{\delta\le r<2,\ |y|<1\}
		\ge1-C\delta^{\sigma-\varepsilon}>\tfrac12.
		\]
		Choose a compactly supported cutoff $0\le\chi\le1$ away from the axis,
		equal to one on $\{\delta\le r\le2,\ |y|\le1\}$.
		Smooth convergence on its support yields
		$\int_{\barC}\chi U^2\dd\mu_{\barC}
		=\lim_i\int_{\R^{n+1}}\chi\dd\widehat\mu_i\ge1/2$; hence $U\ne0$. By passing to the limit in \eqref{eq:cyl-micro-selection}, we obtain
		\begin{equation}\label{eq:cyl-micro-limit-growth}
			\int_{\barC\cap\{r<\delta,\ |y|<R\}}U^2\dd\mu_{\barC}
			\le C_R\delta^{\sigma-\varepsilon},\qquad
			\int_{\barC\cap B_R}U^2\dd\mu_{\barC}\le CR^{k+\sigma+\varepsilon}.
		\end{equation}

		\medskip
		\noindent\textbf{Step 5. Vanishing of the ordinary Jacobi field.}
		Set $p=\sigma-\varepsilon$ and $q=\sigma+\varepsilon$. The choice
		of $\varepsilon$ gives $p>2>\tau_0$ and $q<2+2\nu_0$.
		By \eqref{eq:cyl-micro-limit-growth}, the squared $L^2$ norm of $U$
		is $O_R(\delta^p)$ on bounded tubes and $O(R^{k+q})$ on balls.
		Its equation \eqref{eq:cyl-limit-jacobi-equation} is the ordinary
		Jacobi equation on $\Reg(\barC)$. Lemma~\ref{lem:endpoint-static-jacobi}
		therefore gives $U=0$, contradicting the positive normalization
		obtained in Step~4. This proves the tube estimate.
	\end{proof}

	The next normal-angle argument is analogous to the radial-tilt estimate in \cite[Lemma~2.11]{FH1}. The cylindrical conclusion additionally controls the axial normal component, with a local inverse-radius distance term near the axis.
	\begin{lem}\label{lem:axial-tilt-estimate}
		Let $\barC_a:=\C_a\times\R^k$ be a reduced reference, and align its axis with $\{0\}\times\R^k$. Put $Z:=(x,y)$,
		$r:=|x|$, $r_y:=|y|$, and $j_a(Z):=\dist(Z,\barC_a)$. There are
		fixed constants $r_c>0$, $L_{\rm ax}<\infty$, and $C<\infty$ such
		that every smooth properly embedded boundaryless hypersurface $N$
		with the fixed all-center area-ratio bound satisfies
		\begin{equation}\label{eq:cyl-tilt-local-to-global}
			\int_N|\mathbf n_y|^2\rho\dd\mu_{N}
			\le C\left(D_a(N)^2+G(N)^2+|q(a)|^2
			+\int_{N\cap\{r<2r_c,\ r_y<L_{\rm ax}+1\}}
			\frac{j_a^2}{r^2}\rho\dd\mu_{N}\right).
		\end{equation}
		Here $D_a(N)^2:=\int_N j_a^2\rho\dd\mu_{N}$ and
		$\mathbf n:=(\mathbf n_x,\mathbf n_y)$ is any local unit normal to $N$.
		Now let $\barC_{a_i}$ be fixed references on an interval
		$I$. If
		\[
		\int_I\bigl(D_{a_i}(N_s^i)^2+G(N_s^i)^2+|q(a_i)|^2\bigr)\dd s
		\le C\eta_i^2
		\]
		and there are $\sigma>2$, $\delta_0>0$, $C<\infty$, and $i_0$ such that
		for every $i\ge i_0$ and every $0<\delta\le\delta_0$ we have
		\[
		\int_I\int_{N_s^i\cap\{r<\delta,\ r_y<L_{\rm ax}+1\}}
		j_{a_i}^2\rho\dd\mu_{N_s^i}\dd s\le C\delta^\sigma\eta_i^2,
		\]
		then for every $i\ge i_0$ we have
		\[
		\int_I\int_{N_s^i}|\mathbf n_y|^2
		\rho\dd\mu_{N_s^i}\dd s\le C\eta_i^2 .
		\]
	\end{lem}

	\begin{proof}
		Use $\mathcal T_{a,\varepsilon}$, $l_a$, $\mathbf n_a$, $\vartheta_a$,
		$w$, and the cutoff $\zeta$ from Lemma~\ref{lem:weighted-normal-angle}
		and its proof.
		The Gaussian moments are uniformly bounded, and we may assume $G(N)<\infty$.
		
		Next, we take a smooth ambient cutoff $\chi$
		equal to one on $\{r<r_c,\ r_y<L_{\rm ax}\}$, supported in
		$\{r<2r_c,\ r_y<L_{\rm ax}+1\}$, and test with $\chi^2\zeta^2l_a\mathbf n_a$.
		Here $r_c,L_{\rm ax}$ are any fixed positive constants and will be
		chosen in the proof of \eqref{eq:cyl-tilt-divergence-regions}.
		Apply \eqref{eq:cyl-tilt-tubular-divergence} and the Young inequalities
		in \eqref{eq:cyl-tilt-absorption-terms} with $\chi^2$ in place of $w$. The new cutoff term satisfies
		\[
		2|\chi\zeta^2l_a\,\nabla_N\chi\cdot\nabla_N l_a|
		\le \epsilon\chi^2\zeta^2\vartheta_a^2+C_\epsilon l_a^2|D\chi|^2 .
		\]
		The tube transition now costs $C l_a^2/r^2$ and the
		reference term is bounded by $C(l_a^2/r^2+|q(a)|^2)$ on the
		support of $\chi$. The axis cutoff error again tends to zero:
		its gradient is $O(\delta^{-1})$ on $\delta<r<2\delta$ and
		$|l_a|\le Cr$ there, so the absolute error is at most
		\[
		C\delta^{-1}\int_{N\cap\{\delta<r<2\delta,\ r_y<L_{\rm ax}+1\}}\zeta^2|l_a|\rho\dd\mu_N
		\le C\mu_N\{r<2\delta,\ r_y<L_{\rm ax}+1\}\le C\delta^m\longrightarrow0.
		\]
		Thus, if
		\[
		\mathfrak e_a(N):=D_a(N)^2+G(N)^2+|q(a)|^2,\qquad
		\mathcal H_a:=\int_{N\cap\{r<2r_c,\ r_y<L_{\rm ax}+1\}}r^{-2}j_a^2\rho\dd\mu_{N} ,
		\]
		the cutoff calculation gives
		\[
		\int_{N\cap\mathcal T_{a,\varepsilon}\cap\{r<r_c,r_y<L_{\rm ax}\}}
		|\mathbf n_y|^2\rho\dd\mu_N
		\le\int_{N\cap\mathcal T_{a,\varepsilon}\cap\{r<r_c,r_y<L_{\rm ax}\}}
		\vartheta_a^2\rho\dd\mu_N\le C(\mathfrak e_a(N)+\mathcal H_a).
		\]
		On $\{r>0\}\setminus\mathcal T_{a,\varepsilon}$,
		$|\mathbf n_y|^2\le1\le\varepsilon^{-2}j_a^2/r^2$;
		the axis has zero area by \eqref{eq:raw-tube-area}. Hence we get
		\begin{equation}\label{eq:cyl-tilt-compact-axis}
			\int_{N\cap\{r<r_c,\ r_y<L_{\rm ax}\}}|\mathbf n_y|^2\rho\dd\mu_{N}
			\le C(\mathfrak e_a(N)+\mathcal H_a).
		\end{equation}
		
		For the field defined in \eqref{eq:cyl-tilt-product-field},
		\eqref{eq:cyl-tilt-product-divergence} gives zero weighted divergence
		on every reference cylinder. For a radial function of
		$t$ put $\mathcal L_j h:=jh+t h'-t^2h/2$. Choose $h_x$ a
		smooth compactly supported function of $r^2$, equal to one near
		$r=0$, and put $h_y(r_y):=(1+r_y^2)^{-1}$. Define
		\begin{equation}\label{eq:cyl-tilt-product-field}
			\alpha:=-h_x\mathcal L_kh_y,\qquad
			\beta:=h_y\mathcal L_mh_x,\qquad
			\mathcal X(Z):=\alpha x+\beta y.
		\end{equation}
		The vector field is smooth and bounded. More explicitly, we have
		\[
		\mathcal L_kh_y
		=\frac{k}{1+r_y^2}-\frac{2r_y^2}{(1+r_y^2)^2}
		-\frac{r_y^2}{2(1+r_y^2)}.
		\]
		Thus $\alpha\to1/2$ and $\beta\to0$ as $r_y\to\infty$
		where $h_x=1$, uniformly for bounded $r$.
		Put $\mathfrak a:=x\cdot\mathbf n_x$,
		$\mathfrak b:=y\cdot\mathbf n_y$, and $\mathfrak t:=|\mathbf n_y|^2$. The weighted divergence of $\mathcal X$
		on $N$ is
		\begin{equation}\label{eq:cyl-tilt-product-divergence}
			\operatorname{div}_{N,\rho}\mathcal X
			=(\alpha-\beta)\mathfrak t-\frac{\alpha_r}{r}\mathfrak a^2
			-\frac{\beta_{r_y}}{r_y}\mathfrak b^2
			-\left(\frac{\alpha_{r_y}}{r_y}+\frac{\beta_r}{r}\right)\mathfrak a\mathfrak b .
		\end{equation}
		Indeed, \eqref{eq:cyl-tilt-product-field} gives
		\[
		\mathcal L_m\alpha+\mathcal L_k\beta
		=-(\mathcal L_mh_x)(\mathcal L_kh_y)
		+(\mathcal L_kh_y)(\mathcal L_mh_x)=0,
		\]
		and consequently
		\[
		\begin{aligned}
			\operatorname{div}_{N,\rho}\mathcal X
			&=\operatorname{div}\mathcal X-\tfrac12Z\cdot\mathcal X
			-\langle D\mathcal X\mathbf n,\mathbf n\rangle\\
			&=\alpha-\alpha(1-\mathfrak t)-\beta\mathfrak t
			-\frac{\alpha_r}{r}\mathfrak a^2-\frac{\beta_{r_y}}{r_y}\mathfrak b^2
			-\left(\frac{\alpha_{r_y}}{r_y}+\frac{\beta_r}{r}\right)\mathfrak a\mathfrak b\\
			&=(\alpha-\beta)\mathfrak t
			-\frac{\alpha_r}{r}\mathfrak a^2-\frac{\beta_{r_y}}{r_y}\mathfrak b^2
			-\left(\frac{\alpha_{r_y}}{r_y}+\frac{\beta_r}{r}\right)\mathfrak a\mathfrak b.
		\end{aligned}
		\]
		This proves \eqref{eq:cyl-tilt-product-divergence}.
		The quotients in this formula extend continuously at zero.
		
		Fix $r_c$ so small that $h_x=1$ on $r\le2r_c$ and
		$\mathcal L_m h_x=m-r^2/2>0$ there. On $r<r_c$ we have
		$\alpha_r=0$, $-\beta_{r_y}/r_y\ge0$, and
		$\alpha-\beta\to1/2$ as $r_y\to\infty$ uniformly in $r$.
		Moreover
		$r_y|\alpha_{r_y}/r_y+\beta_r/r|\le C$.
		Since $|\mathfrak b|\le r_y\sqrt{\mathfrak t}$, choose $L_{\rm ax}$ large and use Young's
		inequality to get
		\begin{equation}\label{eq:cyl-tilt-divergence-regions}
			\operatorname{div}_{N,\rho}\mathcal X\ge c\mathfrak t-C\mathfrak a^2
			\quad(r<r_c,\ r_y\ge L_{\rm ax}),\qquad
			\operatorname{div}_{N,\rho}\mathcal X\ge-C(\mathfrak t+\mathfrak a^2)
			\quad(r<r_c,\ r_y<L_{\rm ax}).
		\end{equation}
		On $N\cap\mathcal T_{a,\varepsilon}$, $\mathfrak a^2\le C(l_a^2+r^2\vartheta_a^2)$:
		the transverse radial component of $\mathbf n_a$ is $l_a$ and the
		normal of $N$, with its sign chosen locally, differs from
		$\mathbf n_a$ by at most $C\vartheta_a$. On $N\cap(\{r>0\}\setminus\mathcal T_{a,\varepsilon})$, we have
		$\mathfrak a^2\le r^2\le C j_a^2$. Thus the integral of $\mathfrak a^2$ on
		$r<r_c$ is bounded by $C\mathfrak e_a(N)$ using
		\eqref{eq:cyl-weighted-caccioppoli}. On $r\ge r_c$, all
		coefficients in \eqref{eq:cyl-tilt-product-divergence},
		including those of $\mathfrak b^2$ and $\mathfrak a\mathfrak b$ multiplied by $r_y^2$ and
		$r_y$ respectively, are bounded. They vanish outside a fixed
		transverse compact set. On $N\cap\mathcal T_{a,\varepsilon}$, $\mathfrak t\le\vartheta_a^2$,
		$|\mathfrak b|\le r_y\vartheta_a$, and $\mathfrak a^2\le C(l_a^2+r^2\vartheta_a^2)$.
		On $N\cap(\{r>0\}\setminus\mathcal T_{a,\varepsilon})$ the divergence is bounded and $j_a\ge c r_c$.
		Consequently, we have
		\begin{equation}\label{eq:cyl-tilt-outer-divergence}
			\int_{N\cap\{r\ge r_c\}}|\operatorname{div}_{N,\rho}\mathcal X|\rho\dd\mu_{N}
			\le C\mathfrak e_a(N).
		\end{equation}
		By \eqref{eq:cyl-weighted-caccioppoli}, we also obtain $\int_{N\cap\{r\ge r_c\}}|\mathcal X^\perp|^2\rho\dd\mu_{N}
		\le C\mathfrak e_a(N)$, since $\mathcal X\cdot\mathbf n=\alpha \mathfrak a+\beta \mathfrak b$ and
		$|\beta|r_y$ is bounded. On $r<r_c$ we have
		$|\mathcal X^\perp|^2\le C(\mathfrak a^2+\mathfrak t)$.
		
		Let $T_{\mathrm{far}}:=\int_{N\cap\{r<r_c,\ r_y\ge L_{\rm ax}\}}\mathfrak t\rho\dd\mu_{N}$.
		Young's inequality and the bounds for $|\mathcal X^\perp|^2$ give
		\begin{equation}\label{eq:cyl-tilt-source-absorption}
			\left|\int_N\Phi\cdot \mathcal X\rho\dd\mu_{N}\right|
			\le \epsilon\int_N|\mathcal X^\perp|^2\rho\dd\mu_{N}
			+C_\epsilon G(N)^2
			\le C\epsilon T_{\mathrm{far}}+C_\epsilon(\mathfrak e_a(N)+\mathcal H_a).
		\end{equation}
		Equations \eqref{eq:cyl-tilt-divergence-regions},
		\eqref{eq:cyl-tilt-outer-divergence},
		\eqref{eq:cyl-tilt-compact-axis}, and
		\eqref{eq:cyl-tilt-source-absorption} yield
		\[
		cT_{\mathrm{far}}-C(\mathfrak e_a(N)+\mathcal H_a)
		\le\int_N\operatorname{div}_{N,\rho}\mathcal X\rho\dd\mu_{N} =-\int_N\Phi\cdot \mathcal X\rho\dd\mu_{N}
		\le C\epsilon T_{\mathrm{far}}+C_\epsilon(\mathfrak e_a(N)+\mathcal H_a).
		\]
		The equality is Gaussian first variation; the outer cutoff is
		removed using boundedness of $\mathcal X$,
		\eqref{eq:gaussian-moment-bound}, and $G(N)<\infty$.
		With $C\epsilon\le c/2$, we obtain
		$(c/2)T_{\mathrm{far}}\le C_\epsilon(\mathfrak e_a(N)+\mathcal H_a)$.
		For
		$r\ge r_c$, the integral of $\mathfrak t$ is bounded by $C\mathfrak e_a(N)$, using
		\eqref{eq:cyl-weighted-caccioppoli} on $\mathcal T_{a,\varepsilon}$ and
		$\mathfrak t\le1\le C j_a^2$ on $\{r>0\}\setminus\mathcal T_{a,\varepsilon}$. This proves
		\eqref{eq:cyl-tilt-local-to-global}.
		
		Choose a fixed $\delta_0<2r_c$ in the range of the relative
		tube estimate. The integral of $j_{a_i}^2/r^2$ on
		$\delta_0\le r<2r_c$ is bounded by
		$\delta_0^{-2}\int_I D_{a_i}(N_s^i)^2\dd s$. On the smaller annuli, put $\delta_j:=2^{-j}\delta_0$.
		The assumed relative tube bound gives
		\begin{equation}\label{eq:cyl-relative-hardy}
			\begin{aligned}
				\int_I\int_{N_s^i\cap\{r<\delta_0,\ r_y<L_{\rm ax}+1\}}
				r^{-2}j_{a_i}^2\rho\dd\mu_{N_s^i}\dd s
				&\le\sum_{j\ge0}\frac4{\delta_j^2}
				\int_I\int_{N_s^i\cap\{r<\delta_j,\ r_y<L_{\rm ax}+1\}}
				j_{a_i}^2\rho\dd\mu_{N_s^i}\dd s\\
				&\le C\eta_i^2\sum_{j\ge0}\delta_j^{\sigma-2}
				=\frac{C\delta_0^{\sigma-2}}{1-2^{-(\sigma-2)}}\eta_i^2.
			\end{aligned}
		\end{equation}
		This sum is finite precisely because $\sigma>2$.
		By \eqref{eq:cyl-tilt-local-to-global},
		\eqref{eq:cyl-relative-hardy}, and the assumed squared-distance bound,
		\[
		\begin{aligned}
			\int_I\int_{N_s^i}|\mathbf n_y|^2\rho\dd\mu_{N_s^i}\dd s
			&\le C\eta_i^2+C\delta_0^{-2}\int_I D_{a_i}(N_s^i)^2\dd s
			+\frac{C\delta_0^{\sigma-2}}{1-2^{-(\sigma-2)}}\eta_i^2\\
			&\le C\eta_i^2.
		\end{aligned}
		\]
		This proves the final assertion.
	\end{proof}
	
	The following result supplies the cylindrical counterpart of the relative Gaussian tail control in \cite[Lemma~2.9]{FH1}. The axial-tilt bound is an additional input needed to control the full ambient second moment.
	\begin{lem}\label{lem:cyl-gaussian-tightness}
		Let $I:=[s_1,s_2]\Subset J$ be fixed bounded intervals with
		$s_1<s_2$. Fix a reference $\barC_z$ on $J$, and let $W_z$ be
		its degree-two squared-distance function, comparable to
		$j_z^2$ from \eqref{eq:reference-distance} and satisfying
		$|DW_z|^2\le CW_z$. Align the axis of $\barC_z$ with
		$\{0\}\times\R^k$ by a fixed rotation, and let $q(a)$ be
		its link residual. Suppose that for a rescaled flow $(N_s)_{s\in J}$,
		\begin{equation}\label{eq:cyl-moment-input}
			\int_I D_z(N_s)^2\dd s+\int_J G(N_s)^2\dd s
			+|q(a)|^2
			+\int_I\int_{N_s}|\mathbf n_y|^2\rho\dd\mu_{N_s}\dd s
			\le C_0\eta^2 .
		\end{equation}
		Then
		\begin{equation}\label{eq:cyl-relative-tail}
			\int_I\int_{N_s}|Z|^2W_z\rho\dd\mu_{N_s}\dd s\le C\eta^2,\qquad
			\int_I\int_{N_s\setminus B_R}W_z\rho\dd\mu_{N_s}\dd s
			\le CR^{-2}\eta^2\quad(R\ge1).
		\end{equation}
		The constant depends on the area bound, the intervals, the fixed
		reference neighborhood, and $C_0$.
	\end{lem}
	\begin{proof}
		Write $p(s):=\int_{N_s}W_z\rho\dd\mu_{N_s}$ and
		$\Delta:=\int_J G(N_s)^2\dd s$. We first bound $p$ at one
		time before $I$. Choose $s'\in I$ with
		$D_z(N_{s'})^2\le C\eta^2$. In a fixed subinterval of
		$J$ to the left of $I$, choose $s_0$ with
		$G(N_{s_0})^2\le C\Delta$. Set $j(Z):=\dist(Z,\barC_z)$,
		$f:=j/\sqrt{1+j^2}$, and $b(s):=\int_{N_s}f^2\rho\dd\mu_{N_s}$.
		Here $0\le f\le1$ and $f$ has Lipschitz constant at most one.
		Gaussian first variation gives
		\[
		-2\sqrt b\,G-G^2\le b'\le 2\sqrt b\,G
		\]
		in integrated form. Indeed, mollify $f$ in the ambient space,
		preserving both bounds, apply first variation with an outer cutoff,
		and remove the cutoff using the Gaussian area bound and the finite
		time integral of $G^2$. Uniform convergence of the mollifications
		gives the displayed inequalities for $b$. If
		$M:=\sup_{[s_0,s']}b$, integration from a time of maximum to $s'$
		gives $M\le b(s')+2\sqrt{M|J|\Delta}+\Delta$. Young's inequality
		and $b(s')\le D_z(N_{s'})^2$ give $b(s_0)\le C\eta^2$.
		
		To recover the full distance on $N=N_{s_0}$, put $v:=f^2$.
		Gaussian first variation with the vector field $v^2 Z$ gives
		\[
		\tfrac12\int_N |Z|^2v^2\rho\dd\mu_N
		=\int_N\bigl(nv^2+2v\langle Z^\top,\nabla_N v\rangle+
		v^2\langle\Phi,Z\rangle\bigr)\rho\dd\mu_N .
		\]
		On each compact patch we approximate the restriction of $v$
		in $W^{1,2}(N)$ by uniformly bounded smooth functions. We then
		remove the outer cutoff using the Gaussian moment bounds and $G(N)<\infty$.
		Young's inequality bounds
		the left-hand integral by a constant times
		$\int_N(v^2+|\nabla_N v|^2+v^2|\Phi|^2)\rho\dd\mu_N$.
		Since $v^2\le v$, $|\nabla_N v|^2\le4v$, and
		$j^2=v+(1+j^2)v^2\le2v+|Z|^2v^2$, this proves
		$D_z(N_{s_0})^2\le C(b(s_0)+G(N_{s_0})^2)\le C\eta^2$.
		Comparability of $W_z$ and $j^2$ now gives $p(s_0)\le C\eta^2$.
		
		We next estimate $\int_I\int_{N_s}W_z|\Phi|^2\rho\dd\mu_{N_s}\dd s$
		by differentiating the Gaussian integral of a truncated positional function. With an ambient cutoff $\chi_R$ equal to one on
		$B_R$ and supported on $B_{2R}$, with $|D\chi_R|\le C/R$,
		put $W_{z,R}:=\chi_R^2W_z$ and $p_R(s):=\int_{N_s}W_{z,R}\rho\dd\mu_{N_s}$.
		Since $0\le\chi_R\le1$, $W_z\le C|Z|^2$, and
		$|DW_z|^2\le CW_z$, the product rule gives
		\begin{equation}\label{eq:cyl-truncated-position-gradient}
			|DW_{z,R}|^2\le8\chi_R^2W_z^2|D\chi_R|^2
			+2\chi_R^4|DW_z|^2\le C\chi_R^2W_z=CW_{z,R}.
		\end{equation}
		The constant is independent of $R$. Gaussian first variation
		\cite[Section~2.1]{FH1}, applied to the time-independent function
		$W_{z,R}$, and the Cauchy--Schwarz inequality now give
		\begin{equation}\label{eq:cyl-truncated-position-transport}
			\begin{aligned}
				p_R'(s)+\int_{N_s}W_{z,R}|\Phi|^2\rho\dd\mu_{N_s}
				&=\int_{N_s}DW_{z,R}\cdot\Phi\rho\dd\mu_{N_s}\\
				&\le\left(\int_{N_s}|DW_{z,R}|^2\rho\dd\mu_{N_s}\right)^{1/2}
				\left(\int_{N_s}|\Phi|^2\rho\dd\mu_{N_s}\right)^{1/2}\\
				&\le C\sqrt{p_R(s)}\,G(N_s),
			\end{aligned}
		\end{equation}
		where the last line uses \eqref{eq:cyl-truncated-position-gradient}.
		The squared-distance function is $C^{1,1}$ across the axis, so it is an admissible
		first-variation test by smooth approximation. We integrate the
		upper differential inequality for $\sqrt{p_R+\varepsilon}$
		and let $\varepsilon\downarrow0$ to obtain
		\[
		\sup_{[s_0,s_2]}\sqrt{p_R(s)}
		\le\sqrt{p_R(s_0)}+C\int_{s_0}^{s_2}G(N_s)\dd s
		\le C\eta .
		\]
		Integrating \eqref{eq:cyl-truncated-position-transport} from
		$s_0$ to $s_2$ and using $p_R(s_2)\ge0$, we obtain
		\[
		\begin{aligned}
			\int_{s_0}^{s_2}\int_{N_s}W_{z,R}|\Phi|^2\rho\dd\mu_{N_s}\dd s
			&\le p_R(s_0)-p_R(s_2)
			+C\sup_{[s_0,s_2]}\sqrt{p_R}\int_{s_0}^{s_2}G(N_s)\dd s\\
			&\le C\eta^2+C\eta\sqrt{s_2-s_0}\sqrt\Delta\le C\eta^2.
		\end{aligned}
		\]
		Fatou's lemma as $R\to\infty$, followed by restriction to $I$, gives
		\begin{equation}\label{eq:cyl-distance-function-weighted-dissipation}
			\int_I\int_{N_s}W_z|\Phi|^2\rho\dd\mu_{N_s}\dd s\le C\eta^2 .
		\end{equation}
		
		At almost every time $s\in I$, both $G(N_s)^2$ and
		$\int_{N_s}W_z|\Phi|^2\rho\dd\mu_{N_s}$ are finite by
		\eqref{eq:cyl-moment-input} and
		\eqref{eq:cyl-distance-function-weighted-dissipation}.
		Write $\mathfrak a:=x\cdot\mathbf n_x$, $\mathfrak b:=y\cdot\mathbf n_y$, and
		$\mathfrak t:=|\mathbf n_y|^2$. Since $Z\cdot DW_z=2W_z$ by homogeneity
		and $n:=m+k$ is the dimension of $N_s$, by the first variation formula and $\operatorname{div}_{N_s}(W_z Z)
		=nW_z+DW_z\cdot Z^\top$, we get
		\begin{equation}\label{eq:cyl-distance-function-moment-identity}
			\frac12\int_{N_s}|Z|^2W_z\rho\dd\mu_{N_s}
			=(n+2)p(s)-\int_{N_s}(DW_z\cdot\mathbf n)(\mathfrak a+\mathfrak b)\rho\dd\mu_{N_s}
			+\int_{N_s}W_z\Phi\cdot Z\rho\dd\mu_{N_s} .
		\end{equation}
		
		For every $\epsilon>0$, Young's inequality gives
		\[
		|DW_z\cdot\mathbf n|\,|\mathfrak b|
		\le C\sqrt{W_z}\,|y|\sqrt{\mathfrak t}
		\le \epsilon |y|^2W_z+C_\epsilon \mathfrak t ,
		\qquad
		W_z|\Phi\cdot Z|
		\le \epsilon |Z|^2W_z+C_\epsilon W_z|\Phi|^2 .
		\]
		For the transverse term, fix $\varepsilon>0$ sufficiently small
		that the signed normal coordinates from
		Lemma~\ref{lem:distance-function} are defined on
		\begin{equation}\label{eq:cyl-moment-normal-tube}
			\mathcal T_{z,\varepsilon}:=\{Z=(x,y):r>0,\ j_z(Z)<\varepsilon r\},
		\end{equation}
		and that \eqref{eq:cyl-weighted-caccioppoli} applies there after
		alignment of the reference axis. On $\mathcal T_{z,\varepsilon}$,
		let $l_z$ be signed distance to $\barC_z$, and put
		$\vartheta_z^2:=1-(\mathbf n\cdot D l_z)^2$, as in
		Lemma~\ref{lem:weighted-normal-angle}. Then $|l_z|=j_z$, $W_z=l_z^2$,
		and $|\mathfrak a|\le |l_z|+Cr\vartheta_z$. If $r\le1$,
		we get $\sqrt{W_z}|\mathfrak a|\le C(W_z+r^2\vartheta_z^2)$. If $r>1$,
		we get $\sqrt{W_z}|\mathfrak a|\le \epsilon r^2W_z+
		C_\epsilon(W_z+\vartheta_z^2)$.
		On $\{r>0\}\setminus\mathcal T_{z,\varepsilon}$,
		\eqref{eq:cyl-moment-normal-tube} gives $j_z\ge\varepsilon r$.
		Thus $|\mathfrak a|\le r\le\varepsilon^{-1}j_z\le C\sqrt{W_z}$
		by Lemma~\ref{lem:distance-function}; on the axis $\mathfrak a=0$.
		Thus \eqref{eq:cyl-weighted-caccioppoli} gives
		\begin{equation}\label{eq:cyl-transverse-moment-bound}
			\int_{N_s}|DW_z\cdot\mathbf n|\,|\mathfrak a|\rho\dd\mu_{N_s}
			\le \epsilon\int_{N_s}r^2W_z\rho\dd\mu_{N_s}
			+C_\epsilon\bigl(D_z(N_s)^2+G(N_s)^2+|q(a)|^2\bigr).
		\end{equation}
		Indeed, the weight $w:=r^2/(1+r^2)$ in
		\eqref{eq:cyl-weighted-caccioppoli} satisfies
		\[
		\tfrac12r^2\le w\le r^2\quad(0<r\le1),\qquad
		\tfrac12\le w\le1\quad(r>1).
		\]
		Consequently, \eqref{eq:cyl-weighted-caccioppoli} bounds the
		integrals of $r^2\vartheta_z^2$ on
		$N_s\cap\mathcal T_{z,\varepsilon}\cap\{r\le1\}$ and of
		$\vartheta_z^2$ on $N_s\cap\mathcal T_{z,\varepsilon}\cap\{r>1\}$.
		Together with $\int_{N_s}W_z\rho\dd\mu_{N_s}\le CD_z(N_s)^2$,
		this proves \eqref{eq:cyl-transverse-moment-bound}.
		
		Put $M_z(s):=\int_{N_s}|Z|^2W_z\rho\dd\mu_{N_s}$. Substitution in
		\eqref{eq:cyl-distance-function-moment-identity} gives
		\[
		\tfrac12M_z(s)
		\le 3\epsilon M_z(s)+C_\epsilon\bigl(D_z(N_s)^2+G(N_s)^2
		+|q(a)|^2\bigr)+C_\epsilon\int_{N_s}(\mathfrak t+W_z|\Phi|^2)\rho\dd\mu_{N_s}.
		\]
		Take $\epsilon\le1/12$, so that
		$\tfrac12M_z(s)-3\epsilon M_z(s)\ge\tfrac14M_z(s)$,
		and move $3\epsilon M_z(s)$ to the left.
		We conclude that
		\begin{equation}\label{eq:cyl-relative-moment-slice}
			M_z(s)\le C\left(D_z(N_s)^2+G(N_s)^2+|q(a)|^2
			+\int_{N_s}(\mathfrak t+W_z|\Phi|^2)\rho\dd\mu_{N_s}\right).
		\end{equation}
		Equations \eqref{eq:cyl-relative-moment-slice},
		\eqref{eq:cyl-moment-input}, and
		\eqref{eq:cyl-distance-function-weighted-dissipation} give
		\[
		\begin{aligned}
			\int_I M_z(s)\dd s
			&\le C\eta^2+C\int_I\int_{N_s}W_z|\Phi|^2\rho\dd\mu_{N_s}\dd s
			\le C\eta^2,\\
			\int_I\int_{N_s\setminus B_R}W_z\rho\dd\mu_{N_s}\dd s
			&\le R^{-2}\int_I M_z(s)\dd s\le CR^{-2}\eta^2.
		\end{aligned}
		\]
		This proves \eqref{eq:cyl-relative-tail}.
	\end{proof}
	
	\section{Uniqueness, \L{}ojasiewicz inequality and rigidity}\label{sec:main-proofs}
	
	\subsection{Fitting a time-independent reference and the local discrete inequality}\label{sec:discrete}
	Write $\Theta:=F(\barC)$ and $D_0(N):=D_{z=0}(N)$, as in
	\eqref{eq:reference-distance}.
	
	\begin{lem}\label{lem:cyl-distance-entry}
		Let $(N_s^i)_{s\in J}$ be smooth rescaled flows on a fixed compact
		interval $J:=[s_1,s_2]$, with $s_1<s_2$, satisfying the uniform area-ratio bound.
		Suppose
		\begin{equation}\label{eq:cyl-distance-entry-input}
			\sup_{s\in J}\bigl(D_0(N_s^i)+|F(N_s^i)-\Theta|\bigr)\longrightarrow0.
		\end{equation}
		Then the area measures $\mu_{N_s^i}$ converge locally weakly as
		Radon measures to $\mu_{\barC}$, uniformly for $s\in J$. That is, for every
		$\varphi\in C_c(\R^{n+1})$,
		\begin{equation}\label{eq:cyl-uniform-radon-convergence}
			\sup_{s\in J}\left|\int_{N_s^i}\varphi\dd\mu_{N_s^i}
			-\int_{\barC}\varphi\dd\mu_{\barC}\right|
			\longrightarrow0.
		\end{equation}
		On compact subsets of $\operatorname{Reg}(\barC)\times(s_1,s_2)$,
		the convergence is smooth and single-sheeted. Moreover, on every fixed
		relatively compact open set $\mathcal U$ whose closure is
		disjoint from the axis, for $s$ in a compact subinterval of $(s_1,s_2)$ and all large $i$,
		$N_s^i\cap\mathcal U=\operatorname{graph}(u_i(\cdot,s))\cap\mathcal U$,
		where the graph is defined on a regular reference patch containing
		$\overline{\mathcal U}\cap\barC$ in its interior and tends smoothly to zero.
	\end{lem}
	\begin{proof}
		By the monotonicity formula and \eqref{eq:cyl-distance-entry-input}, we have
		\begin{equation}
			\Delta_i:=\int_{s_1}^{s_2}G(N_s^i)^2\dd s
			=F(N_{s_1}^i)-F(N_{s_2}^i)
			\le2\sup_{s\in J}|F(N_s^i)-\Theta|\longrightarrow0.
		\end{equation}
		Choose $t_i\in J$ with $G(N_{t_i}^i)^2\le\Delta_i/(s_2-s_1)$.
		The relation $\mathbf H=\Phi-Z^\perp/2$ and the area bound give
		uniform local $L^2$ bounds for the mean curvature at $t_i$.
		Integral-varifold compactness \cite{Allard} therefore gives a
		subsequential integral limit $V$. By \eqref{eq:cyl-distance-entry-input},
		$\operatorname{spt}\|V\|\subset\barC$.
		The area bound gives zero mass to the axis by
		\eqref{eq:raw-tube-area}. Almost every tangent plane of $V$ is
		therefore the tangent plane of the regular cylinder, so $Z^\perp=0$
		almost everywhere. We pass to the limit in the Gaussian first variation formula
		and use $G(N_{t_i}^i)\to0$ to conclude that the Gaussian first variation of the limit
		vanishes. Since $Z^\perp=0$ almost everywhere, this is equivalent
		to ordinary stationarity of $V$.
		The constancy theorem \cite{SimonGMT} on the connected regular cylinder
		now gives $V=j|\barC|$ for an integer $j\ge0$.
		By the uniform Gaussian tails in \eqref{eq:gaussian-moment-bound}, we have
		\[
		j\Theta=F(V)=\lim_iF(N_{t_i}^i)=\Theta.
		\]
		Thus $j=1$, which identifies the limit at the selected times.
		
		For $\varphi\in C_c^\infty(\R^{n+1})$, Gaussian first variation and
		the Cauchy--Schwarz inequality give
		\begin{equation}
			\sup_{s\in J}\left|
			\int_{N_s^i}\varphi\rho\dd\mu_{N_s^i}
			-\int_{N_{t_i}^i}\varphi\rho\dd\mu_{N_{t_i}^i}\right|
			\le C_\varphi\bigl((s_2-s_1)^{1/2}\Delta_i^{1/2}+\Delta_i\bigr)
			\longrightarrow0.
		\end{equation}
		Consequently $\mu_{N_s^i}$ converge locally weakly as Radon measures
		to $\mu_{\barC}$, uniformly for $s\in J$, as stated in
		\eqref{eq:cyl-uniform-radon-convergence}.
		Every subsequence has this same limit, so the assertion holds for
		the full sequence. We can then apply Brakke's local regularity theorem
		\cite[Section~6.11]{Brakke}, in the local smooth-convergence form of
		\cite[Lemma~2.1]{Schulze}, to the corresponding ordinary flows on
		compact regular spacetime sets in the interior of $J$.
		This gives smooth single-sheeted convergence on regular patches.
		If further points persisted in a compact complement of those patches
		at positive distance from the cylinder, a subsequence of their
		centers and times would converge. Their backward Gaussian integrals
		are at least one by monotonicity, while at a sufficiently small
		fixed backward time the limiting integral is less than one by
		positive distance, local convergence, and the uniform Gaussian
		tails. This excludes such points and proves the asserted equality
		in that open set.
	\end{proof}
	
	Recall the squared positional function $W_z$ from
	\eqref{eq:positional-function-definition} and the reduced link residual
	$q(a)$ from \eqref{eq:reduced-link-data}.
	\begin{prop}\label{prop:joint-full}
		Assume the cone hypotheses of Theorem~\ref{thm:main},
		and fix the area-ratio constant $\Lambda$ in \eqref{equ:areabd}.
		There are constants $\epsilon_0>0$ and $C<\infty$ with the
		following property. Suppose a smooth rescaled flow with this area-ratio bound on $J:=[t,t+2]$ satisfies
		\begin{equation}\label{eq:cyl-fitted-smallness}
			\sup_{s\in J}\bigl(D_0(N_s)+|F(N_s)-\Theta|\bigr)<\epsilon_0.
		\end{equation}
		For $I:=[t+\tfrac12,t+\tfrac32]$,
		\begin{equation}\label{eq:cyl-fitted-interval}
			\min_{z\in\mathcal Z}\left\{\frac1{|I|}\int_I\int_{N_s}W_z\rho\dd\mu_{N_s}\dd s
			+|q(a)|^2\right\}
			\le C\int_J G(N_s)^2\dd s.
		\end{equation}
		The minimum is taken over
		\begin{equation}\label{eq:cyl-fitting-parameter-ball}
			\mathcal Z:=\{z=(a,\mathsf B)\in K\times\mathscr E:
			|a|^2+|\mathsf B|^2\le\varepsilon_{\rm par}^2\},
		\end{equation}
		where $K$ and $\mathscr E$ are the spaces in \eqref{eq:augmented-family},
		$|a|$ is the fixed $L^2(\Sigma)$ norm, and $|\mathsf B|$ is the
		Euclidean Hilbert--Schmidt norm. The radius $\varepsilon_{\rm par}>0$
		is fixed sufficiently small that the augmented family is defined
		and locally injective on a neighborhood of $\mathcal Z$.
	\end{prop}
	\begin{proof}
		\noindent\textbf{Step 1. The minimizing parameters and local position bounds.}
		If the assertion fails, take flows with the left side of
		\eqref{eq:cyl-fitted-smallness} tending to zero and with a positive
		fitted minimum exceeding $i\int_JG(N_s^i)^2\dd s$ in
		\eqref{eq:cyl-fitted-interval}. Shift time so that $J=[0,2]$.
		Lemma~\ref{lem:cyl-distance-entry} shows that $\mu_{N_s^i}$ converge
		locally weakly as Radon measures to $\mu_{\barC}$, uniformly for $s\in J$. Choose a minimizing parameter
		$z_i=(a_i,\mathsf B_i)\in\mathcal Z$, which exists by continuity of
		the objective on the compact set \eqref{eq:cyl-fitting-parameter-ball},
		and let $\eta_i^2>0$ be the minimum value. Evaluation at $z=0$ and
		$W_0\le C\operatorname{dist}(\,\cdot\,,\barC)^2$ give $\eta_i\to0$. To identify the limit of the minimizing parameters, pass to a subsequence
		with $z_i\to z_*\in\mathcal Z$. Choose a nonnegative ambient cutoff $\chi$
		that equals one near the support of the cutoff in \eqref{eq:parameter-gram-matrix}
		and has compact support away from the axis. Since $W_{z_i}\to W_{z_*}$
		uniformly on $\operatorname{spt}\chi$, the uniform measure convergence
		\eqref{eq:cyl-uniform-radon-convergence} and the definition of $\eta_i$ give
		\[
		\int_{\barC}\chi W_{z_*}\rho\dd\mu_{\barC}
		=\lim_{i\to\infty}\frac1{|I|}\int_I\int_{N_s^i}\chi W_{z_i}\rho\dd\mu_{N_s^i}\dd s
		\le\lim_{i\to\infty}\|\chi\|_\infty\eta_i^2=0 .
		\]
		Thus $W_{z_*}$ vanishes on the reference patch where $\chi=1$.
		Since its zero set is $\barC_{z_*}$, the normal graph height $h_{z_*}$
		vanishes on the support of the cutoff defining $\mathcal F$.
		Hence $\mathcal F(z_*)=\mathcal F(0)$. Local injectivity, established using \eqref{eq:parameter-gram-matrix}, gives $z_*=0$.
		Every subsequential limit is therefore zero, so $z_i\to0$ and the
		minimizers lie in the interior of $\mathcal Z$ for all large $i$. Write
		\[
		\mathfrak q_i:=|q(a_i)|\le\eta_i,\qquad
		\Delta_i:=\int_J G(N_s^i)^2\dd s=o(\eta_i^2).
		\]
		Apply the inverse of the fixed fitted rotation to the flow and to
		its entire comparison family. More precisely, put
		\[
		\widetilde N_s^i:=R_{\mathsf B_i}^{-1}N_s^i,\qquad
		\widetilde{\barC}_z^{\,i}:=R_{\mathsf B_i}^{-1}\barC_z,\qquad
		\widetilde W_z^{\,i}(X):=W_z(R_{\mathsf B_i}X).
		\]
		At the minimizing parameter these definitions give
		$\widetilde{\barC}_{z_i}^{\,i}=\barC_{a_i}$ and
		$\widetilde W_{z_i}^{\,i}=W_{(a_i,0)}$.
		Orthogonal invariance of $\rho$ and the induced measure preserves
		the fitting objective and its minimizing parameter. In differentiating with respect to $z$, hold $R_{\mathsf B_i}$ fixed. Thus
		\[
		D_z\widetilde W_z^{\,i}(X)=D_zW_z(R_{\mathsf B_i}X).
		\]
		Since $R_{\mathsf B_i}\to\mathrm{Id}$, the limiting parameter normal
		velocities remain those in
		\eqref{eq:parameter-normal-velocities} at $z=0$.
		For the rest of this proof we suppress the tildes and the additional
		index $i$ on the transported comparison family, and write
		$\barC_i:=\barC_{a_i}$.
		
		For a fixed compact cutoff $\chi$, put
		$\mathcal P_i^\chi(s):=\int_{N_s^i}\chi^2W_{z_i}\rho\dd\mu_{N_s^i}$. On its support the
		squared-distance function is bounded and $|D(\chi^2W_{z_i})|^2\le C_\chi\chi^2W_{z_i}$.
		The rescaled first-variation formula \cite[Section~2.1]{FH1}
		therefore gives
		\[
		|(\mathcal P_i^\chi)'|\le C_\chi\sqrt{\mathcal P_i^\chi}G(N_s^i)
		+C_\chi G(N_s^i)^2.
		\]
		There is a time $s_i^\chi\in I$ with
		$\mathcal P_i^\chi(s_i^\chi)\le C\eta_i^2$.
		Define $M_i^\chi:=\sup_{s\in J}\mathcal P_i^\chi(s)$.
		By integrating the differential inequality from $s_i^\chi$ to any
		$s\in J$ in either time direction and applying the Cauchy--Schwarz
		inequality in time, we obtain
		\begin{equation}
			\begin{aligned}
				\mathcal P_i^\chi(s)
				&\le C_\chi\eta_i^2+C_\chi\sqrt{M_i^\chi}
				\int_JG(N_\tau^i)\dd\tau+C_\chi\Delta_i\\
				&\le C_\chi\eta_i^2+C_\chi\sqrt{M_i^\chi}
				\left(|J|\int_JG(N_\tau^i)^2\dd\tau\right)^{1/2}
				+C_\chi\Delta_i\\
				&\le\tfrac12M_i^\chi+C_\chi(\eta_i^2+\Delta_i).
			\end{aligned}
		\end{equation}
		The last line uses Young's inequality and the fixed length of $J$.
		We take the supremum in $s$, absorb $M_i^\chi/2$, and use
		$\Delta_i=o(\eta_i^2)$ to obtain
		\begin{equation}
			\sup_{s\in J}\int_{N_s^i}\chi^2W_{z_i}\rho\dd\mu_{N_s^i}
			\le C_\chi\eta_i^2.
		\end{equation}
		By \eqref{eq:reduced-cone-curvature}, the reference curvature is bounded
		by $C|q(a_i)|/r\le C\eta_i/r$. Lemmas~\ref{lem:cyl-tube},
		\ref{lem:axial-tilt-estimate}, and \ref{lem:cyl-gaussian-tightness},
		with \eqref{eq:cyl-relative-tube}, \eqref{eq:cyl-tilt-local-to-global},
		and \eqref{eq:cyl-relative-tail}, give, for some $\sigma>2$, every fixed
		$R\ge1$, and sufficiently small $\delta>0$,
		\begin{equation}
			\begin{aligned}
				\frac1{|I|\eta_i^2}\int_I\int_{N_s^i}|\mathbf n_y|^2
				\rho\dd\mu_{N_s^i}\dd s&\le C,\\
				\limsup_{i\to\infty}\frac1{|I|\eta_i^2}
				\int_I\int_{N_s^i\cap\{r<\delta,\ |Z|<R\}}
				W_{z_i}\rho\dd\mu_{N_s^i}\dd s&\le C_R\delta^\sigma,\\
				\limsup_{i\to\infty}\frac1{|I|\eta_i^2}
				\int_I\int_{N_s^i\setminus B_R}
				W_{z_i}\rho\dd\mu_{N_s^i}\dd s&\le CR^{-2}.
			\end{aligned}
		\end{equation}
		
		\medskip
		\noindent\textbf{Step 2. The normalized graph limit and its mass.}
		On every fixed compact regular patch, we apply
		Lemma~\ref{lem:cyl-distance-entry}.
		The area measures of the ordinary flows converge locally weakly
		as Radon measures to $\mu_{\barC}$, uniformly on compact time
		intervals, because $\barC$ is dilation invariant. It follows
		that on compact subintervals of the interior of $J$ the rescaled flows
		are complete single graphs over the reduced references, with
		heights $u_i$ tending smoothly to zero. The graph equation for
		$u_i/\eta_i$ has uniformly parabolic coefficients converging
		smoothly to those of $\partial_s-L_{\barC}$. The inhomogeneous term
		is the reference curvature divided by $\eta_i$, which converges
		on compact regular sets to $-r^{-1}\xi$ after choosing a
		subsequence $q(a_i)/\eta_i\to\xi\in K$.
		The normalized spacetime $L^2$ bound on each compact patch
		comes from the compact transport estimate and comparability
		of graph measures. Interior estimates give smooth convergence
		on smaller regular cylinders. The normal velocity relation
		$\partial_s u_i\langle \mathbf n_{\barC_i},\mathbf n_{N_s^i}\rangle
		=\langle\Phi_i,\mathbf n_{N_s^i}\rangle$ and
		$\Delta_i/\eta_i^2\to0$ show that the limit $U$ is independent
		of time. Therefore we obtain
		\begin{equation}\label{eq:cyl-forced-limit}
			L_{\barC} U=r^{-1}\xi,\qquad
			U\in L^2_\rho(\barC),\qquad
			\|U\|_{L^2_\rho}^2+|\xi|^2=1 .
		\end{equation}
		We verify the normalization in \eqref{eq:cyl-forced-limit}. For $R>1$ and $\delta>0$,
		let $\mathcal K_{\delta,R}:=\{|Z|<R,\ r>\delta\}\cap \barC$. On this set the
		squared-distance function equals squared graph height for large $i$, so smooth graph
		convergence gives convergence of its normalized spacetime integral
		to $\int_{\mathcal K_{\delta,R}}U^2\rho\dd\mu_{\barC}$. On the complement,
		\eqref{eq:cyl-relative-tube} and \eqref{eq:cyl-relative-tail} give
		\[
		\limsup_{i\to\infty}\frac1{|I|\eta_i^2}
		\int_I\int_{N_s^i\cap(\{r<\delta,\ |Z|<R\}\cup\{|Z|>R\})}
		W_{z_i}\rho\dd\mu_{N_s^i}\dd s
		\le C_R\delta^\sigma+CR^{-2}.
		\]
		First let $\delta\downarrow0$ at fixed $R$, and then let $R\to\infty$.
		Since the objective divided by $\eta_i^2$ equals one and
		$q(a_i)/\eta_i\to\xi$, this proves the last equality in
		\eqref{eq:cyl-forced-limit}. Also, we have
		\[
		\|r^{-1}\xi\|_{L^2_\rho(\barC)}^2
		=c_m\left(\int_0^\infty r^{m-3}e^{-r^2/4}\dd r\right)
		\|\xi\|_{L^2(\Sigma)}^2<\infty,
		\]
		since $m>2$.
		
		\medskip
		\noindent\textbf{Step 3. Gaussian moments of the limit.}
		The relative tube estimate \eqref{eq:cyl-relative-tube} passes to
		$U$ on each fixed bounded region. Indeed, first restrict to
		$\delta'<r<\delta$ and a compact regular spacetime set in the
		middle interval, where the normalized graphs converge smoothly.
		Then exhaust that interval and let $\delta'\downarrow0$.
		Time independence gives
		\[
		\int_{\barC\cap\{r<\delta,\ |Z|<R\}}|U|^2\rho\dd\mu_{\barC}
		\le C_R\delta^\sigma,\qquad \sigma>2 .
		\]
		By \eqref{eq:cyl-forced-limit}, $L_{\barC}U=r^{-1}\xi$ on the
		regular cylinder and $U\in L^2_\rho$. The right hand side is
		square integrable by \eqref{eq:residual-square-integrability}.
		Thus Lemma~\ref{lem:all-root-moment-test} applies with $p=\sigma$.
		By testing the equation against $r\xi$, we obtain $\xi=0$. Proposition~\ref{prop:gaussian-bordered}
		then gives $U\in\mathcal G_{\barC}$ under the angular hypotheses
		of the theorem.
		
		\medskip
		\noindent\textbf{Step 4. The fitting equation and the contradiction.}
		The interior minimizing parameter satisfies
		\begin{equation}\label{eq:cyl-fit-euler-equation}
			\frac1{|I|}\int_I\int_{N_s^i}D_zW_{z_i}\rho\dd\mu_{N_s^i}\dd s
			+2D_zq(z_i)^\ast q(a_i)=0.
		\end{equation}
		Differentiation under the integral follows from the uniform
		parameter derivative bound and the Gaussian moment bound.
		On compact regular sets Lemma~\ref{lem:distance-function}
		expresses the parameter derivative of the squared distance by
		the graph height and the normal parameter velocity. After division
		by $\eta_i$ its quadratic error tends to zero, since
		$\int_I\int_{N_s^i} W_{z_i}\rho\dd\mu_{N_s^i}\dd s
		\le C\eta_i^2$.
		The uniform bound $|D_zW_z|\le C|Z|\sqrt{W_z}$ gives
		\begin{equation}
			\begin{aligned}
				\frac1{\eta_i}\int_I\int_{N_s^i\cap\{|Z|>R\}}
				|D_zW_{z_i}|\rho\dd\mu_{N_s^i}\dd s
				&\le C\left(\frac1{\eta_i^2}\int_I\int_{N_s^i}
				W_{z_i}\rho\dd\mu_{N_s^i}\dd s\right)^{1/2}\\
				&\quad{}\times\left(\int_I\int_{N_s^i\cap\{|Z|>R\}}
				|Z|^2\rho\dd\mu_{N_s^i}\dd s\right)^{1/2}.
			\end{aligned}
		\end{equation}
		The last factor tends to zero uniformly in $i$ by
		\eqref{eq:gaussian-moment-bound}. On $\{|Z|<R,r<\delta\}$ the
		bound from the same Cauchy--Schwarz inequality is $O_R(\delta^{m/2})$ by
		\eqref{eq:raw-tube-area}. Finally,
		\[
		\eta_i^{-1}|D_zq(z_i)^\ast q(a_i)|
		\le\|D_zq(z_i)\|\,|q(a_i)|/\eta_i
		\le\|D_zq(z_i)\|\longrightarrow0 .
		\]
		Here $D_zq(0)=0$ by Proposition~\ref{prop:cyl-link-reduction}
		and the independence of $q$ from
		$\mathsf B$. Divide \eqref{eq:cyl-fit-euler-equation} by
		$\eta_i$, pass to the limit on regular compact sets, then let
		$\delta\downarrow0$ and $R\to\infty$. Each limiting normal
		parameter velocity $g_{0,j}$ satisfies
		\[
		\int_{\barC} U g_{0,j}\rho\dd\mu_{\barC}=0.
		\]
		By \eqref{eq:parameter-normal-velocities} these velocities
		span $\mathcal G_{\barC}$. Since $U$ belongs to that space,
		$U=0$. Together with $\xi=0$ this contradicts the normalization
		in \eqref{eq:cyl-forced-limit} and proves the interval estimate.
	\end{proof}
	
	\subsection{Continuation and uniqueness}\label{sec:uniqueness}
	
	\begin{proof}[Proof of Theorems~\ref{thm:main} and \ref{thm:energy-loja}]
		By the prescribed tangent subsequence and the uniform area bound,
		Gaussian area converges to $F(\barC)$ along a sequence of times going to
		infinity. Huisken's monotonicity formula \cite{Huisken} implies $E\ge0$ and $E\to0$, as
		well as $E'(s)=-G(N_s)^2$. Choose $r_D,\epsilon_E>0$ with
		$r_D+\epsilon_E<\epsilon_0$, where $\epsilon_0$ is from
		Proposition~\ref{prop:joint-full}. By taking $r_D$ and $\epsilon_E$ sufficiently small, we ensure that
		the complete graphs required by Proposition~\ref{prop:anchorarea}
		exist on the middle interval whenever $D_0(N_\tau)<r_D$ and
		$0\le E(\tau)<\epsilon_E$ for every $\tau\in[s,s+2]$.
		This follows by contradiction from Lemma~\ref{lem:cyl-distance-entry};
		the minimizing references tend to $\barC$ by Step~1 of the proof of
		Proposition~\ref{prop:joint-full}, so these graphs can also be taken
		over the fitted references.
		
		On an interval $[s,s+2]$ satisfying these distance and energy
		thresholds, put
		$\Delta:=E(s)-E(s+2)$. By \eqref{eq:cyl-fitted-interval} and the comparison with squared distance in
		Lemma~\ref{lem:distance-function}, the fixed minimizing reference satisfies
		\[
		\int_{s+1/2}^{s+3/2}
		\bigl(D_z(N_\tau)^2+G(N_\tau)^2+|q(a)|^2\bigr)\dd\tau
		\le C\Delta+\int_s^{s+2}G(N_\tau)^2\dd\tau
		\le C\Delta.
		\]
		Thus there is a middle time $s_{\rm mid}\in[s+\tfrac12,s+\tfrac32]$ such that
		\begin{equation}\label{eq:cyl-good-middle-time}
			D_z(N_{s_{\rm mid}})^2+G(N_{s_{\rm mid}})^2+|q(a)|^2\le C\Delta.
		\end{equation}
		Apply Proposition~\ref{prop:anchorarea} in coordinates aligned with
		$\barC_z$. Discarding the nonnegative term in \eqref{eq:anchorarea-statement} and using
		\eqref{eq:energy-dissipation} and \eqref{eq:cyl-good-middle-time}, we obtain
		\begin{equation}\label{eq:cyl-energy-area-comparison}
			\begin{aligned}
				E(s)&=E(s_{\rm mid})+\int_s^{s_{\rm mid}} G(N_\tau)^2\dd\tau\\
				&\le F(N_{s_{\rm mid}})-F(\barC_z)+F(\barC_z)-F(\barC)+\Delta\\
				&\le |F(\barC_z)-F(\barC)|+C\Delta.
			\end{aligned}
		\end{equation}
		By \eqref{eq:cyl-energy-area-comparison} and the subadditivity
		of $t\mapsto t^{1-\theta}$ on $[0,\infty)$, we have
		\[
		E(s)^{1-\theta}
		\le C\bigl(|F(\barC_z)-F(\barC)|^{1-\theta}
		+\Delta^{1-\theta}\bigr).
		\]
		By orthogonal invariance and \eqref{eq:cyl-product-functionals},
		\[
		F(\barC_z)-F(\barC)
		=F(\C_a)-F(\C)
		=\frac{f(a)-f(0)}{\sigma_{m-1}},
		\qquad
		\sigma_{m-1}=\mathcal H^{m-1}(S^{m-1}).
		\]
		The finite-dimensional gradient inequality \eqref{eq:reduced-link-gradient}
		therefore gives
		\[
		|F(\barC_z)-F(\barC)|^{1-\theta}
		=\sigma_{m-1}^{-(1-\theta)}|f(a)-f(0)|^{1-\theta}
		\le C|q(a)|.
		\]
		By \eqref{eq:cyl-good-middle-time}, $|q(a)|^2\le C\Delta$.
		Since $1-\theta\ge1/2$, we obtain
		\begin{equation}\label{eq:cyl-local-loja}
			E(s)^{1-\theta}
			\le C\bigl(|q(a)|+\Delta^{1-\theta}\bigr)
			\le C\bigl(\Delta^{1/2}+\Delta^{1-\theta}\bigr)
			\le C\Delta^{1/2}
			\qquad(0\le\Delta\le1).
		\end{equation}
		This proves the first inequality of \eqref{eq:cyl-loja} on every
		late interval satisfying these thresholds.
		
		For the fixed reference, let $U(Z):=\dist(Z,\barC)$.
		It is one-Lipschitz and $U(0)=0$. First use the bounded functions
		$U_K:=\min(U,K)$. Approximate them by bounded nonnegative smooth
		one-Lipschitz functions and apply Gaussian first variation to their
		squares. We drop the negative dissipation term to obtain
		$(\int_{N_s} U_K^2\rho\dd\mu_{N_s})'\le
		2(\int_{N_s} U_K^2\rho\dd\mu_{N_s})^{1/2}G(N_s)$ in integrated form.
		The outer cutoff is removed before passing to the Lipschitz limit,
		using the Gaussian moment bound and finite dissipation.
		We regularize the square root and integrate the resulting inequality.
		We then let $K\to\infty$ and use monotone convergence to obtain
		\begin{equation}\label{eq:cyl-fixed-distance-transport}
			D_0(N_v)\le D_0(N_u)+\int_u^v G(N_\tau)\,\dd\tau
			\qquad(u\le v).
		\end{equation}
		Choose $s_*$ sufficiently far along the prescribed tangent subsequence that
		\begin{equation}\label{eq:cyl-continuation-entry}
			E(s_*)<\epsilon_E,\qquad D_0(N_{s_*})<r_D/4,\qquad
			CE(s_*)^\theta+\sqrt{2E(s_*)}<r_D/4.
		\end{equation}
		Indeed, choose a fixed negative time of slice convergence along the
		parabolic dilations defining the tangent flow. After passing to the
		corresponding logarithmic times, we obtain local convergence of the
		rescaled slices to the unit cylinder by its dilation invariance. Equation
		\eqref{eq:gaussian-moment-bound} then gives $D_0(N_{s_*})\to0$
		along such choices, as well as convergence of the Gaussian areas.
		Energy monotonicity preserves the energy threshold.
		
		For complete two-unit blocks preceding a first exit $T_*$ from
		$D_0<r_D$, write $e_j:=E(s_*+2j)$ and
		$d_j:=\sqrt{e_j-e_{j+1}}$. The local inequality
		\eqref{eq:cyl-local-loja} gives $e_j^{1-\theta}\le C d_j$.
		If $e_j>0$, concavity implies
		\[
		e_j^\theta-e_{j+1}^\theta\ge
		\theta e_j^{\theta-1}d_j^2\ge(\theta/C)d_j .
		\]
		When $e_j=0$ subsequent lengths vanish by dissipation.
		The Cauchy--Schwarz inequality bounds the length of the $j$th block by
		$\sqrt{2}\,d_j$. Summing over $N_{\rm blk}$ such intervals gives
		\[
		\int_{s_*}^{s_*+2N_{\rm blk}}G(N_s)\dd s\le C E(s_*)^\theta
		\qquad(s_*+2N_{\rm blk}<T_*).
		\]
		The remaining interval has length at most two and contributes at most
		$\sqrt{2E(s_*)}$ by the Cauchy--Schwarz inequality and
		\eqref{eq:energy-dissipation}. Hence
		\eqref{eq:cyl-fixed-distance-transport} and
		\eqref{eq:cyl-continuation-entry} imply
		\[
		D_0(N_{T_*})
		\le D_0(N_{s_*})+CE(s_*)^\theta+\sqrt{2E(s_*)}
		<r_D/2,
		\]
		contrary to the definition of $T_*$. Smoothness and properness of the
		flow and the uniform Gaussian tails give continuity of $D_0$.
		Thus both thresholds persist, and \eqref{eq:cyl-local-loja}
		holds at every sufficiently late starting time.
		
		By starting the same block sum at any $s\ge s_*$ and letting the number
		of blocks tend to infinity, we obtain $\int_s^\infty G(N_\tau)\dd\tau\le CE(s)^\theta$.
		For $0<\theta<1/2$ put $d=1-2\theta>0$ and write
		$e_j=E(s_*+2j)$. Then $e_j-e_{j+1}\ge c e_j^{1+d}$.
		As long as $e_{j+1}>0$,
		\[
		e_{j+1}^{-d}-e_j^{-d}
		=d\int_{e_{j+1}}^{e_j}v^{-d-1}\dd v
		\ge d e_j^{-d-1}(e_j-e_{j+1})\ge cd .
		\]
		Summation and monotonicity between grid times yield
		$E(s)\le C(1+s-s_*)^{-1/(1-2\theta)}$. The bound also holds
		after the energy first becomes zero. The length estimate holds for
		$\theta=1/2$. At that endpoint, set $E_j:=E(s_*+2j)$;
		\eqref{eq:cyl-local-loja} gives $E_j-E_{j+1}\ge c E_j$ with $0<c<1$,
		and hence $E_j\le(1-c)^jE_0$. Monotonicity between grid times
		and the length bound give the asserted exponential decay.
		If $E(s_0)=0$, then \eqref{eq:energy-dissipation} implies
		$G(N_s)=0$ for every $s\ge s_0$.
		
		For $\psi\in C_c^1(\R^{n+1})$ and $v>u\ge s_*$, Gaussian first
		variation gives
		\[
		\left|\int_{N_v}\psi\rho\dd\mu_{N_v}
		-\int_{N_u}\psi\rho\dd\mu_{N_u}\right|
		\le C_\psi\int_u^v G(N_s)\dd s+
		\|\psi\|_\infty\bigl(E(u)-E(v)\bigr).
		\]
		The right side tends to zero as $u\to\infty$ by finite Gaussian
		length. The prescribed tangent subsequence identifies the common
		limit with $\int_{\barC}\psi\rho\dd\mu_{\barC}$. Since $\rho$ is
		positive and the local masses are bounded, approximation gives
		local weak convergence $\mu_{N_s}\rightharpoonup\mu_{\barC}$ in
		the original coordinates. For arbitrary $\lambda_i\downarrow0$, the ordinary
		dilated flows are $\sqrt{-\tau}\,N_{-2\log\lambda_i-\log(-\tau)}$.
		By homogeneity and the uniform area bound, the area measures
		$\mu_{\sqrt{-\tau}\,N_{-2\log\lambda_i-\log(-\tau)}}$ converge locally
		weakly as Radon measures to $\mu_{\barC}$, uniformly for $\tau$ in
		compact negative time intervals. Local regularity
		\cite[Section~6.11]{Brakke} gives smooth convergence away from the axis.
		
		On every bounded axial
		region, \eqref{eq:raw-tube-area} gives
		\[
		\mu_{\sqrt{-\tau}\,N_{-2\log\lambda_i-\log(-\tau)}}
		\{r<\delta,\ |y|<L\}\le C_L\delta^m,
		\]
		uniformly for $\tau$ in compact negative time intervals, by the
		scale-invariant area-ratio bound. The same estimate holds for
		$|\barC|$. For any compactly supported varifold test, the contributions
		from this tube are therefore $O(\delta^m)$, whereas the contributions
		away from the axis converge by smooth convergence. First let
		$i\to\infty$ and then $\delta\downarrow0$. This proves varifold
		convergence across the axis and identifies every tangent flow with
		the static cylinder $|\barC|$.
	\end{proof}
	
	\subsection{Rigidity among smooth self-shrinkers}\label{sec:shrinker-rigidity}

	\begin{proof}[Proof of Corollary~\ref{cor:shrinker-rigidity}]
		Suppose that smooth normalized shrinkers $S_i$ converge locally as
		varifolds with multiplicity one to $|\barC|$. The Gaussian tail and
		area estimates in the proof of \cite[Corollary~1.4]{FH1} apply in
		dimension $n=m+k$, since $\dist(Z,\barC)\le |Z|$. Hence we get
		\[
		F(S_i)\to\Theta,\qquad D_0(S_i)\to0,\qquad
		\mathcal H^n(S_i\cap B_R(x))\le C R^n
		\quad(x\in\R^{n+1},\ R>0),
		\]
		with $C$ independent of $i,x,R$.
		
		Regard $S_i$ as a constant rescaled flow on an interval of length two.
		For large $i$, \eqref{eq:cyl-fitted-smallness} holds, and $G(S_i)=0$.
		Thus Proposition~\ref{prop:joint-full} gives a parameter
		$z_i=(a_i,\mathsf B_i)$ such that
		\[
		\int_{S_i}W_{z_i}\rho\dd\mu_{S_i}+|q(a_i)|^2=0.
		\]
		Since $W_{z_i}$ is comparable to squared distance, $S_i\subset\barC_{z_i}$.
		The regular part of $\barC_{z_i}=R_{\mathsf B_i}(\C_{a_i}\times\R^k)$
		is connected, and its axis has dimension $k<n$. Thus $S_i$ meets the
		regular part, and the open and closed argument in the proof of
		\cite[Corollary~1.4]{FH1} gives $S_i=\barC_{z_i}$.
		The nearby transverse links remain nonflat, so this cylinder is singular
		along its axis, contrary to the smoothness of $S_i$.
	\end{proof}
	
	\section{Examples}\label{sec:examples}
	In this section, we verify the hypotheses for quadratic cone cylinders in transverse
	dimension $m\ge8$. We first compute the link spectrum and identify
	its Jacobi fields. Then we prove weighted stationary rigidity by
	comparison with transverse barriers whose scale varies in the axial
	directions. 
	
	\begin{prop}\label{prop:quadratic-cylinder-example}
		For $p,q\ge1$ and $m:=p+q+1\ge8$, the cone $\C_{p,q}$ satisfies
		strict Hardy stability, $\ker\mathcal A_\Sigma=T_\C$, and
		cylindrical Gaussian nonresonance. Its reduced link area is locally
		constant.
	\end{prop}
	\begin{proof}
		The link is $S^p(\sqrt{p/(m-1)})\times S^q(\sqrt{q/(m-1)})$.
		With either choice of normal, the two principal curvatures
		have absolute values $\sqrt{q/p}$ and $\sqrt{p/q}$ with
		multiplicities $p$ and $q$. Thus $|A_\Sigma|^2=m-1$.
		Tensor products of spherical harmonics of degrees $i,j$ give a complete
		basis, with link eigenvalues
		\[
		\mu_{i,j}=(m-1)\left(\frac{i(i+p-1)}p+
		\frac{j(j+q-1)}q-1\right).
		\]
		The constant mode gives $\mu_{0,0}=-(m-1)$. Thus
		$\nu_0^2=(m-2)^2/4-(m-1)>0$ for $m\ge8$, proving
		\eqref{eq:cyl-S1}. All remaining
		eigenvalues are nonnegative.
		To locate the larger root when
		$m\ge8$, put $p_m(\gamma):=\gamma(\gamma+m-2)+(m-1)$. Then
		\[
		p_m(-2)=7-m<0,\qquad p_m(-1)=2>0,
		\qquad p_m'(\gamma)=2\gamma+m-2\ge m-6>0
		\quad(-2\le\gamma\le-1).
		\]
		The larger principal root therefore lies in $(-2,-1)$.
		For every other angular mode, $\mu_{i,j}\ge0$, so
		$\gamma_{i,j}^+\ge0$. Hence no larger indicial root is a negative
		integer, which proves \eqref{eq:cyl-S3}.
		
		The zero eigenspace has precisely bidegrees
		$(1,0),(0,1)$, which are the components of the normal
		$(\sqrt{q/p}\,x_1,-\sqrt{p/q}\,x_2)$. Thus $\ker\mathcal A_\Sigma=T_\C$.
		
		The eigenvalue $m-1$ has precisely bidegree $(1,1)$, because
		$j(j+d-1)/d>2$ for every $j\ge2$. These fields are exactly
		the normal velocities of transverse rotations mixing the two
		sphere factors. To see this, their velocities are nonzero
		constant multiples of $\langle\mathsf E x_2,x_1\rangle$ for arbitrary
		linear maps $\mathsf E$ between the factors. The derivative of the
		map taking such a rotated link to the projection of its
		spherical graph height onto $K$ is therefore invertible.
		The inverse function theorem parametrizes the rotated links
		by every sufficiently small $a\in K$. Since each is minimal,
		uniqueness in the complementary equation defining
		\eqref{eq:reduced-link-data} identifies it with
		$\Sigma_a$. This proves that $f(a)=f(0)$ there.
		
	\end{proof}

	\begin{prop}\label{prop:quadratic-static}
		Let $p,q\ge1$, $m:=p+q+1\ge8$, and $k\ge1$.
		There is $\epsilon_{p,q}>0$ such that
		$\barC:=\C_{p,q}\times\R^k$ satisfies \eqref{eq:cyl-S4} whenever
		\[
		0<B-k/2<\epsilon_{p,q}.
		\]
	\end{prop}
	\begin{proof}
		We fix $p,q\ge1$, $m:=p+q+1\ge8$, and let
		$\C:=\C_{p,q}=\{q|x_1|^2=p|x_2|^2\}$ and $\barC:=\C\times\R^k$,
		where $x=(x_1,x_2)$. We prove rigidity for the finite moment
		\begin{equation}\label{eq:quad-cyl-moment}
			\int_{\R^{n+1}}\frac{\operatorname{dist}(Z,\barC)^2}
			{(1+r^2)(1+|y|^2)^B}\dd\|V\|(Z)<\infty,\qquad
			r:=|x|,\quad B>k/2,
		\end{equation}
		provided $B-k/2$ is sufficiently small, depending on $p,q$.
		The varifold $V$ is stationary integral, has finite all-center
		area ratios, and has full multiplicity-one blowdown $|\barC|$.
		
		\noindent\textbf{Step 1. Decay of the support in bounded axial sectors.}
		There is a smooth homogeneous degree-two transverse function
		$Q_0$ on $\R^{m+1}\setminus\{0\}$, vanishing transversally
		precisely on $\C\setminus\{0\}$, such that
		$|Q_0|\asymp r\operatorname{dist}(x,\C)$.
		We take the degree-two extension of a smooth signed defining function
		of the link on the sphere and set $Q:=Q_0^2$.
		For any $n$-plane $S$ in the ambient space,
		\begin{equation}\label{eq:cyl-static-trace-bound}
			\operatorname{tr}_S D^2Q\ge -C r^{-2}Q.
		\end{equation}
		To see this, we use homogeneity to work at $r=1$. Near the link we set $\mathbf n_0:=DQ_0/|DQ_0|$, and put $\tau:=|\operatorname{proj}_S \mathbf n_0|$.
		Minimality on $Q_0=0$ and smoothness give
		\[
		|\operatorname{tr}_S D^2Q_0|
		\le |\operatorname{tr}_{\mathbf n_0^\perp}D^2Q_0|+C\tau
		\le C|Q_0|+C\tau,\qquad |DQ_0|\ge c>0.
		\]
		Consequently,
		\[
		\begin{aligned}
			\operatorname{tr}_S D^2Q
			&=2Q_0\operatorname{tr}_S D^2Q_0+2|\operatorname{proj}_S DQ_0|^2\\
			&\ge -C Q_0^2-C|Q_0|\tau+2c^2\tau^2\\
			&\ge -C Q_0^2-c^2\tau^2+2c^2\tau^2\ge-CQ.
		\end{aligned}
		\]
		Away from the link $|Q_0|$ is bounded below, so \eqref{eq:cyl-static-trace-bound} at $r=1$ follows
		from compactness and boundedness of $D^2Q$ on $r=1$. Since $Q$ is
		homogeneous of degree four and $D^2Q$ of degree two, dilation gives
		\eqref{eq:cyl-static-trace-bound} at every $r>0$.
		
		We write $\omega_n:=\operatorname{Vol}_n(B_1^n)$ for the volume of the unit ball in $\mathbb R^n$.
		For a nonnegative smooth function $w$ on a ball $B_R(Z_0)$
		satisfying $\operatorname{tr}_S D^2 w\ge-\kappa w$, stationarity
		gives the weighted mean-value bound
		\begin{equation}\label{eq:cyl-static-mean-value}
			w(Z_0)\le \omega_n^{-1} e^{\kappa R^2/4}R^{-n}
			\int_{B_R(Z_0)}w\dd\|V\|
			\qquad(Z_0\in\operatorname{spt}\|V\|).
		\end{equation}
		Indeed, writing $A_w(t):=\int_{B_t(Z_0)}w\dd\|V\|$ and
		$\varrho:=|Z-Z_0|$, first variation with $w\zeta(\varrho)(Z-Z_0)$
		gives, in measures in $t$,
		\[
		t\dd A_w(t)-n A_w(t)\dd t\ge
		\left(\int_{B_t(Z_0)}(Z-Z_0)\cdot\operatorname{proj}_S Dw
		\dd\|V\|\right)\dd t .
		\]
		Here $|\operatorname{proj}_S D\varrho|^2\le1$ in the term
		from the decreasing cutoff. The integral on the right is,
		by testing with $\frac12(t^2-\varrho^2)_+ Dw$, equal to
		$\int_{B_t(Z_0)}(t^2-\varrho^2)\operatorname{tr}_S D^2w/2\dd\|V\|$
		for almost every $t$. It is at least $-\kappa t^2 A_w(t)/2$.
		Thus $e^{\kappa t^2/4}t^{-n}A_w(t)$ is nondecreasing.
		As $t\downarrow0$, its limit is the product of $w(Z_0)$, the density at $Z_0$, and $\omega_n$.
		The density is at least one at each support point, by stationary
		monotonicity \cite{SimonGMT} at nearby points of positive integer density and
		then a limit of the centers. This proves \eqref{eq:cyl-static-mean-value}.
		
		Now we apply \eqref{eq:cyl-static-mean-value} with $w=Q$, $R=r/4$,
		and $\kappa=Cr^{-2}$, where $Z=(x,y)$ is the support point and
		$r:=|x|\ge1$. The transverse radii $r'$ on this ball satisfy
		$3r/4<r'<5r/4$, and $|y'|\le |y|+r/4$. Put
		\[
		\mathcal T(r):=\int_{\{Z'\in\R^{n+1}:r'>3r/4\}}
		\frac{\dist(Z',\barC)^2}{(1+(r')^2)(1+|y'|^2)^B}\dd\|V\|(Z').
		\]
		By \eqref{eq:quad-cyl-moment}, $\mathcal T(r)\to0$ as $r\to\infty$.
		Applying the mean-value estimate and using $Q\asymp r^2\dist(\cdot,\barC)^2$, we obtain
		\begin{equation}
			\begin{aligned}
				r^2\dist(Z,\barC)^2
				&\le C r^{-n}\int_{B_{r/4}(Z)}Q(Z')\dd\|V\|(Z')\\
				&\le C r^{2-n}\int_{B_{r/4}(Z)}\dist(Z',\barC)^2\dd\|V\|(Z')\\
				&\le C r^{4-n}(r+|y|)^{2B}\mathcal T(r).
			\end{aligned}
		\end{equation}
		By taking square roots, we obtain the following bound uniformly along support points,
		\[
		\dist(Z,\barC)\le C r^{1-n/2}(r+|y|)^B\mathcal T(r)^{1/2}
		=o\bigl(r^{1-n/2}(r+|y|)^B\bigr).
		\]
		In particular, with $\varepsilon:=B-k/2$ and any fixed $L$,
		\begin{equation}\label{eq:cyl-static-sector}
			\operatorname{dist}(Z,\barC)
			=o(r^{1-m/2+\varepsilon})
			\quad\text{as }r\to\infty,\quad |y|\le Lr .
		\end{equation}
		
		\medskip
		\noindent\textbf{Step 2. A local strict barrier test.}
		Suppose a smooth function $h$ has nonzero gradient at a
		point of the stationary support, reaches a local maximum there,
		and has strictly positive Hessian trace on its level tangent
		plane there. This is impossible. Indeed, we subtract the maximum
		value and a small multiple of the squared distance from the
		point. The strict trace persists nearby, and the modified
		function is strictly negative on the support away from the
		contact point. On a sufficiently small ball we set
		$\mathfrak v:=1/|Dh|$ and $\mathbf n_h:=Dh/|Dh|$. If the normal to an $n$-plane is
		$c\mathbf n_h+\tau$, where $\tau\perp\mathbf n_h$, direct differentiation
		gives
		\[
		\operatorname{tr}_S D(\psi(h)\mathfrak v\mathbf n_h)
		=\psi(h)\mathfrak v\operatorname{div}\mathbf n_h
		+\psi'(h)|\tau|^2+
		\psi(h)\bigl(D_{\mathbf n_h}\mathfrak v\,|\tau|^2
		-\mathfrak v\langle D_\tau\mathbf n_h,\tau\rangle\bigr).
		\]
		The mixed terms cancel because
		$D_T \mathfrak v+\mathfrak v\langle D_{\mathbf n_h}\mathbf n_h,T\rangle=0$
		for $T\perp\mathbf n_h$, by symmetry of $D^2h$.
		The first term is strictly positive wherever $\psi>0$.
		We choose $\psi$ to vanish below a small negative threshold and to satisfy
		$\psi'>(C+1)\psi$ above it, so that the second term controls
		the last. A compact ambient cutoff equal to one on the support
		superlevel causes no error, since the contact is isolated.
		The strictly positive trace on positive support mass contradicts
		stationarity. By applying the same argument to $-h$, we obtain the corresponding statement
		for a local minimum with strictly negative trace.
		
		\medskip
		\noindent\textbf{Step 3. Construction of transverse barriers.}
		We use the angular coordinates
		\[
		|x_1|=r\sin\vartheta,\qquad |x_2|=r\cos\vartheta,\qquad
		0\le\vartheta<\vartheta_*:=\arctan\sqrt{p/q},\qquad
		H:=p\cot\vartheta-q\tan\vartheta
		\]
		on the side adjacent to the $x_2$-axis. The other side is obtained
		by interchanging the factors. For a radial profile
		$\mathcal R'/\mathcal R=\cot\alpha$, we set $\lambda(x):=r/\mathcal R(\vartheta)$.
		The transverse unit normal in the increasing $\lambda$ direction
		is $\mathbf n_\lambda=\sin\alpha\,e_r-\cos\alpha\,e_\vartheta$.
		The polar volume density is proportional to
		$r^m\sin^p\vartheta\cos^q\vartheta$, and hence
		\begin{equation}\label{eq:quad-cyl-div}
			\operatorname{div}\mathbf n_\lambda=
			\frac{\sin\alpha}{r}(m+\alpha'-H\cot\alpha).
		\end{equation}
		We set $\alpha:=\arctan(H/d)$ with $d>0$ constant and $\mathcal R(0)=1$. We also set $s:=p+q$ and $\beta:=d/(2s)$. The profile is explicit:
		\begin{equation}\label{eq:cyl-explicit-barrier-profile}
			\begin{aligned}
				H(\vartheta)&=\frac{p-s\sin^2\vartheta}{\sin\vartheta\cos\vartheta},\\
				\mathcal R(\vartheta)
				&=\exp\left(\int_0^\vartheta\frac{d}{H(\tau)}\dd\tau\right)
				=\left(\frac{p}{p-s\sin^2\vartheta}\right)^\beta.
			\end{aligned}
		\end{equation}
		In particular, $\mathcal R$ is even analytic at the focal axis,
		and the phase is smooth there for $r>0$.
		
		We choose $d$ so that the bracket in \eqref{eq:quad-cyl-div}
		is positive, with a uniform positive lower bound up to both
		endpoints, and so that $d>4(p+q)/m$. Its sign is that of
		$\mathscr N:=(m-d)(d^2+H^2)+dH'$. For $d=m/2$, the sign reduces to
		\begin{equation}\label{eq:quad-phase-sign}
			\frac{m^2}{4}-2pq-(p+q)
			+p(p-1)\cot^2\vartheta+q(q-1)\tan^2\vartheta .
		\end{equation}
		For $p,q\ge2$, we set $v:=pq\ge2(s-2)$.
		The arithmetic--geometric mean inequality bounds \eqref{eq:quad-phase-sign} from below by
		\[
		(s-1)\left(\frac{s-1}{4}
		-\frac{2}{1+\sqrt{1-(s-1)/v}}\right)>0 ,
		\]
		since $s\ge6$ and $\sqrt{1-(s-1)/v}\ge\sqrt{3/8}>3/5$.
		For $q=1$, $p\ge6$, we use $\cot^2\vartheta\ge1/p$ to obtain
		the positive lower bound $(p^2-4p-4)/4$.
		At the focal axis the bracket tends to $m-d-d/p>0$ in both
		these cases. At $\vartheta_*$ it tends to
		$m-d-2(p+q)/d>0$, since $m^2>8(p+q)$.
		
		For $p=1$, $q\ge8$, the expression at $d=m/2$ is
		$q(q-8)/4+q(q-1)\tan^2\vartheta$. We choose $d$ slightly less
		than $m/2$ to obtain a positive focal limit. In detail, the numerator $\mathscr N$ is
		$\mathfrak a_0\cot^2\vartheta+\mathfrak b_0 \tan^2\vartheta+\mathfrak c_0$ with $\mathfrak a_0:=m-2d$, $\mathfrak b_0:=q^2(m-d)-qd$,
		and $\mathfrak c_0:=(m-d)(d^2-2q)-d(q+1)$. At $d=m/2$, $\mathfrak b_0>0$, $\mathfrak c_0\ge0$.
		For $d=m/2-\delta$, the lower bound $2\sqrt{\mathfrak a_0\mathfrak b_0}+\mathfrak c_0$
		is positive for small $\delta>0$, since $\mathfrak a_0=2\delta$ and
		$\mathfrak c_0\ge-C\delta$. The strict inequality $d>4(p+q)/m$ persists.
		For $p=1$, $q=6,7$, we first choose $d=7/2,4$, respectively.
		With $\zeta:=q\tan^2\vartheta$ the polynomials $\zeta \mathscr N$ are
		$6-187\zeta/8+47\zeta^2/2$ and $7-22\zeta+31\zeta^2$, with negative
		discriminants $-1127/64$ and $-384$, respectively, and positive leading coefficients.
		Their focal and outer bracket limits are positive.
		For $q=6$ we increase $d$ slightly so that $d>4(p+q)/m$,
		preserving the strict signs. This covers all sides when $m\ge8$.
		
		By \eqref{eq:cyl-explicit-barrier-profile}, as $\vartheta\uparrow\vartheta_*$,
		$p-s\sin^2\vartheta\asymp\vartheta_*-\vartheta$ and
		$\dist(Z,\barC)=r\sin(\vartheta_*-\vartheta)$. Thus
		\begin{equation}\label{eq:cyl-barrier-phase-decay}
			\lambda=\frac r{\mathcal R(\vartheta)}
			\asymp r(\vartheta_*-\vartheta)^\beta
			\asymp r^{1-\beta}\dist(Z,\barC)^\beta
			=o\bigl(r^{1-\beta(m/2-\varepsilon)}\bigr)
		\end{equation}
		along support points with $r\to\infty$ and $|y|\le Lr$,
		where the last equality uses \eqref{eq:cyl-static-sector}.
		Our choice $d>4s/m$ gives $\beta>2/m$. We choose the positive
		$\varepsilon$ in \eqref{eq:cyl-static-sector} smaller than
		$m/2-1/\beta$ for both sides. Then
		$1-\beta(m/2-\varepsilon)<0$, so \eqref{eq:cyl-barrier-phase-decay}
		gives $\lambda\to0$ in each bounded axial sector.
		Also $\lambda\le r$ and $\lambda$
		extends continuously by zero to the cone and the axis.
		
		\medskip
		\noindent\textbf{Step 4. Axially varying barriers and the maximum argument.}
		Suppose that the support meets this open side of the cone. We set
		$L_{\mathfrak s}(y):=(\mathfrak s^2+|y|^2)^{1/2}$. For each $\mathfrak s>0$ the phase quotient
		$\lambda/L_{\mathfrak s}$ tends to zero along support points tending to
		infinity in the side. If $r/|y|\to0$, the decay follows from
		$0\le\lambda/L_{\mathfrak s}\le r/|y|\to0$.
		Otherwise we pass to a subsequence with $|y|\le Lr$ and use the
		phase decay \eqref{eq:cyl-barrier-phase-decay}. These alternatives prove $\lambda/L_{\mathfrak s}\to0$
		along every sequence of support points tending to infinity in the side.
		Consequently the quotient has a positive maximum $t_{\mathfrak s}$
		at a point $Z_{\mathfrak s}$ of the support in the open side.
		At any fixed support point in the side,
		$\mathfrak s\lambda/L_{\mathfrak s}$ converges to its positive phase. Thus
		$t_{\mathfrak s} \mathfrak s\ge c_0>0$ for large $\mathfrak s$. In particular, the phases at
		the maximizing points are bounded below. We have
		$r_{\mathfrak s}/L_{\mathfrak s}(y_{\mathfrak s})\to0$: otherwise a subsequence would have
		$r_{\mathfrak s}\ge c \mathfrak s$ and $|y_{\mathfrak s}|\le c^{-1}r_{\mathfrak s}$, contradicting the
		sector estimate for the phase.
		
		At the maximum, we use the level function
		$h:=\lambda(x)-t_{\mathfrak s} L_{\mathfrak s}(y)$. We set $w:=|D_x\lambda|$,
		$\mathbf v_\lambda:=t_{\mathfrak s}D_yL_{\mathfrak s}$, and $\mathbf n_\lambda:=D_x\lambda/w$.
		The normal to its level is $(w\mathbf n_\lambda,-\mathbf v_\lambda)/(w^2+|\mathbf v_\lambda|^2)^{1/2}$.
		The Hessian trace on that tangent plane is bounded from below by
		\[
		w\operatorname{div}\mathbf n_\lambda-t_{\mathfrak s}\Delta L_{\mathfrak s}
		-\frac{|\mathbf v_\lambda|^2}{w^2}\,|D^2_x\lambda(\mathbf n_\lambda,\mathbf n_\lambda)| ,
		\]
		where the omitted axial normal correction is nonnegative by
		convexity of $L_{\mathfrak s}$. Uniformly in the angular variable we have
		\begin{equation}\label{eq:cyl-barrier-phase-derivatives}
			w\operatorname{div}\mathbf n_\lambda\ge\frac{c}{r\mathcal R},\qquad
			\frac{|D^2_x\lambda|}{w^2}\le C \mathcal R/r .
		\end{equation}
		The first follows from \eqref{eq:quad-cyl-div} and
		$w\sin\alpha=\mathcal R^{-1}$. For the second inequality we set $\mathcal S:=\mathcal R^{-1}$. Near
		$\vartheta_*$ we have $|\mathcal S'|\asymp \mathcal S/(\vartheta_*-\vartheta)$
		and $|\mathcal S''|\le C\mathcal S/(\vartheta_*-\vartheta)^2$.
		Since $\lambda=r\mathcal S$, its transverse Hessian is
		$r^{-1}(\nabla^2_{S^m}\mathcal S+\mathcal S g_{S^m})$ in spherical directions,
		and $w^2=\mathcal S^2+|\mathcal S'|^2$. These identities give the second inequality in \eqref{eq:cyl-barrier-phase-derivatives} near $\vartheta_*$. Away from this endpoint, smoothness at the focal axis,
		positivity of $\mathcal S$, and $w\ge \mathcal S$ prove the second inequality in \eqref{eq:cyl-barrier-phase-derivatives} on the remaining angular region.
		Now $|D L_{\mathfrak s}|\le1$, $\Delta L_{\mathfrak s}\le k/L_{\mathfrak s}$, and
		$r_{\mathfrak s}=t_{\mathfrak s} L_{\mathfrak s} \mathcal R$ at contact. At the contact point, we write $r:=r_{\mathfrak s}$, $L:=L_{\mathfrak s}(y_{\mathfrak s})$, and use
		$t_{\mathfrak s}=r/(L\mathcal R)$. The two error terms satisfy
		\begin{equation}
			\begin{aligned}
				t_{\mathfrak s}\Delta L_{\mathfrak s}&\le\frac{k t_{\mathfrak s}}{L}
				=\frac{kr}{L^2\mathcal R},\\
				\frac{|\mathbf v_\lambda|^2}{w^2}
				|D_x^2\lambda(\mathbf n_\lambda,\mathbf n_\lambda)|
				&\le Ct_{\mathfrak s}^2\frac{\mathcal R}{r}
				=\frac{Cr}{L^2\mathcal R}.
			\end{aligned}
		\end{equation}
		Consequently the trace on the level tangent plane obeys
		\[
		\operatorname{tr}_{(Dh)^\perp}D^2h
		\ge\frac{c}{r\mathcal R}-\frac{Cr}{L^2\mathcal R}
		=\frac1{r\mathcal R}\left(c-C\frac{r^2}{L^2}\right)>0
		\]
		for all sufficiently large $\mathfrak s$, because $r_{\mathfrak s}/L_{\mathfrak s}(y_{\mathfrak s})\to0$.
		Furthermore $D_xh=D_x\lambda\ne0$, and maximality of
		$\lambda/L_{\mathfrak s}$ says $h\le0$ on the support with equality at $Z_{\mathfrak s}$.
		The strict barrier test from Step 2 therefore gives a contradiction. We repeat
		the argument on the other side and conclude that the support lies in $\barC$.
		The axis has dimension $k<n$ and carries no mass by
		rectifiability. The constancy theorem \cite{SimonGMT} gives a constant
		integer multiplicity on the connected regular part of $\barC$,
		and the prescribed blowdown fixes it to be one.
		
	\end{proof}
	
	\begin{proof}[Proof of Corollary~\ref{cor:examples}]
		Propositions~\ref{prop:quadratic-cylinder-example} and
		\ref{prop:quadratic-static} verify the cone hypotheses of
		Theorem~\ref{thm:main} for sufficiently small $B-k/2>0$. Hence the tangent flow is unique.
		
		By Proposition~\ref{prop:quadratic-cylinder-example}, $f(a)=f(0)$ throughout the reduced family. We may therefore take
		$\theta=1/2$ in \eqref{eq:reduced-link-gradient}. With this choice,
		$2(1-\theta)=1$, so \eqref{eq:cyl-loja} gives
		\[
		E(s)-E(s+2)\ge cE(s)\qquad(s\ge s_0)
		\]
		for some $0<c<1$. Iterating,
		we obtain $E(s_0+2j)\le(1-c)^jE(s_0)$. Monotonicity of $E$ between
		these times gives $E(s)\le C E(s_0)e^{-\lambda(s-s_0)}$ for all
		$s\ge s_0$, with $\lambda=-\tfrac12\log(1-c)>0$.
	\end{proof}

	Hanbing Fang, Department of Mathematics, The Chinese University of Hong
	Kong, Shatin, N.T., Hong Kong; E-mail: hanbingfang@cuhk.edu.hk\\

	Yongheng Han, Institute of Mathematics, Academy of Mathematics and Systems Science,
	Chinese Academy of Sciences, Beijing 100190, China; E-mail: hyh98@amss.ac.cn
	
\end{document}